\documentclass[11pt]{article}
\usepackage{graphicx}
\usepackage{amssymb,amsfonts,amsthm}
\usepackage{enumerate}
\usepackage[T2B]{fontenc}
\usepackage[cp1251]{inputenc}

\usepackage{enumitem}
\usepackage{amsmath,amsthm,amssymb}
\usepackage{amscd}
\usepackage{amsfonts}

\usepackage{xcolor}

\newcommand{\la}{\langle}
\newcommand{\ra}{\rangle}

\newtheorem{theorem}{Theorem}[section]
\newtheorem{lemma}[theorem]{Lemma}

\theoremstyle{definition}
\newtheorem{df}[theorem]{Definition}

\newtheorem{rk}[theorem]{Remark}

\newcommand{\area}{\mathrm{Area}}

\newcommand{\me}{\medskip}

\newcommand{\Lab}{{\mathrm{Lab}}}

\newcommand{\tool}{\stackrel{\ell}{\too} }

\newcommand{\rhh}{\mathbf{LR}}
\newcommand{\rhr}{\mathbf{RL}}
\newcommand{\mmm}{\mathbf{M}}
\newcommand{\mtt}{\mathbf{\Theta}}
\renewcommand{\lll}{{\mathcal L}}

\newcommand{\Z}{\mathbb Z}

\newcommand{\bmm}{\mathbf{\overline{\mmm}}}

\newcommand{\ttt}{{\mathcal T}}
\renewcommand{\tt}{{\tilde t}}

\newcommand{\bb}{{\mathcal B}}
\newcommand{\topp}{\mathbf{top}}

\newcommand{\bott}{\mathbf{bot}}
\newcommand{\vk}{van Kampen }

\newcommand{\ccc}{{\mathcal C}}

\newcommand{\iv}{^{-1}}

\newcommand{\too}{\to }

\newcommand{\rrr}{{\mathcal R} }

\newcommand{\sss}{{\mathbf S} }

\usepackage{imakeidx}
\makeindex[name=g,title=Subject index]
\usepackage[columns=2]{idxlayout}
\usepackage[pagebackref, breaklinks]{hyperref}

    \renewcommand*{\backref}[1]{}
    \renewcommand*{\backrefalt}[4]{%
    \ifcase #1 %
        (Not cited).
    \or
        (Cited on page~#2).%
    \else
        (Cited on pages~#2).
    \fi}

 \usepackage{fancyhdr,lipsum}
\begin{document}
\renewcommand{\theequation}{\thesection.\arabic{equation}}
\bigskip

\title{Algorithmic problems in groups with quadratic Dehn function}

 \author{A.Yu. Olshanskii, M.V. Sapir\thanks{Both authors were supported in part by the NSF grant DMS-1901976.}}

\date{}
\maketitle

\renewcommand\rightmark{[Short Title]}

\begin{abstract}
We construct and study finitely presented groups with quadratic Dehn function (QD-groups) and present the following applications of the method developed in our recent papers.
(1) The isomorphism problem is undecidable in the class of QD-groups.
(2) For every recursive function $f$, there is a QD-group $G$ containing
a finitely presented subgroup $H$ whose Dehn function grows faster than $f$.
(3) There exists a group with undecidable
conjugacy problem but decidable power conjugacy problem; this group is QD.
\end{abstract}

{\bf Key words:} generators and relations in groups, finitely presented groups, the Dehn function of a group, algorithm,
$S$-machine,  van Kampen
diagram.

\medskip

{\bf AMS Mathematical Subject Classification:} 20F05, 20F06,  20F65,  03D10.

\setcounter{tocdepth}{2}

\tableofcontents

\section{Introduction}

Every group given by a presentation $G=\la X \mid R \ra$ is a factor group $F/N$ of the free group $F=F(X)$
with the set of free generators $X$ over the normal closure $N=\la\la R\ra\ra^F$ of the set of relators $R$. Therefore every word $w$ over
the alphabet $X^{\pm 1}$ vanishing in $G$ represents
an element of $N$, and so in $F$, $w$ is a product $v_1\dots v_m$
of factors $v_i=u_ir_i^{\pm 1}u_i^{-1}$ which are  conjugate of the relators $r_i\in R$ or their inverses.

The minimal number of factors $m=m(w)$ is called the {\it area of the word $w$} with respect to the presentation
$G=\la X | R \ra $. M. Gromov \cite{GrHyp, Gr} introduced this concept in geometric group theory, because $m$ is equal to the
minimal number of $2$-cells (counting with multiplicities) used in a $0$-homotopy of the path $\bf p$ labeled by $w$ in the
Cayley complex of the presentation of $G$ (or the $0$-homotopy of a
singular disk with boundary $\bf p$).

In other words, given equality $w=1$ in $G$, one can construct a van Kampen diagram, that is a finite, connected graph
on Euclidean plane with $m$ bounded regions, where every edge has label from $X^{\pm 1}$, the boundary path of
every region (= $2$-cell) is therefore labelled,  the label of it belongs in $R^{\pm 1}$, and the boundary of  the whole
map is labelled by $w$. (See more details for this visual definition of area and van Kampen diagram in Section \ref{md}.)

The \index[g]{Dehn function of a finitely presented group}
Dehn function of a finitely presented group $G=\la X\mid R\ra$ is the smallest function $f(n)$ such that for every word $w$ of length at most $n$ in the alphabet $X^{\pm 1}$, which is equal to $1$ in $G$, the area of $w$ is at most $f(n)$.

It is well known (see \cite{MO}) that the Dehn functions of different finite presentations of the same finitely presented group are equivalent, where we call two functions $f, g$ equivalent if $f\preceq g$ and $g\preceq f$. Here $f\preceq g$ means that there is a constant $c>0$ such that $f(n)\le cg(cn)+ cn$ for every $n=1,2,\dots$.

The Dehn function is an important invariant of a group for the following reasons.

\begin{itemize}
\item  It easily follows from the definition that if $G$ is the fundamental group of a compact Riemannian manifold $M$ then the Dehn function of $G$ is equivalent to the smallest isoperimetric function of the universal cover $\tilde M$.

\item From the Computer Science point of view, the Dehn function of a group $G$ is equivalent to the time function of a non-deterministic Turing machine "solving" the word problem in $G$ (see \cite[Introduction]{SBR} for details). Moreover as was  shown in \cite{BORS}:

 \begin{quote}
 A not necessarily finitely presented finitely generated group has word problem in $\mathbf{NP} $ if and only if it is a subgroup of a finitely presented group with at most polynomial Dehn function (a similar result holds for other computational complexity classes \cite{BORS}).
\end{quote}

\item From the geometric point of view the Dehn function measures the "curvature" of the group: linear Dehn functions correspond to negative curvature, quadratic Dehn function correspond to non-positive curvature, etc.
\end{itemize}

More precisely, a finitely presented group is hyperbolic if and only if it has a subquadratic (hence linear) Dehn function \cite{GrHyp, Bow, Ol91}. In particular, the conjugacy problem in such groups is decidable \cite{GrHyp}.
In contrast, we recently constructed a group with quadratic Dehn function and undecidable conjugacy problem \cite{OS20}. This result answers Rips' question of 1994. The present paper is based on the constructions of groups with small
Denn functions from \cite{O18}, \cite{OS20} as well as on the application of S-machines introduced in \cite{SBR}.

The affirmative solution of the isomorphism problem was obtained in \cite{Sel} for the class of torsion free hyperbolic groups and in \cite{DG} for the class of all hyperbolic groups. This means that there exists an algorithm recognizing whether two hyperbolic groups $G_1$ and $G_2$ are isomorphic or not, provided $G_1$ and $G_2$ are given by their finite presentations. We show that the linearity is the only
possible restriction of Dehn functions providing a positive solution of the
isomorphism problem.

\begin{theorem}\label{isom} In the class $QD$ of finitely presented groups
with quadratic Dehn function, the isomorphism problem is undecidable.
Moreover, one can select a $QD$-group $\bar G$ such that there exist
no algorithm deciding if a $QD$-group $G$ and $\bar G$ are isomorphic or not.
 \end{theorem}

It is known that the Dehn function of a finitely presented subgroup
can grow faster than the Dehn function of the entire group. For
example, the group $SL(5, Z)$ has quadratic Dehn function \cite{Y},
but it contains subgroups with exponential Dehn function. Here we prove the following

\begin{theorem}\label{subgr} For every recursive function $f$, there exists a pair of finitely presented groups $H\le G$, such that $f\preceq d_H$, where $d_H$ is the Dehn function of the subgroup $H$, while the Dehn function of $G$ is quadratic.
\end{theorem}

For a group with presentation $G=\la X \mid R \ra$, the power conjugacy problem is to determine, given words
$u,v\in F(X)$ whether or not there exist non-zero integers $k$ and $l$ such that $u^k$ is  conjugate to $v^l$ in $G$. The power-conjugacy problem has been the subject of extensive research, see  \cite{LM, AS, C, BM, BK, P, KM, Be, BDR}. However to the best of our knowledge, the interconnection of this problem and the classical conjugacy problem has not been studied yet.

\begin{theorem}\label{pc} (1) There is a finitely presented group $G$ with undecidable conjugacy problem but decidable power conjugacy problem.
Moreover $G$ has quadratic Dehn function.

(2) There is a finitely presented group $H$ with undecidable power conjugacy problem but decidable conjugacy problem.
\end{theorem}

Notice that for the group $G$ from Theorem \ref{pc}(1)  and \cite{OS20}, there
exists no algorithm recognizing the conjugacy of some {\it nontrivial} powers
of two elements (see Remark \ref{Wkk}) since elements of finite and
infinite orders behave differently in $G$. Although $G$ has undecidable
conjugacy problem, this problem is decidable in $G$ for elements of infinite order.
The following property of $G$ is used in the proof of Theorem \ref{pc} (1), and
it is also interesting in itself.

\begin{theorem}\label{c} For the group $G$ from Theorem \ref{pc} (1), there is an algorithm, which recognizes whether two elements $g$ and $h$ are conjugate in $G$ or not, provided the orders of $g$ and $h$ are infinite. The order of every element of $G$
can be also computed effectively.
\end{theorem}

Theorems \ref{isom}, \ref{subgr}, \ref{pc}, and \ref{c} are proved in Sections
\ref{iso},\ref{sub}, \ref{pow}, and \ref{Co}, respectively. The information needed
for understanding the proofs, has been selected from earlier papers and placed
in Sections \ref{mp} and \ref{gdp}.

\section{Machine preliminaries}\label{mp}

\subsection{S-machines}

Here we will use  definitions which equivalent to the definitions  used in \cite{OS06}) and \cite{O12}.

A \index[g]{S@$S$-machine!hardware of an $S$-machine}"hardware" of an $S$-machine $\bf S$ is a pair $(Y,Q),$ where $Q=\sqcup_{i=0}^nQ_i$ and $Y= \sqcup_{i=1}^n Y_i$ for some $n\ge 1$. Here and below $\sqcup$ denotes the disjoint union of sets.

We always set $Y_{n}=Y_0=\emptyset$ and if $Q_{n}=Q_0$ (i.e., the indices of $Q_i$ are counted $\mod n$), then we say that $\bf S$ is a \index[g]{S@$S$-machine!circular} {\em circular} $S$-machine.

The elements from $Q$ are called \index[g]{S@$S$-machine!state letters of an $S$-machine}{\it state letters}, the elements from $Y$ are \index[g]{S@$S$-machine!tape letters of an $S$-machine}{\it tape letters}. The sets $Q_i$ (resp.
$Y_i$) are called \index[g]{S@$S$-machine!parts of state and tape letters of an $S$-machine}{\em parts} of $Q$ (resp. $Y$).

The {\it language of \index[g]{S@$S$-machine!admissible words of an $S$-machine} admissible words}
consists of reduced words $W$ of the form
\begin{equation}\label{admiss}
q_1u_1q_2\dots u_s q_{s+1},
\end{equation}
where every $q_i$  is a state letter from some part $Q_{j(i)}^{\pm 1}$, $u_i$ are reduced group words in the alphabet of tape letters of the part $Y_{k(i)}$ and for every $i=1,...,s$ one of the following holds:

\begin{itemize}
\item If $q_i$ is from $Q_{j(i)}$ then $q_{i+1}$ is either from $Q_{j(i)+1}$ or is equal to $q_i\iv$; moreover $k(i)=j(i)+1$.

\item If $q_i\in Q_{j(i)}\iv$ then $q_{i+1}$ is either from $Q_{j(i)-1}^{-1}$ or is equal to $q_i\iv$; moreover $k(i)=j(i)$.
\end{itemize}
Every subword $q_iu_iq_{i+1}$ of an admissible word (\ref{admiss}) will be called the \index[g]{S@$S$-machine!admissible words of an $S$-machine!sector of an admissible word} {\em $Q_{j(i)}^{\pm 1}Q_{j(i+1)}^{\pm 1}$-sector} of that word. An admissible word may contain many $Q_{j(i)}^{\pm 1}Q_{j(i+1)}^{\pm 1}$-sectors.

We denote by $||W||$ the length of a word $W$. For every word $W$, if we delete all non-$Y^{\pm 1}$ letters from $W$ we get the \index[g]{Y@$Y$-projection of a word} $Y$-projection of the word $W$. The length of the $Y$-projection of $W$ is called the  \index[g]{Y@$Y$-length of a word} $Y$-{\em length} and is denoted by $|W|_Y$.
Usually parts of the set $Q$ of state letters are denoted by capital letters. For example, a part $P$ would consist of
letters $p$ with various indices.

If an admissible word $W$
has the form (\ref{admiss}), $W=q_1u_1q_2u_2...q_s,$
and $q_i\in Q_{j(i)}^{\pm 1},$
$i=1,...,s$, $u_i$  are
group words in tape letters, then we shall say that the \index[g]{S@$S$-machine!admissible words of an $S$-machine!base of an admissible word}{\em base} of $W$ is the word
$Q_{j(1)}^{\pm 1}Q_{j(2)}^{\pm 1}...Q_{j(s)}^{\pm 1}$. Here $Q_i$ are just symbols which denote the corresponding parts of the set of state letters. Note that, by the definition of admissible words, the base is not necessarily a reduced word.

Instead of saying that the parts of the set of state letters of $\sss$ are $Q_0, Q_1, ... , Q_n$ we will write
that the \index[g]{S@$S$-machine!standard base of an $S$-machine}{\em the standard base}
of the $S$-machine is $Q_0...Q_n$.

The \index[g]{S@$S$-machine!software of an $S$-machine}\emph{software} of an $S$-machine with the standard base $Q_0...Q_n$ is a set of \index[g]{S@$S$-machine!rule of an $S$-machine}{\em rules} $\Theta$.  Every $\theta\in \Theta$
is a sequence $[q_0\to a_0q_0'b_0,...,q_n\to a_nq_n'b_n]$ and a subset $Y(\theta)=\sqcup Y_j(\theta)$, where $q_i\in Q_i$,
$a_i$ is a reduced word in the alphabet $Y_{i-1}(\theta)$, $b_i$ is a reduced word in $Y_i(\theta)$, $Y_i(\theta)\subseteq Y_i$, $i=0,...,n$ (recall that $Y_{0}=Y_n=\emptyset$).

Each component $q_i\to a_iq_i'b_i$ is called a \index[g]{S@$S$-machine!rule of an $S$-machine!part of a rule}{\em part} of the rule.
In most cases the sets $Y_j(\theta)$ will be equal to
either $Y_j$ or $\emptyset$. By default $Y_j(\theta)=Y_j$.

Every rule $$\theta=[q_0\to a_0q'_0b_0,...,q_n\to a_nq'_nb_n]$$ has an inverse
$$\theta\iv=[q_0'\to a_0\iv q_0b_0\iv,..., q_n'\to a_n\iv q_nb_n]$$ which is also a rule of $\sss$. It is always the case that
$Y_i(\theta\iv)=Y_i(\theta)$ for every $i$. Thus the set of rules
$\Theta$ of an $S$-machine is divided into two disjoint parts, $\Theta^+$ and
$\Theta^-$ such that for every $\theta\in \Theta^+$, $\theta\iv\in
\Theta^-$ and for every $\theta\in\Theta^-$, $\theta\iv\in\Theta^+$ (in particular $\Theta\iv=\Theta$, that is any $S$-machine is symmetric).

The rules from $\Theta^+$ (resp. $\Theta^-$) are called \index[g]{S@$S$-machine@rule of an $S$-machine!positive or negative}{\em
positive} (resp. {\em negative}).

 To \index[g]{S@$S$-machine!rule of an $S$-machine!application of a rule}apply a rule  $\theta=[q_0\to a_0q'_0b_0,...,q_n\to a_nq'_nb_n]$ as above to an admissible word $p_1u_1p_2u_2...p_s$ (\ref{admiss})
where each $p_i\in Q_{j(i)}^{\pm 1}$ means
\begin{itemize}

\item check if $u_i$ is a word in the alphabet $Y_{j(i)+1}(\theta)$ when $p_i\in Q_{j(i)}$ or if it is a word
    in $Y_{j(i)}(\theta)$ when $p_i\in Q_{j(i)}^{-1}$ ($i=1,\dots,s-1$); and if this property holds,

\item replace each $p_i=q_{j(i)}^{\pm 1}$ by $(a_{j(i)}q'_{j(i)}b_{j(i)})^{\pm 1}$,
\item if the resulting word is not reduced or starts (ends) with $Y$-letters, then reduce the word and trim the first and last $Y$-letters to obtain an admissible word again.
\end{itemize}

 If a rule $\theta$ is applicable to an admissible word $W$ (i.e., $W$ belongs to the \index[g]{S@$S$-machine!rule of an $S$-machine!domain of a rule}{\em domain} of $\theta$) then we say that $W$ is a \index[g]{S@$S$-machine!$\theta$-admissible word} $\theta$-{\it admissible word} and denote the result of application of $\theta$ to $W$ by $W\cdot \theta$.
 Hence each rule defines an invertible partial map from the set of admissible words to itself, and one can consider an $S$-machine as an inverse semigroup of partial bijections of the set of admissible words.

We call an admissible word with the standard base a  \index[g]{S@$S$-machine!configuration of an $S$-machine} \emph{configuration} of an $S$-machine.

We usually assume that every part $Q_i$ of the set of state letters contains a \index[g]{S@$S$-machine!start state letter of an $S$-machine} {\em start state letter} and an
\index[g]{S@$S$-machine!end state letter of an $S$-machine} {\em end state letter}. Then a configuration is called a \index[g]{S@$S$-machine!start configuration of an $S$-machine}\index[g]{S@$S$-machine!end configuration of an $S$-machine} \emph{start} (\emph{end}) configuration if all state letters in it are start (end) letters. As Turing machines, some $S$-machines are \index[g]{S@$S$-machine!recognizing a language} {\em recognizing a language}. In that case we choose an \index[g]{S@$S$-machine!recognizing a language!input sector of an admissible word of an $S$-machine} {\em input} sector, usually the $Q_0Q_1$-sector, of every configuration. The $Y$-projection of that sector is called the \index[g]{S@$S$-machine!recognizing a language!input of a configuration of an $S$-machine recognizing a language} {\em input} of the configuration. In that case, the end configuration with empty $Y$-projection is called \index[g]{S@$S$-machine!recognizing a language!accept configuration of an $S$-machine recognizing a language} the \emph{accept} configuration. If the $S$-machine (viewed as a semigroup of transformations as above) can take an input configuration with input $u$ to the accept configuration, we say that $u$ is \index[g]{S@$S$-machine!recognizing a language!accepted input word} {\em accepted} by the $S$-machine. We define \index[g]{S@$S$-machine!recognizing a language!accepted configuration of an $S$-machine} {\em accepted} configurations (not necessarily start configurations) similarly.

A \index[g]{S@$S$-machine!computation of an $S$-machine}{\em computation} of  \index[g]{S@$S$-machine!computation of an $S$-machine!length of a computation} {\it length} $t\ge 0$ is a sequence of admissible words $$W_0\stackrel{\theta_1}{\to} \dots\stackrel{\theta_{t}}{\to} W_t$$
such that for
every $i=0,..., t-1$ the $S$-machine passes from $W_i$ to $W_{i+1}$ by applying
the rule $\theta_i$ from $\Theta$.  The word $H=\theta_1\dots\theta_{t}$ is called the \index[g]{S@$S$-machine!computation of an $S$-machine!history of computation} {\it history}
of the computation, and the word $W_0$ is called \index[g]{S@$S$-machine!$H$-admissible word} $H$-{\it admissible}. Since $W_t$ is determined by $W_0$ and the history $H$, we use notation $W_t=W_0\cdot H$.

A computation is called \index[g]{S@$S$-machine!computation of an $S$-machine!reduced} {\em reduced} if its history is a reduced word.

Note, though, that in \cite{OS20} and in this paper, unlike the previous ones, we consider non-reduced computations too because these may correspond to reduced van Kampen diagrams (tra\-pe\-zia) under our present interpretation of $S$-machines in groups.

If for some rule $\theta=[q_0\to a_0q_0'b_0,...,q_n\to a_nq'_nb_n]\in \Theta$ of an $S$-machine $\sss$ the set  $Y_{i+1}(\theta)$ is empty (hence in every admissible word in the domain of $\theta$ every $Q_iQ_{i+1}$-sector has no $Y$-letters) then we say that $\theta$ \index[g]{S@$S$-machine!rule of an $S$-machine!locking a sector} locks the $Q_iQ_{i+1}$-sector. In that case we always assume that $b_i, a_{i+1}$ are empty and we denote the $i$-th part of the rule $q_i\tool a_iq_i'$. If the $Q_iQ_{i+1}$-sector is locked by $\theta$ then we also assume that $a_{i+1}$ is empty too.

\begin{rk} \label{tool}
For the sake of brevity, the substitution $[q_i\tool aq_i', q_{i+1}\to q_{i+1}'b]$
can be written in the form $[q_iq_{i+1}\to aq'_iq'_{i+1}b]$.
\end{rk}

The above definition of $S$-machines resembles the definition of multi-tape Turing machines (see \cite{SBR}).
The main differences are that every state letter of an $S$-machines is blind: it does not "see" tape letters next to it (two state letters can see each other if they stay next to each other). Also $S$-machines are symmetric (every rule has an inverse), can work with words containing negative letters, and words with "non-standard" order of state letters.

It is important that $S$-machines can simulate the work of Turing machines. This non-trivial fact, especially if one tries to get a polynomial time simulation, was first proved in \cite{SBR}. But we do not need a restriction on time, and it would be more convenient for us to use an easier $S$-machine from \cite{OS06}.

Let \index[g]{Turing machine $\mmm_0$}$\mmm_0$ be a deterministic Turing machine accepting a
non-empty
language $\lll$ of words in the one-letter alphabet $\{\alpha\}$. In different sections we will use two versions  of an equivalent S-machine ${\bf M}_1$. The first version is borrowed from \cite{OS06}, where Lemmas 3.25 and 3.27 of \cite{OS06} imply the following
property.

\begin{lemma} \label{negat} There is a recognizing $S$-machine $\mmm_1$ whose language of accepted input words is $\lll$. In every input configuration of $\mmm_1$, there is exactly one input sector, the first sector of the word, and all other sectors are empty of $Y$-letters.

If a non-empty reduced computation $C_0\to\dots \to C_t$
of ${\bf M}_1$  starts with an
input configuration containing a negative letter
then $C_t$ is neither an input nor the accept configuration.
\end{lemma} $\Box$

The following statement can be found in \cite{O12} (Lemmas 4.15 and 4.16(a)), although below we denote the machine by $\mmm_1$ instead of $M_2$ in \cite{O12}.

\begin{lemma} \label{S}  There is a recognizing $S$-machine $\mmm_1$ whose language of accepted input words is $\lll$. In every input configuration of $\mmm_1$, there is exactly one input sector, the first sector of the word, and all other sectors are empty of $Y$-letters.

For every reduced computation $W_0\to\dots\to W_t$ of  $\mmm_1$  with the standard
base and a non-empty history $H$, we have $W_t\ne W_0$.
\end{lemma} $\Box$

Lemma \ref{wi} is a modified formulation of Lemma 2.8 from \cite{OS20}.
(In \cite{OS20}, it was formulated for reduced computations, but the proof
did not use that the history was reduced.)

\begin{lemma} \label{wi} Suppose that a computation $W_0\to W_1\to\dots \to W_t$ of an $S$-machine $\sss$
has a $2$-letter base
and the history of the form  $H \equiv  H_1H_2^k H_3$ ($k\ge 0$). Then for every $i=0,1,\dots,t$,
we have the inequality $$||W_i||\le ||W_0|+ ||W_t||+2||H_1||+ 3||H_2||+2||H_3||$$.
\end{lemma} $\Box$

Recall that a word $w$ is called a {\it periodic word} with period $v$
if $w$ is a subword of some power of $v$.

\begin{lemma} \label{perio} There is an exponential function $f$ with the following property.

Suppose a computation ${\cal C}: W_0\to W_1\to\dots \to W_t$ of an $S$-machine $\sss$
has a periodic history with period $H$.
Assume that $\cal C$  has no subcomputations $W_i\to\dots \to W_j$
with history $H$ and  $W_i\equiv W_j$. Then $t\le f(||W_0||(||W_0||+||W_t||+||H||))$.
\end{lemma} $\Box$

\proof Since the history is $H$-periodic, there are words $W_{i_1},\dots, W_{i_s}$, where $i_{k+1}-i_{k}=||H||$ ($k=1,\dots, s-1$), the history
of every subcomputations $W_{i_k}\to W_{i_k+1}\dots\to W_{i_{k+1}}$ is $H$, and $s\ge t ||H||^{-1} - 1$.

Assume that $W_{i_k}\equiv W_{i_l}$ for some $l>k$. Then we have
$V_{i_k}\equiv V_{i_l}$ for arbitrary restriction of $\cal C$ to a
subbase $B$ of length $2$. Arbitrary computation with base $B$ and
history $H$ multiply the  $Y$-projection $v$ from the left by a word $a$
and from the right by a word $b$, where the words $a$ and $b$ depend
on $B$ and $H$ only. Therefore for the $Y$-projection $v$
of the equal words $V_{i_k}$ and  $V_{i_l}$, we obtain the equality $v=a^mvb^m$, where $m=l-k\ge 1$. Hence we have $(v^{-1}a^{-1}v)^m=b^m$, which implies in the free group that $v^{-1}a^{-1}v=b$, i.e. $avb=v$.
Hence $V_{i_k}\equiv V_{i_{k+1}}$ for every 2-letter subbase $B$.
It follows that $W_{i_k}\equiv W_{i_{k+1}}$, contrary to the lemma assumption.

Therefore we obtain $s$ different admissible words in the computation
$\cal C$.  Lemma \ref{wi} bounds their lengths  by a liner function of $||W_0||(||W_0||+||W_t||+||H||)$ since every word $W_i$ is covered by
at most $||W_0||$ admissible subwords with 2-letter bases. Hence the number $s$ and the number $t\le(s+1)||H||$ are bounded from above by an exponential function.
\endproof

\subsection{Running state letters} \label{pm}

For every alphabet $Y$ we define a "running state letters" $S$-machine \index[g]{S@$S$-machine!L@$\rhh$}$\rhh(Y)$. We will omit $Y$ if it is obvious or irrelevant. The standard base of $\rhh(Y)$ is $Q^{(1)}PQ^{(2)}$ where $Q^{(1)}=\{q^{(1)}\}$, $P=\{p^{(i)}, i=1,2\}$, $Q^{(2)}=\{q^{(2)}\}$. The state letter $p$ with indices runs from the state letter $q^{(2)}$ to the state letter $q^{(1)}$ and back.  The $S$-machine $\rhh$ will be used to check the "structure" of a configuration (whether the state letters of a configuration are in the appropriate order), and to recognize a computation by its history.

The alphabet of tape letters $Y$ of $\rhh(Y)$ is $Y^{(1)}\sqcup Y^{(2)}$, where $Y^{(2)}$
is a (disjoint) copy of $Y^{(1)}$.
The positive rules of $\rhh$ are defined
as follows.
\begin{itemize}

\item $\zeta^{(1)}(a)=[q^{(1)}\to q^{(1)}, p^{(1)}\to a\iv p^{(1)}a', q^{(2)}\to q^{(2)}]$,
where $a$ is any positive letter from $Y=Y^{(1)}$ and $a'$ is the corresponding letter in the copy $Y^{(2)}$ of $Y^{(1)}$.
\me

{\em Comment.} The state letter $p^{(1)}$ moves left replacing letters $a$ from $Y^{(1)}$ by their copies $a'$
from $Y^{(2)}$.

\me

\item $\zeta^{(12)}=[q^{(1)}p^{(1)}\to  q^{(1)}p^{(2)}, q^{(2)}\to q^{(2)}]$.

\me

{\em Comment.} When $p^{(1)}$ meets $q^{(1)}$, $p^{(1)}$  turns into $p^{(2)}$.

\me

\item $\zeta^{(2)}(a) =[q^{(1)}\to q^{(1)}, p^{(2)}\to ap^{(2)}(a')^{-1}, q^{(2)}\to q^{(2)}]$

{\em Comment.} The state letter $p^{(2)}$ moves right towards $q^{(2)}$ replacing letters $a'$ from $Y^{(2)}$ by their copies $a$
from $Y^{(1)}$.
\end{itemize}

The start (resp. end) state letters of $\rhh$ are $\{q^{(1)},p^{(1)}, q^{(2)}\}$ (resp. $\{q^{(1)},p^{(2)}, q^{(2)}\}$).

\begin{lemma}\label{Hprim} (\cite{OS20}, Lemma 3.14 (b))
Let $\ccc\colon
W_0\to\dots\to W_t$ be
a reduced computation of  $\mmm_3$ consisting of rules of one of the copies of $\rhh$ or $\rhr$ with standard base. Then
$ t\le ||W_0||+||W_t||-2$.
\end{lemma} $\Box$

\begin{rk} \label{lrm} For some large integer $m$, we will also need the $S$-machine \index[g]{S@$S$-machine!L@$\rhh_m$}$\rhh_m$ from \cite{OS20}, that repeats the
work of $\rhh$ $m$ times. That is the $S$-machine $\rhh_m$ runs the state letter $p$ back and forth between $q^{(2)}$ and $q^{(1)}$ $m$ times. Every time $p$ meets $q^{(1)}$ or $q^{(2)}$, the upper index of $p$ increases by $1$ after the application of the rule $\zeta^{(i,i+1)}$ ($i=1,\dots, 2m-1$), so the highest upper index of $p$ is $(2m)$.
\end{rk}

\begin{rk}\label{para}
$m$ is one of the  big constants used in \cite{OS20}.
In this paper we will use just few of them. Here they are:
\begin{equation}\label{const} \begin{array}{l}  m, \; N \ll c_4\ll L \end{array}
  \end{equation}
where $\ll$ means "much smaller".
\end{rk}

\subsection{Adding history sectors} \label{hs}

We will add new (history) sectors to an $S$-machine $\mmm_1$ provided
by Lemma \ref{negat} or by Lemma \ref{S}. The history sectors split the base letters of $\bf M_1$. (See the definition below). If we ignore the new sectors, in essence, we get the hardware and the  software of the $S$-machine $\mmm_1$. The new $S$-machine $\mmm_2$ will start with a configuration where in every history sector a copy of the history $H$ of a computation of $\mmm_1$ is written. Then it will execute $H$ on the other (working) sectors simulating the work of $\mmm_1$, while  in the history sector, state letters scan the history, one symbol at a time. Thus if a computation of $\mmm_2$ with the standard base starts with a configuration $W$ and ends with configuration $W'$, then the length of the computation does not exceed $||W||+||W'||$.

Here is a precise definition of $\mmm_2$. Let the $S$-machine $\mmm_1$
have hardware $(Q,Y)$, where
$Q=\sqcup_{i=0}^n Q_i$, and the set  of rules $\Theta$. The new
$S$-machine \index[g]{S@$S$-machine!M2a@$\mmm_2$} $\mmm_2$ has hardware
$$Q_{0,r}\sqcup Q_{1,\ell}\sqcup Q_{1,r}\sqcup Q_{2,\ell}\sqcup Q_{2,r}\sqcup\dots \sqcup Q_{n,\ell},\;\; Y_h=  Y_1
\sqcup X_1\sqcup Y_2\sqcup\dots \sqcup X_{n-1}\sqcup Y_n$$
where $Q_{i,\ell}$ and $Q_{i,r}$ are (left and right) copies of $Q_i$  and $X_i$ is a disjoint union of two copies of $\Theta^+$, namely $X_{i,\ell}$ and $X_{i,r}$. (The sets $Q_{0,\ell}$, $Q_{n,r}$ %$X_{0,.}$, $X_{n,.}$
are empty.) Every letter $q$ from $Q_i$ has two copies $q^{(\ell)}\in Q_{i,\ell}$ and $q^{(r)}\in Q_{i,r}$. By definition, the start (resp. end) state letters  of $\mmm_2$ are copies of the corresponding start (end) state letters of $\mmm_1$.  The $Q_{0,r}Q_{1,\ell}$-sectors are the \index[g]{S@$S$-machine!M2a@$\mmm_2$!input sector of a configuration of $\mmm_2$} {\em input sectors} of configurations of $\mmm_2$.

The positive rules $\theta_h$ of $\mmm_2$ are in one-to-one correspondence with the positive rules $\theta$ of $\mmm_1$. If $\theta=[q_0\to a_0q_0'b_0,...,q_n\to a_nq_n'b_n]$ is a positive rule of $\mmm_1$, then each part $q_i\to a_iq_i'b_i$ is replaced in $\theta_h$ by two parts $$q_{i,\ell}\to a_iq_{i,\ell}'h_{\theta,i}\iv $$ and
$$q_{i,r}\to \overline h_{\theta,i}q_{i,r}' b_i,$$
where $h_{\theta,i}$ (resp., $\overline h_{\theta,i}$) is a copy of $\theta$ in the alphabet $X_{i,\ell}$ (in $X_{i,r}$, respectively).

If $\theta$ is the start (resp. end) rule of $\mmm_1$, then for any word in the domain of $\theta_h$ (resp. $\theta_h\iv$) all $Y$-letters in history sectors are from $\sqcup_i X_{i,\ell}$ (resp. $\sqcup X_{i,r}$).

Thus for every rule $\theta$ of $\mmm_1$, the rule $\theta_h$ of $\mmm_2$ acts in the $Q_{i,r}Q_{i+1,\ell}$-sector in the same way as $\theta$ acts in the $Q_iQ_{i+1}$-sector. In particular, $Y$-letters which can appear in the $Q_{i,r}Q_{i+1,\ell}$-sector of an admissible word in the domain of $\theta_h$ are the same as the $Y$-letters  that can appear in the $Q_iQ_{i+1}$-sector of an admissible word in the domain of $\theta$. Hence if $\theta$ locks $Q_iQ_{i+1}$-sectors, then $\theta_h$ locks $Q_{i,r}Q_{i+1,\ell}$-sectors.

\begin{rk}
Every computation of the $S$-machine $\mmm_2$ with history $H$ and the standard base coincides with the a computation of $\mmm_1$ whose history is a copy of $H$
if
one observes it only in sectors $Q_{i,r}Q_{i+1,l}$.
\end{rk}

Let \index[g]{S@$S$-machine!M1@$\mmm_1$!I@$I_1(\alpha^k)$ - a start configuration of $\mmm_1$}$I_1(\alpha^k)$ be a start configuration of $\mmm_1$ (an input configuration in the domain of the start rule of $\mmm_1$) with $\alpha^k$ written in the input sector (all other sectors do not contain $Y$-letters), and $H$ be a word in the alphabet
of rules of ${\bf M}_1$. Then the corresponding start configuration \index[g]{S@$S$-machine!M2a@$\mmm_2$!I@$I_2(\alpha^k,H)$ - a start configuration of $\mmm_2$} $I_2(\alpha^k,H)$ of $\mmm_2$ is obtained by first replacing each state letter $q$ by the product of two corresponding letters $q^{(\ell)}q^{(r)}$, and then inserting a copy of $H$ in the \index[g]{X@$X_{i,\ell}$, a left alphabet} {\em left alphabet} $X_{i,\ell}$ in every history $Q_{i,\ell}Q_{i,r}$-sector. End configurations \index[g]{S@$S$-machine!M2a@$\mmm_2$!A@$A_2(H)$ - an end configuration of $\mmm_2$} $A_2(H)$ of $\mmm_2$ are defined similarly, only the $Y$-letters in the history sectors must be from the \index[g]{X@$X_{i,r}$, a right alphabet} {\em right alphabet} $X_{i,r}$.
\label{r:09}

\subsection{Adding running state letters}\label{csl}

Our next $S$-machine will be a composition of $\mmm_2$ with $\rhh$ and $\rhr$. The running state letters will control the work of $\mmm_3$.

First we replace every part $Q_i$ of the state letters in the standard base of $\mmm_2$ by three parts $P_iQ_iR_i$ where $P_i, R_i$ contain the running state letters. Thus if $Q_0...Q_s$ is the standard base of $\mmm_2$ then the standard base of  \index[g]{S@$S$-machine!M2b@$\bmm_2$} $\bmm_2$ is
\begin{equation}\label{basec}
P_0Q_0R_0P_1Q_1R_1\dots P_sQ_sR_s,
\end{equation}
where
$P_i$ (resp., $R_i$) contains
copies of running $P$-letters (resp. $R$-letters) of $\rhh$ (resp. $\rhr$), $i=0,\dots, s$.

For every rule $\theta$ of  $\mmm_2$,
its $i$-th part $[q_i\to a_iq_i'b_i]$ is replaced in $\bmm_2$ with
\begin{equation}\label{3p}
[p^{(i)}q_ir^{(i)}\to a_ip^{(i)}q_i'r^{(i)}b_i],  (i=0,\dots, s),
\end{equation}
 where $ p^{(i)}\in P_i, r^{(i)}\in R_i$ do not depend on $\theta$.

 \medskip

{\it Comment.} Thus, the sectors $P_iQ_i$ and $Q_iR_i$ are always locked. Of course, such a modification is useless for solo work of $\mmm_2$. But it will be helpful when one constructs
a composition of $\bmm_2$ with $\rhh$ and $\rhr$ which will be turned on after certain rules of $\bmm_2$ are applied.

If $Q_{i-1}Q_i$ is an input sector of configurations of the machine ${\bf M}_2$, then $R_{i-1}P_i$ is an
input sector of the configurations of ${\bf\bar M}_2$.

\subsection{The machine $M_3$}\label{M3M5}

The next $S$-machine \index[g]{S@$S$-machine!M3@$\mmm_3$} $\mmm_3$ is the composition of the $S$-machine $\mathbf{\overline M}_2$ with $\rhh$ and $\rhr$.
The $S$-machine $\mmm_3$ has the same base as $\bmm_2$, although the parts of this base have more
state letters than the corresponding parts of $\mathbf{\overline M}_2$. It works as follows. Suppose that $\mmm_3$ starts with a start configuration of $\bmm_2$, a word $\alpha^k$ in the input $R_0P_1$-sector, copies of a history word $H$ in the alphabets $X_{i,\ell}$ in the history sectors, all other sectors empty of $Y$-letters. Then $\mmm_3$ first executes $\rhr$ in all history sectors (moves the running state letter from $R_i$ in the history sectors right and left), then it executes the history $H$ of $\bmm_2$. After that the $Y$-letters in the history sectors are in $X_{i,r}$ and $\mmm_3$ executes copies of $\rhh$ in the history sectors (moves the running state letters left then right). After that $\mmm_3$ executes a copy of $H$ backwards, getting to a copy of the same start configuration of $\bmm_2$, runs $\rhr$, executes a copy of the history $H$ of $\bmm_2$, runs a copy of $\rhh$, etc. It stops after $m$ times running $\rhr, \bmm_2,\rhh, \bmm_2\iv$ and running $\rhr$ one more time.

Thus the $S$-machine $\mmm_3$ is a concatenation of $4m+1$ $S$-machines $\mmm_{3,1}-\mmm_{3,4m+1}$. After one of these $S$-machines terminates, a transition rule changes its end state letters to the start state letters of the next $S$-machine. All these $S$-machines have the same standard bases as $\bmm_2$.

The configuration $I_3(\alpha^k,H)$ of ${\bf M}_3$ is obtained from $I_2(\alpha^k,H)$
by adding the control state letters $r^{(1)}_i$ and $p^{(1)}_i$ according to (\ref{basec}) in Section \ref{csl}.

\medskip

{\bf Set $\mmm_{3,1}$} is a copy of the set of rules of the $S$-machine $\rhr$, with \index[g]{S@$S$-machine!L@$\rhh$!parallel work of $\rhh$ or $\rhr$ in several sectors} \emph{parallel work} in all history sectors, i.e., every subword $Q_{i-1}R_{i-1}P_i$ of the standard base, where $Q_{i-1}Q_i$ is a history sector of $\mmm_2$,
is treated as the base of a copy of $\rhr$, that is $R_{i-1}$ contain the running state letters which run between state letters from $Q_{i-1}$ and $P_i$. Each rule of Set $\mmm_{3,1}$ executes the corresponding rule of $\rhr$ simultaneously in each history sector of $\mmm_2$. The partition of the set of state letters of these copies of $\rhr$ in each history sector is $X_{i,\ell}\sqcup X_{i,r}$ for some $i$ (that is state letters from $R_{i-1}$ first run right replacing letters from $X_{i,\ell}$ by the corresponding letters of $X_{i,r}$ and then run left replacing letters from $X_{i,r}$ by the corresponding letters of $X_{i,\ell}$.

The  transition rule $\chi(1,2)$ changes the state letters to the state letters of start configurations of $\bmm_2$. The admissible words in the domain of $\chi(1,2)^{\pm 1}$
have all $Y$-letters from the left alphabets $X_{i,\ell}$. The rule $\chi(1,2)$ locks all sectors except the history sectors $R_{i-1}P_i$ and the input sector. It does not apply to admissible words containing $Y$-letters from right alphabets.

\medskip

{\bf Set  $\mmm_{3,2}$} is a copy of the set of rules of the $S$-machine $\bmm_2$.

The transition rule $\chi(2,3)$ changes the state letters of the stop configuration of $\bmm_2$ to  their copies in a different alphabet. The admissible words in the domain of $\chi(2,3)^{\pm 1}$
have no $Y$-letters from the left alphabets $X_{i,\ell}$. The rule $\chi(2,3)$ locks all sectors except for the history sectors $R_{i-1}P_i$. It does not apply to admissible words containing $Y$-letters from right alphabets.

\medskip

{\bf  Set $\mmm_{3,3}$} is a copy of the set of rules of the $S$-machine $\rhh$, with parallel work in the same sectors as  $\mmm_{3,1}$ (and the same partition of $Y$-letters in each history sector $X_{i,r}\sqcup X_{i,\ell}$).

The transition rule $\chi(3,4)$ changes the state letters of the stop configuration of $\bmm_2$ to  their copies in a different alphabet. The admissible words in the domain of $\chi(3,4)^{\pm 1}$
have no $Y$-letters from the left alphabets $X_{i,l}$. The rule $\chi(3,4)$ locks all non-history sectors.

\medskip

{\bf Set $\mmm_{3,4}$.} The positive rules of Set $\mmm_{3,4}$ are the copies of the negative rules of the $S$-machine  $\bmm_2$.

The transition rule $\chi(4,5)$ changes the state letters of the start configuration of $\bmm_2$ to  their copies in a different alphabet. The admissible words in the domain of $\chi(4,5)^{\pm 1}$
have no $Y$-letters from the right alphabets $X_{i,r}$. The rule $\chi(4,5)$ locks all  non-history and non-input sectors.

\medskip

{\bf Sets $\mmm_{3,5}, \dots, \mmm_{3,8}$} consist of rules that are  copies of the rules of the Sets $\mmm_{3,1},\dots,$ $\mmm_{3,4}$, respectively.

$\dots$

{\bf Sets $\mmm_{3,4m-3},\dots, \mmm_{3,4m}$} consist of copies of the steps $\mmm_{3,1}, \dots, \mmm_{3,4}$, respectively.

\medskip

{\bf Set $\mmm_{3,4m+1}$} is a  copy of Set $\mmm_{3,1}$. The end configuration for Set $\mmm_{3,4m+1}$,
$A_3(H)$, is obtained from a copy of $A_2(H)$ by inserting the control letters according to (\ref{basec}).

The transition rules $\chi(i,i+1)$ are called \index[g]{S@$S$-machine!M3@$\mmm_3$!Chi@$\chi$-rules of $\mmm_3$}
$\chi$-rules.

\begin{lemma}\label{M31} (\cite{OS20}, Lemma 3.15) Let $\ccc\colon W_0\to\dots\to W_t$ be a reduced computation
of $\mmm_3$ with the standard base. Then for every $i$, there is at most one occurrence of the rules
$\chi(i,i+1)^{\pm 1}$ in the history $H$ of $\ccc$.
\end{lemma} $\Box$

\subsection{\texorpdfstring{$\mmm_4$}{M4} and \texorpdfstring{$\mmm_5$}{M5}}\label{45}

Let $B_3$ be the standard base of $\mmm_3$ and $B_3'$ be its disjoint copy. By \index[g]{S@$S$-machine!M4@$\mmm_4$} $\mmm_4$ we denote
the $S$-machine with standard base $B_3(B'_3)^{-1}$ and rules $\theta(\mmm_4)=[\theta, \theta]$,
where $\theta\in \Theta$ and $\Theta$ is the set of rules of $\mmm_3$. So the rules of $\Theta(\mmm_4)$ are the same
for $\mmm_3$-part of $\mmm_4$ and for the mirror copy of $\mmm_3$. Therefore we will denote $\Theta(\mmm_4)$ by
$\Theta$ as well.  The sector between the last state letter of $B_3$ and the first state letter of $(B_3')^{-1}$ is locked by any rule from $\Theta$.
(The 'mirror' symmetry of the base  is used in \cite{OS20} for the upper estimate of the Dehn function.)

\medskip

The $S$-machine \index[g]{S@$S$-machine!M5@$\mmm_5$}$\mmm_5$ is a circular analog of $\mmm_4$. We add one more base letter $\tt$ to the hardware of $\mmm_4$. So the standard base $B$ of $\mmm_5$ it $\{\tt\}B_3(B_3')^{-1}\{\tt\}$, where the part $\{\tt\}$ has only one letter $\tt$
and the first part $\{\tt\}$ is identified with the last part. For example, $\{\tt\}B_3(B_3')^{-1}\{\tt\}B_3(B_3')^{-1}$ can be a base of an admissible word for $\mmm_5$. Furthermore,  sectors involving $\tt^{\pm 1}$ are
 locked by every rule from $\Theta$. For $\mmm_5$, we have the start and stop words $I_5(\alpha^k,H)$ and $A_5(H)$ similar to the configurations
$I_3(\alpha^k,H)$ and $A_3(H)$.

\subsection{The main machine \texorpdfstring{$\mmm$}{M}} \label{M6}

We use the $S$-machine $\mmm_5$ from Section \ref{45}, $\rhh_m$ from Section \ref{pm} and three more easy $S$-machines to compose
the main circular $S$-machine \index[g]{S@$S$-machine!M@$\mmm$} $\mmm$ needed for this paper. The standard base of $\mmm$ is the same as the standard base of $\mmm_5$, i.e.,$\{\tt\}B_3(B_3')^{-1}$, where $B_3$ has the form (\ref{basec}).
We will use ${\tilde Q}_0$ instead of $Q_0$, ${\tilde R}_1$ instead of $R_1$ and so on to denote parts of the set of state letters since $\mmm$ has more
state letters in every part of its hardware.

The rules of $\mmm$ will be partitioned into five sets ($S$-machines) $\mtt_i$ ($i=1,\dots,5$) with transition rules $\theta(i,i+1)$ \index[g]{S@$S$-machine!M@$\mmm$!transition rules $\theta(i,i+1)$} connecting $i$-th and $i+1$-st sets.
The state letters are also disjoint for different sets $\mtt_i$. It will be clear that  ${\tilde Q}_0$ is the disjoint union of 5 disjoint sets including $Q_0$, ${\tilde R}_1$ is the disjoint union of five disjoint sets including $R_1$, etc.

By default, every transition rule $\theta(i,i+1)$ of $\mmm$ locks a sector if this sector
is locked by all rules from $\mtt_i$ or if it is locked by all rules from $\mtt_{i+1}$.
It also changes the end state letters of $\mtt_i$ to the start state letters of $\mtt_{i+1}$.

The \index[g]{start configuration $W_{st}$ of $\bf M$} {\it start configuration} $W_{st}$ of $\mtt$ is $\tilde t b_3(b'_3)^{-1}$,
where $b_3$ and $b'_3$ are obtained by replacing every base letter of $B_3$ and $B'_3$
by special start letter. The start rule $\theta_1$ of $\bf M$ changes the letters from $b_3$ and $b'_3$ to
their copies and starts the work of the rules from the set {\bf $\mtt_1$}.

\medskip

 Set {\bf $\mtt_1$} inserts input words in the input sectors. The set contains only one positive rule inserting the letter $\alpha$ in the input sector next to the left of a letter $p$ from ${\tilde P}_1$. It also inserts a copy $\alpha\iv$  next to the right of the corresponding letter $(p')\iv$ (the similar
mirror symmetry is assumed in the definition of all other rules.)
So the positive rule of $\mtt_1$ has the form $$[q_0\tool q_0, r_1\to r_1, p_1\tool \alpha p_1, ..., (p_1')^{-1}\to  (p_1')^{-1}\alpha^{-1}, (r_1')^{-1}\tool (r_1')^{-1}, t\tool t]$$

The rules of $\mtt_1$ do not change state letters, so it has one state letter in each part of its hardware.

The connecting rule $\theta(12)$ changes the state letters of $\mtt_1$ to their copies in a disjoint alphabet. It locks all sectors except for the input sector ${\tilde R}_0{\tilde P}_1$ and the mirror copy of this sector.

\medskip

 Set {\bf $\mtt_2$} is a copy of the $S$-machine $\rhh_m$ working in the input sector and its mirror image in parallel, i.e.,we identify the standard base
of $\rhh_m$ with ${\tilde R}_0 {\tilde P}_1{\tilde Q}_1$. The connecting rule
$\theta(23)$ locks all sectors except for the input sector ${\tilde R}_0{\tilde P}_1$ and its mirror image.
\medskip

 Set {\bf $\mtt_3$} inserts history in the history sectors.
This set of rules is a copy of each of the left alphabets $X_{i,l}$ of
the $S$-machine $\mmm_2$. Every positive rule of $\mtt_3$ inserts a copy of the corresponding positive letter in every
history sector ${\tilde R}_i{\tilde P}_{i+1}$ next to the right of a state letter from ${\tilde R}_i$.

Again, $\mtt_3$ does not change the state letters, so each part of its hardware contains one letter.

The transition rule $\theta(34)$ changes the state letters to their copies from Set $\mmm_{5,1}$ of $\mmm_5$.
It locks all sectors except for the input sectors and the history sectors. The history sectors in admissible words from the domain of $\theta(34)$ have  $Y$-letters from the left
alphabets $X_{i,l}$ of the $S$-machine $\mmm_5$.

\medskip

 Set {\bf $\mtt_4$} is a copy of  the $S$-machine $\mmm_5$. The transition rule
$\theta(45)$ locks all sectors except for history ones. The admissible words in the domain of
$\theta(45)$ have no letters from right alphabets.

\medskip

Set {\bf $\mtt_5$.} The positive rules from $\mtt_5$ simultaneously
erase the letters of the history sectors from the  right  of the state letter from ${\tilde R}_i$. That is, parts of the rules are of the form $r\to ra\iv $ where $r$ is a state letter from ${\tilde R}_i$, $a$ is a letter from the left alphabet of the history sector.

\medskip

Finally the accept rule $\theta_0$ (regarded as a transition rule) from $\mmm$ can be applied when all the sectors are empty, so it locks all the sectors and changes the end state letters of $\mmm_5$ to the corresponding end state letters of $\mmm$.
Thus, the main $S$-machine $\mmm$ has unique accept (or stop) configuration which we will denote by \index[g]{S@$S$-machine!M@$\mmm$!W@$W_{ac}$, the accept word of $\mmm$} $W_{ac}$.

.

\begin{lemma} \label{resto} (\cite{OS20}, Lemma 4.4) Let the history of
a reduced computation  $\ccc\colon W_0\to\dots\to W_t$
have a subword $\chi(i-1,i)H'\chi(i,i+1)$  (i.e. the S-machine $\bf M$ works as ${\bf M}_3$ with rules from $\bf\Theta_4$) or the subword  $\zeta^{(i-1,i)} H'\zeta^{(i,i+1)}$ (i.e. it works as ${\bf LR}_m$ with rules from $\bf\Theta_2$).
Then the base of the computation $\cal C$ is a reduced word.
and
 all configurations of $\ccc$ are uniquely defined by the history
 $H$ and the base of $\ccc$.
 \end{lemma} $\Box$

We say that the history $H$ of a computation of  $\mmm$ (and the computation itself) is \index[g]{S@$S$-machine!M@$\mmm$!eligible computation of $\mmm$} \index[g]{S@$S$-machine!M@$\mmm$!eligible history of computation of $\mmm$}{\it eligible} if it has no
neighboring mutually inverse letters except possibly for the subwords $\theta(23)\theta(23)^{-1}$.
(The subword $\theta(23)^{-1}\theta(23)$ is not allowed.)
Considering eligible computations instead of just reduced computations is necessary for our interpretation of $\mmm$ in a group.

The history $H$ of an eligible computation of $\mmm$ can be factorized so that every
factor is either a transition rule $\theta(i,i+1)^{\pm 1}$ or a maximal non-empty product of rules of one of the sets $\mtt_1 - \mtt_5$. If, for example, $H=H'H''H'''$,
where $H'$ is a product of rules from $\mtt_2$, $H''$ has only one rule
$\theta(23)$ and $H'''$ is a product of  rules from $\mtt_3$, then we say
that the \index[g]{S@$S$-machine!M@$\mmm$!eligible computation of $\mmm$!step history of a computation of $\mmm$} {\it step history} of the computation is $(2)(23)(3)$.

Thus the step history of a computation is a word in the alphabet $\{(1),(2), (3), (4),(5)$, $(12), (23), (34), (45), (21), (32), (43), (54)\}$, where $(21)$
is used for the rule $\theta(12)^{-1}$ an so on.
For brevity, we can omit some transition symbols, e.g. we may use $(2)(3)$ instead of $(2)(23)(3)$ since the only rule connecting Steps 2 and 3 is $\theta(23)$.

\begin{lemma} \label{212} (\cite{OS20}, Lemma 4.2 (1)) There are no reduced computations $\ccc$ of $\mmm$ with standard base and step history
$(34)(4)(43)$ or $(54)(4)(45)$.
\end{lemma} $\Box$

If the step history of a computation consists of only one letter $(i)$, $i=1,\dots,5$, then we call it a \index[g]{S@$S$-machine!M@$\mmm$!one step computation of $\mmm$} {\em one step computation}.
The computations with step histories $(i)(i,i\pm 1)$, $(i\pm 1, i)(i)$  and $(i\pm 1, i)(i)(i,i\pm 1)$ are also considered
as one step computations. Any eligible one step computation is always reduced by definition.

If the step history of a computation consists of only one letter $(i)$, $i=1,\dots,5$, then we call it a \index[g]{S@$S$-machine!M@$\mmm$!one step computation of $\mmm$} {\em one step computation}.
The computations with step histories $(i)(i,i\pm 1)$, $(i\pm 1, i)(i)$  and $(i\pm 1, i)(i)(i,i\pm 1)$ are also considered
as one step computations.
Any eligible one step computation is always reduced by definition.

By definition,  the rule $\theta(23)$ locks all history sectors of the standard base of $\mmm$
except for the input sector ${\tilde R}_0{\tilde P}_1$ and its mirror copy.
Hence every admissible word in the domain of
$\theta(23)^{-1}$ has the form \index[g]{W@$W(k,k')$ - a word in the domain of $\theta(23)$} $W(k,k')\equiv w_1\alpha^kw_2(\alpha')^{-k'}w_3$,
where $(\alpha')^{-1}$ is the mirror copy of $\alpha$, $k$ and $k'$ are integers, and $w_1, w_2, w_3$ are fixed words in state letters;
$w_1$ starts with $\tt$. Recall that $W_{ac}$ is the accept word of $\mmm$.

\begin{lemma} \label{I6A6} (\cite{OS20}, Lemma 4.6) (1) If the word $\alpha^k$ is accepted by the Turing machine $\mmm_0$, then there
is a reduced computation of $\mmm$, $W(k,k)\to\dots\to W_{ac}$ whose history has no rules of $\mtt_1$ and $\mtt_2$.

(2) If the history of a computation $\ccc\colon W(k,k)\to\dots\to W_{ac}$ of $\mmm$ has no rules of $\mtt_1$ and $\mtt_2$, then the word $\alpha^k$ is accepted by $\mmm_0$.
\end{lemma} $\Box$

A configuration $W$ of $\mmm$ \index[g]{S@$S$-machine!M@$\mmm$!accessible configuration of $\mmm$} is called {\it accessible} if there is a \index[g]{S@$S$-machine!M@$\mmm$!accessible computation of $\mmm$} $W$-\emph{accessible computation}, i.e., either an
accepting computation starting with $W$ or a  computation
$s_1(\mmm)\to\dots\to W$, where \index[g]{S@$S$-machine!M@$\mmm$!s@$s_1(\mmm)$ - the start configuration of $\mmm$} $s_1(\mmm)$ is the start configuration of $\mmm$ (i.e., the configuration where all state letters are start state letters of $\Theta_1$ and the $Y$-projection is empty).

The base of a computation is called {\it revolving} if it starts and ends with the same letter and has no proper subwords with this property. If this base $xvx$ is a reduced word, then it follows from the definition of admissible words that the cyclic order of letters in the word $xv$ is the same as in the standard base. i.e. $xv$  is a cyclic permutation of the standard base.

\begin{lemma}\label{narrow} (\cite{OS20}, Lemmas 4.8 and 4.12) Suppose the base $xvx$ of an eligible computation $\ccc\colon W_0\to\dots\to W_t$ is revolving. Then one of the following statements holds:

(1) we have inequality $||W_j||\le c_4\max(||W_0||, ||W_t||)$, for every $j=0,\dots,t$  or

(2) the base $xvx$ is reduced and if $xv$ is the standard base, then the words $W_0$ and $W_t$ without the last $x$-letters are  accessible words;
the brief history of $\cal C$ contains a subword $(34)(4)(45)$ or a subword $(12)(2)(23)$.
\end{lemma} $\Box$

\begin{rk} \label{les} By Lemma 3.15 \cite{OS20} (by Remark 3.7 \cite{OS20}), a computation with standard base and brief history $(34)(4)(45)$ (resp., with brief history $(12)(2)(23)$) has a subword $\chi(i-1,i)H'\chi(i,i+1)$ (subword $\zeta^{(i-1,i)} H'\zeta^{(i,i+1)}$), as in Lemma \ref{resto}.
\end{rk}

\begin{lemma} \label{xvx}. Suppose $\ccc\colon W_0\to\dots\to W_t$ is an eligible computation, with a base $xvx$. Then
either $(xv)^{\pm 1}$ is a power of a cyclic permutation of the standard base or $|W_j|_Y\le c_4\max(|W_0|_Y, |W_t|_Y)$, for every $j=0,\dots,t$.
\end{lemma}

\proof Note that the base of $\cal C$ has a revolving subword $yv'y$.
Let ${\cal D}\colon V_0\to\dots\to V_t$ be the computation $\cal C$ restricted to this subbase. It has the same history $H$. By Lemma \ref{narrow},
the  base of $\cal D$ is a reduced word and so
$yv'$ is a cyclic permutation of the standard base
or $|V_j|_Y\le c_4\max(|V_0|_Y, |V_t|_Y)$, for every $j=0,\dots,t$.

In the latter case, let us remove the subwords with the base $yv'$, obtaining a computation ${\cal E}: U_0\to\dots\to U_t$  with a shorter base. Arguing by induction, we have either $|U_j|_Y\le c_4\max(|U_0|_Y, |U_t|_Y)$ for every $j=0,\dots,t$, which implies $|W_j|_Y\le c_4\max(|W_0|_Y, |W_t|_Y)$,
or the base of $\cal E$ is a power of a cyclic permutation of the standard base and by Lemma \ref{narrow},
the brief history of $\cal C$ contains a subword $(3)(34)(4)$ or a subword $(12)(2)(23)$. Then by Remark \ref{les}, one can apply Lemma \ref{resto}, and since the computation $\cal D$ has the same history as $\cal E$, the base $yv'y$ must be reduced.
Therefore $yv'$ is a cyclic permutation of the standard base, and so
$xv$ is a power of a cyclic permutation of the standard base.

If $|U_j|_Y\le c_4\max(|U_0|_Y, |U_t|_Y)$ for every $j$, but $yv'$ is a cyclic permutation of the standard base, then the dual
argument
%the word $yv'y$ is reduced, which
implies that the base of $\cal E$ and the base of $\cal C$ are reduced words. Hence $xv$ is a power of a cyclic permutation of the standard base.
\endproof

\section{Group and diagram preliminaries}\label{gdp}

\subsection{The groups}\label{MG}

Every $S$-machine can be simulated by a finitely presented group (see \cite{SBR}, \cite{OS04}, \cite{OS06},  etc.). Here we present the  construction from \cite{OS20}.
To simplify formulas, it is convenient to change the notation. From now on we shall denote by \index[g]{parameters used in the paper!n@$N$ - the length of the standard base of the $S$-machine $\mmm$}$N$ the length of
the standard base of $\mmm$.

Thus the set of state letters is $Q=\sqcup_{i=0}^{N-1}Q_i$ (we set $Q_{N}=Q_0=\{\tt\}$), $Y= \sqcup_{i=1}^{N} Y_i,$ and $\Theta$ is the set of rules of the $S$-machine $\mmm$.

The finite set of generators of the group $M$ \index[g]{M@group $M$}\index[g]{M@group $M$!generators of the group $M$}  consists of {\em $q$-letters},
{\em $Y$-letters} and
{\em $\theta$-letters} defined as follows.

For every letter $q\in Q$ the set of generators of $M$ contains
$L$ copies $q^{(i)}$ of it,  $i=1,\dots, L$, if the letter $q$ occurs in the rules of $\mtt_1$ or $\mtt_2$. \index[g]{parameters used in the paper!l@$L$ - the number of generators $q^{(i)}$ of the group $G$ for each state letter $q$ of $\mmm$, the order of $W_{ac}$ in $G$}(The  number $L$ is one of the parameters from (\ref{const}).)  Otherwise only the letter $q$ is included
in the generating set of $M$.

For every letter $a\in Y$ the set of generators of $M$ contains $a$
and $L$ copies $a^{(i)}$ of it $i=1,..., L$.

For every $\theta\in \Theta^+$ we have $N$ generators $\theta_0,\dots,\theta_N$ in $M$ (here $\theta_{N}\equiv\theta_0$)
if $\theta$ is a rule of $\Theta_3$ (excluding $\theta(23)$) or $\Theta_4$, or $\Theta_5$. For $\theta$ from $\Theta_1$ or $\Theta_2$ (including $\theta(23)$),
we introduce $LN$ generators $\theta_j^{(i)}$, where $j=0,\dots, N$,
$i=1,\dots, L$ and $\theta_N^{(i)}=\theta_0^{(i+1)}$ (the superscripts are taken modulo $L$).

The \index[g]{M@group $M$!relations of the group $M$} relations of the group $M$ correspond to the rules of the $S$-machine $\mmm$ as follows.
For every rule $\theta=[U_0\to V_0,\dots U_{N}\to V_{N}]\in \Theta^+$ of sets $\Theta_1$ or $\Theta_2$, we have

\begin{equation}\label{rel1}
U_j^{(i)}\theta_{j+1}^{(i)}=\theta_j^{(i)} V_j^{(i)},\,\,\,\, \qquad \theta_j^{(i)} a^{(i)}=a^{(i)}\theta_j^{(i)}, \,\,\,\, j=0,...,N,\;\;
i=1,\dots L,
\end{equation}
for all $a\in Y_j(\theta)$, where $U_j^{(i)}$ and $V_j^{(i)}$ are obtained from $U_j$ and $V_j$
by adding the superscript ${(i)}$ to every letter.

For $\theta=\theta(23)$, we introduce relations

\begin{equation}\label{rel11}
U_j^{(i)}\theta_{j+1}^{(i)}=\theta_j^{(i)} V_j,\,\,\,\, \qquad a^{(i)}\theta_j^{(i)}=\theta_j^{(i)} a,
\end{equation}
for all $a\in Y_j(\theta)$, i.e.,the superscripts are erased in the words $U_j^{(i)}$ and in the $Y$-letters after an application of (\ref{rel11}).

For every rule $\theta=[U_0\to V_0,\dots U_{N}\to V_{N}]\in \Theta^+$ from $\Theta_3$ or $\Theta_4$, or $\Theta_5$ and $a\in Y_j(\theta)$, we define

\begin{equation}\label{rel111}
U_j\theta_{j+1}=\theta_j V_j,\,\,\,\, \qquad a\theta_j=\theta_j a
\end{equation}

The first type of relations (\ref{rel1} - \ref{rel111}) will be
called $(\theta,q)$-{\em relations}, the second type - \label{thetaar}
$(\theta,a)$-{\em relations}.

Finally, the required \label{groupG} group $G$ is given by the generators and
relations of the group $M$ and by two more additional
relations, namely the \label{hubr} {\it hub}-relations

\begin{equation}\label{rel3}
W_{st}^{(1)}\dots W_{st}^{(L)}=1\;\; and \;\; (W_{ac})^L=1,
\end{equation}
where the word $W_{st}^{(i)}$ is a copy with superscript $(i)$ of the start word $W_{st}$ (of length $N$) of the $S$-machine $\mmm$
and $W_{ac}$
is  the accept word of $\mmm$.

Note that, as usual, $M$ is a multiple HNN extension  of the free group generated by all $Y$- and $q$-letters, because by Tietze transformations using $(\theta,q)$-relations, all $\theta$-letters, except for one for every rule
$\theta$, can be eliminated.

\subsection{Van Kampen diagrams}\label{md}

Recall that a \index[g]{van Kampen diagram} van Kampen {\it diagram} $\Delta $ over a presentation
$P=\langle A | \rrr\rangle$ (or just over the group $P$)
is a finite oriented connected and simply--connected planar 2--complex endowed with a
\index[g]{van Kampen diagram!labeling function} {\em labeling function} $\Lab \colon E(\Delta )\to A^{\pm 1}$, where $E(\Delta
) $ denotes the set of oriented edges of $\Delta $, such that $\Lab
(e^{-1})\equiv \Lab (e)^{-1}$. Given a \index[g]{van Kampen diagram!cell} {\em cell} (that is a 2-cell) $\Pi $ of $\Delta $,
we denote by $\partial \Pi$ the boundary of $\Pi $; similarly, \index[g]{van Kampen diagram!boundary $\partial(\Delta)$}
$\partial \Delta $ denotes the boundary of $\Delta $. The labels of
$\partial \Pi $ and $\partial \Delta $ are defined up to cyclic
permutations. An additional requirement is that the label of any
cell $\Pi $ of $\Delta $ is equal to (a cyclic permutation of) a
word $R^{\pm 1}$, where $R\in \rrr$. The label and the \index[g]{combinatorial length of a path} combinatorial length $||\bf p||$ of
a path $\bf p$ are defined as for Cayley graphs.

The van Kampen Lemma \cite{LS,book, Sbook} states that a word $W$ over the alphabet $A^{\pm 1}$
represents the identity in the group $P$ if and only
if there exists a diagram $\Delta
$ over $P$ such that
$\Lab (\partial \Delta )\equiv W,$ in particular, the combinatorial perimeter $||\partial\Delta||$ of $\Delta$ equals $||W||.$
(\cite{LS}, Ch. 5, Theorem 1.1; our formulation is closer to
Lemma 11.1 of \cite{book}, see also \cite[Section 5.1]{Sbook}). The word $W$ representing $1$ in $P$ is freely equal
to a product of conjugates of the words from $R^{\pm 1}$.
The minimal number
of factors in such products is called the \index[g]{area of a word} {\em area} of the word $W.$ The \index[g]{van Kampen diagram!area} {\it area}
of a diagram $\Delta$ is the number of cells in it.
The proof of the van Kampen Lemma \cite{book, Sbook} shows that  $\area(W)$ is equal
to the area of a van Kampen diagram having
the smallest number of cells among all van Kampen diagrams with  boundary label $\Lab (\partial \Delta )\equiv W.$

The definition of {\it annular diagram} $\Delta$ over a group $G$ is similar to the definition of van Kampen diagram, but the complement of $\Delta$ in the plane has two
connected components. So $\Delta$ has two boundary components.
By the van Kampen-Schupp lemma (see \cite{LS}, Lemma 5.2 or \cite{book}, Lemma 11.2) there is an  annular diagram $\Delta$ whose boundary components ${\mathbf p}_1$ and ${\mathbf p}_2$ have clockwise labels $W$ and $W'$ if and only if the words $W$ and $W'$ are conjugate in $G$.

We will study diagrams over the group presentations of $M$ and $G$. The edges labeled by state
letters ( = $q$-{\it letters}) will be called
\label{qedge} $q$-{\it edges}, the edges labeled by tape
letters (= $Y$-{\it letters}) will be called \label{Yedge} $Y$-{\it edges}, and the edges labeled by
$\theta$-letters are \label{thedge} $\theta$-{\it edges}.

We denote by $|{\bf p}|_Y$ (by $|{\bf p}|_{\theta}$, by
$|{\bf p}|_q$)
the \label{alength} $Y$-{\it length} (resp., the \label{thlength} $\theta$-{\it length}, the \label{qlength} $q$-length) of a path $\bf p,$ i.e., the number of
$Y$-edges (the number of $\theta$-edges, the number of $q$-edges) in $\bf p.$

The cells corresponding
to relations (\ref{rel3})  are called \index[g]{hub}{\it hubs}, the cells corresponding
to $(\theta,q)$-relations are called \index[g]{t@$(\theta,q)$-cell}$(\theta,q)$-{\it cells} (in particular, there are $(\theta,\tt)$-{\it cells}),
and the cells are called \index[g]{t@$(\theta,a)$-cell}$(\theta,a)$-{\it cells} if they correspond to $(\theta,a)$-relations. A $\theta$-{\it cell} is either $(\theta,q)$- or $(\theta,a)$-cell.

A diagram is \index[g]{van Kampen diagram!reduced}{\em reduced}, if
it does not contain two cells (= closed $2$-cells) that have a
common edge $\bf e$ such that the boundary labels of these two cells are
equal if one reads them starting with $\bf e$
(if such pairs of cells exist, they can be removed to obtain a  diagram of smaller area and with the same boundary label(s)).

To study  diagrams
over the group $G$ we shall use their simpler subdiagrams such as bands.
 Here we repeat one more necessary definition from \cite{OS20}.

\index[g]{band}\begin{df} Let $\mathcal Z$ be a subset of the set of letters in the set of generators of the group $M$. A
$\mathcal Z$-band $\bb$ is a sequence of cells $\pi_1,...,\pi_n$ in a reduced \vk
diagram $\Delta$ such that

\begin{itemize}
\item Every two consecutive cells $\pi_i$ and $\pi_{i+1}$ in this
sequence have a common boundary edge ${\mathbf e}_i$ labeled by a letter from ${\mathcal Z}^{\pm 1}$ (Fig. \ref{bt}).
\item Each cell $\pi_i$, $i=1,...,n$ has exactly two $\mathcal Z$-edges in the boundary $\partial \pi_i$,
${\mathbf e}_{i-1}^{-1}$ and ${\mathbf e}_i$ (i.e.,edges labeled by a letter from ${\mathcal Z}^{\pm 1}$) with the requirement that either
both $\Lab(e_{i-1})$ and $\Lab(e_i)$ are positive letters or both
are negative ones.

\item If $n=0$, then $\bb$ is just a $\mathcal Z$-edge.
\end{itemize}
\end{df}

The counter-clockwise boundary of the subdiagram formed by the
cells $\pi_1,...,\pi_n$ of $\bb$ has the factorization ${\mathbf e}\iv {\mathbf q}_1{\mathbf f} {\mathbf q}_2\iv$
where ${\mathbf e}={\mathbf e}_0$ is a $\mathcal Z$-edge of $\pi_1$ and ${\mathbf f}={\mathbf e}_n$ is an $\mathcal Z$-edge of
$\pi_n$. We call ${\mathbf q}_1$ the {\em bottom} of $\bb$ and ${\mathbf q}_2$ the
{\em top} of $\bb$, denoted $\bott(\bb)$ and $\topp(\bb)$. If the path
${\mathbf e}\iv {\mathbf q}_1{\mathbf f}$ or the path ${\mathbf f} {\mathbf q}_2\iv{\mathbf e}\iv$ is the subpath of the boundary $\partial\Delta$, then $\bb$ is called a {\it rim} band.
Top/bottom paths and their inverses are also called the {\em
sides} of the band. \index[g]{band!sides} The $\mathcal Z$-edges ${\mathbf e}$
and ${\mathbf f}$ are called the {\em start} \index[g]{band!start and end edges} and {\em end} edges of the
band. If $n\ge 1$ but ${\mathbf e}={\mathbf f},$ then the $\mathcal Z$-band is called a  $\mathcal Z$-{\it annulus} \index[g]{band!annulus}.

We consider \index[g]{band!q@$q$-band} $q$-{\it bands}, where $\mathcal Z$ is one of the sets $Q_i$ of state letters
for the $S$-machine $\mmm$, \index[g]{band!t@$\theta$-band}
$\theta$-{\it bands} for every $\theta\in\Theta$, and $Y$-{\it bands}\index[g]{band!a@$Y$-band}, where
${\mathcal Z}=\{a,a^{(1)},\dots, a^{(L)}\}\subseteq Y$.
 The convention is that $Y$-bands do not
contain $(\theta,q)$-cells, and so they consist of $(\theta,a)$-cells  only.

\begin{rk} \label{tb} To construct the top (or bottom) path of a band $\mathcal B$, at the beginning
one can just form a product ${\mathbf x}_1\dots {\mathbf x}_n$ of the top paths ${\mathbf x}_i$-s of the cells $\pi_1,\dots,\pi_n$ (where each $\pi_i$ is a $\mathcal Z$-bands of length $1$).
No  $\theta$-letter is being canceled in the word
$W\equiv\Lab({\mathbf x}_1)\dots\Lab({\mathbf x}_n)$ if $\mathcal B$ is  a $q$- or $Y$-band since
otherwise two neighbor cells of the band would make the diagram non-reduced.
The {\it trimmed} top/bottom label of $\mathcal B$
are the maximal subwords of the top/bottom labels starting and ending
with $q$-letters. \index[g]{band!top path $\topp(\bb)$} \index[g]{band!bottom path $\bott(\bb)$}

However a few cancellations of $Y$-letters  are possible in $W.$ (This can happen if one of $\pi_i, \pi_{i+1}$
is a $(\theta,q)$-cell and another one is a $(\theta,a)$-cell.) We will always assume
that the top/bottom label of a $\theta$-band is a reduced form of the word $W$.
This property  is easy to achieve: by folding edges
with the same labels having the same initial vertex, one can make
the boundary label of a subdiagram in a \vk diagram reduced (e.g., see \cite{book} or
\cite{SBR}).
\end{rk}

\begin{rk} \label{pm1} Since $\theta_N^{(i)}=\theta_0^{(i+1)}$, it follows from (\ref{rel1}) that the superscripts in the $q$-letters of the same $(\theta,q)$-relation
are different if $\theta\in {\bf\Theta}_1\cup {\bf\Theta}_2\cup\{\theta(23)^{\pm 1}\}$
and this relation is a $(\theta,\tt)$-relation. Therefore only the corresponding
cells of a $\theta$-band have different superscripts of the labels of $\theta$-edges,
and this difference is $\pm 1$ modulo $L$.
\end{rk}

We shall call a $\mathcal Z$-band \index[g]{band!maximal band}{\em maximal} if it is not contained in
any other $\mathcal Z$-band.
Counting the number of maximal $\mathcal Z$-bands
 in a diagram we will not distinguish the bands with boundaries
 ${\mathbf e}\iv {\mathbf q}_1{\mathbf f} {\mathbf q}_2\iv$ and ${\mathbf f} {\mathbf q}_2\iv {\mathbf e}\iv {\mathbf q}_1,$ and
 so every $\mathcal Z$-edge belongs to a unique maximal $\mathcal Z$-band.

We say that a ${\mathcal Z}_1$-band and a ${\mathcal Z}_2$-band \index[g]{band!crossing bands}{\em cross} if
they have a common cell and ${\mathcal Z}_1\cap {\mathcal Z}_2=\emptyset.$

Sometimes we specify the types of bands as follows.
A $q$-band corresponding to one
letter $Q$ of the base is called a \index[g]{Q@$Q$-band}$Q$-band. For example, we will consider \index[g]{band!t@$\tt$-band}$\tt$-{\it band} corresponding to the part $\{\tt\}$.

\begin{lemma}\label{NoAnnul}(\cite{OS20}, Lemma 5.6)
A reduced van Kampen diagram $\Delta$ over $M$ has no
$q$-annuli, no $\theta$-annuli, and no  $Y$-annuli.
Every $\theta$-band of $\Delta$ shares at most one cell with any
$q$-band and with any $Y$-band.
\end{lemma} $\Box$

If $W\equiv x_1...x_n$ is a word in an alphabet $X$, $X'$ is another
alphabet, and $\phi\colon X\to X'\cup\{1\}$ (where $1$ is the empty
word) is a map, then $\phi(W)\equiv\phi(x_1)...\phi(x_n)$ is called the
\label{projectw}{\em projection} of $W$ onto $X'$. We shall consider the
projections of words in the generators of $M$ onto
$\Theta$ (all $\theta$-letters map to the
corresponding element of $\Theta$,
all other letters map to $1$), and the projection onto the
alphabet $\{Q_0\sqcup \dots \sqcup Q_{N-1}\}$ (every
$q$-letter maps to the corresponding $Q_i$, all other
letters map to $1$).

\begin{df}\label{dfsides}
{\rm  The projection of the label
of a side of a $q$-band onto the alphabet $\Theta$ is
called the {\em history} of the band. \index[g]{band!history of a $q$-band} The step history of this projection
is the \index[g]{band!step history of a $q$-band}{\it step history} of the $q$-band. The projection of the label
of a side of a $\theta$-band onto the alphabet $\{Q_0,...,Q_{N-1}\}$
is called the {\em base} of the band, i.e., the base of a $\theta$-band
is equal to the base of the label of its top or bottom}\index[g]{band!base of a $\theta$-band}
\end{df}
As in the case of words, we will  use representatives of
$Q_j$-s in base words.

If $W$ is a word in the generators of $M$, then by $W^{\emptyset}$
we denote the projection of this word onto the alphabet of the
$S$-machine $\mmm$, we obtain this projection after deleting all
superscripts in the letters of $W$. In particular, $W^{\emptyset}\equiv W$, if there are no superscripts in the letters of $W$.

We call a word $W$ in $q$-generators and $Y$-generators {\it permissible} \index[g]{permissible word} if the
word $W^{\emptyset}$ is admissible, and the letters of any 2-letter subword of $W$
have equal superscripts (if any), except for the subwords $(q\tt)^{\pm 1}$,
where the letter  $q$ has some superscript $(i)$ and $q^{\emptyset} \in Q_{N-1}$; in this case the superscript of the letter $\tt$ must be $(i+1)$
(modulo $L$).

\begin{rk} \label{perad} It follows from the definition that if $V$ is $\theta$-admissible for a rule $\theta$
of $\{\theta(23)\iv \}\cup  \mtt_3\cup \{\theta(34)\}\cup \mtt_4\cup\{\theta(45)\}\cup \mtt_5$, then there is exactly one permissible
word $W$ such that $W^{\emptyset}\equiv V$, namely, $W\equiv V$.
If $\theta$ is a rule of $\mtt_1\cup\{\theta(12)\}\cup \mtt_2\cup\{\theta(23)\}$, then the permissible word $W$
with property $W^{\emptyset}\equiv V$ exists and it is uniquely defined if one choose arbitrary
superscript for the first letter (or for any particular letter) of $W$.
\end{rk}

\begin{lemma} \label{perm} (\cite{OS20}, Lemma 5.9) (1) The trimmed bottom and top labels $W_1$ and $W_2$ of any reduced $\theta$-band $\mathcal T$ containing at least one $(\theta,q)-cell$ are permissible
and $W_2^{\emptyset}\equiv W_1^{\emptyset}\cdot\theta$.

(2) If $W$ is a $\theta$-admissible word, then for a
permissible word $W_1$ such that $W_1^{\emptyset}\equiv W$ (given by Remark \ref{perad})
one can construct a reduced $\theta$-band with the trimmed bottom
label $W_1$ and the trimmed top label $W_2$, where $W_2^{\emptyset}\equiv
W_1^{\emptyset}\cdot\theta$.
\end{lemma} $\Box$

% This is a LaTeX picture output by TeXCAD.
% File name: [bandtr.pic].
% Version of TeXCAD: 4.3
% Reference / build: 30-Jun-2012 (rev. 105)
% For new versions, check: http://texcad.sf.net/
% Options on the following lines.
%\grade{\on}
%\emlines{\off}
%\epic{\off}
%\beziermacro{\on}
%\reduce{\on}
%\snapping{\off}
%\pvinsert{% Your \input, \def, etc. here}
%\quality{8.000}
%\graddiff{0.005}
%\snapasp{1}
%\zoom{4.0000}
\unitlength 1mm % = 2.845pt
\linethickness{0.4pt}
\ifx\plotpoint\undefined\newsavebox{\plotpoint}\fi % GNUPLOT compatibility
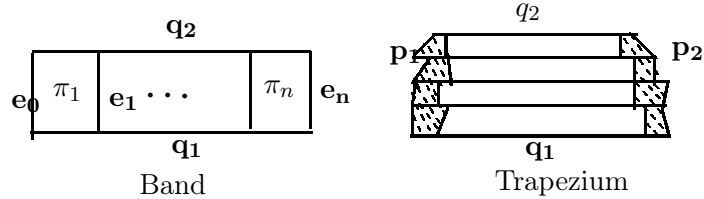
\begin{picture}(98.5,35.75)(0,85.5)
\thicklines
\put(13.5,115.75){\line(0,-1){11.5}}
\put(13.5,116){\line(1,0){37}}
\put(50.5,116){\line(0,-1){10.25}}
\put(13.25,105.25){\line(1,0){37.5}}
\put(22.25,115.75){\line(0,-1){10.25}}
\put(42.5,116){\line(0,-1){10.25}}
\put(16.25,110){$\pi_1$}
\put(44.25,110.5){$\pi_n$}
\put(10.75,108.5){$\bf e_0$}
\put(23.75,108.75){$\bf e_1$}
\put(51.75,109.75){$\bf e_n$}
\put(29.25,110.25){\circle{.5}}
\put(33.75,110){\circle{.5}}
\put(31.5,110.25){\circle{.5}}
\put(32,102.5){$\bf q_1$}
\put(31.25,118.5){$\bf q_2$}
\put(27.75,97){Band}
\put(63.75,104.75){\line(1,0){34.5}}
\put(64.25,108.75){\line(1,0){33.5}}
\put(64.25,112){\line(1,0){34}}
\put(64,115.25){\line(1,0){32.5}}
\put(66.75,118.25){\line(1,0){27.25}}
%\emline(66.75,118)(64.25,115.5)
\multiput(66.75,118)(-.03333333,-.03333333){75}{\line(0,-1){.03333333}}
%\end
\put(68.5,118){\line(0,-1){2.75}}
%\emline(66.25,115)(64,111.75)
\multiput(66.25,115)(-.03358209,-.04850746){67}{\line(0,-1){.04850746}}
%\end
%\emline(64,111.75)(64.25,112)
\multiput(64,111.75)(.03125,.03125){8}{\line(0,1){.03125}}
%\end
%\emline(68.75,115.75)(69.25,111.5)
\multiput(68.75,115.75)(.0333333,-.2833333){15}{\line(0,-1){.2833333}}
%\end
\put(69.25,111.5){\line(0,1){0}}
%\emline(69.25,111.5)(69,112)
\multiput(69.25,111.5)(-.03125,.0625){8}{\line(0,1){.0625}}
%\end
\put(67.5,111.75){\line(0,-1){3}}
%\emline(64.75,111.25)(64.5,111.75)
\multiput(64.75,111.25)(-.03125,.0625){8}{\line(0,1){.0625}}
%\end
%\emline(64.5,112.25)(64,108.75)
\multiput(64.5,112.25)(-.0333333,-.2333333){15}{\line(0,-1){.2333333}}
%\end
\put(64,108.5){\line(0,-1){3.5}}
%\emline(68.75,108.5)(67.25,104.75)
\multiput(68.75,108.5)(-.03333333,-.08333333){45}{\line(0,-1){.08333333}}
%\end
%\emline(93.5,118)(96.5,115)
\multiput(93.5,118)(.03370787,-.03370787){89}{\line(0,-1){.03370787}}
%\end
\put(96.25,115.5){\line(0,-1){4}}
%\emline(98,112)(97.25,108.25)
\multiput(98,112)(-.0326087,-.1630435){23}{\line(0,-1){.1630435}}
%\end
%\emline(97.25,108.25)(98.25,104.5)
\multiput(97.25,108.25)(.0333333,-.125){30}{\line(0,-1){.125}}
%\end
\put(91.75,118){\line(0,-1){3}}
\put(93.25,108.25){\line(0,1){0}}
\put(93.5,115.25){\line(0,-1){7}}
\put(95,108.5){\line(0,-1){3.5}}
\thinlines
%\dashline{1}(65.25,116.5)(68,111.5)
\multiput(65.18,116.43)(.032738,-.059524){12}{\line(0,-1){.059524}}
\multiput(65.965,115.001)(.032738,-.059524){12}{\line(0,-1){.059524}}
\multiput(66.751,113.573)(.032738,-.059524){12}{\line(0,-1){.059524}}
\multiput(67.537,112.144)(.032738,-.059524){12}{\line(0,-1){.059524}}
%\end
%\dashline{1}(66.25,116.75)(68.75,112)
\multiput(66.18,116.68)(.0320513,-.0608974){13}{\line(0,-1){.0608974}}
\multiput(67.013,115.096)(.0320513,-.0608974){13}{\line(0,-1){.0608974}}
\multiput(67.846,113.513)(.0320513,-.0608974){13}{\line(0,-1){.0608974}}
%\end
%\dashline{1}(67,117.75)(68.5,114.25)
\multiput(66.93,117.68)(.033333,-.077778){9}{\line(0,-1){.077778}}
\multiput(67.53,116.28)(.033333,-.077778){9}{\line(0,-1){.077778}}
\multiput(68.13,114.88)(.033333,-.077778){9}{\line(0,-1){.077778}}
%\end
%\dashline{1}(65.75,114)(67.25,110.25)
\multiput(65.68,113.93)(.033333,-.083333){9}{\line(0,-1){.083333}}
\multiput(66.28,112.43)(.033333,-.083333){9}{\line(0,-1){.083333}}
\multiput(66.88,110.93)(.033333,-.083333){9}{\line(0,-1){.083333}}
%\end
%\dashline{1}(65.25,113)(67.5,108.75)
\multiput(65.18,112.93)(.03125,-.059028){12}{\line(0,-1){.059028}}
\multiput(65.93,111.513)(.03125,-.059028){12}{\line(0,-1){.059028}}
\multiput(66.68,110.096)(.03125,-.059028){12}{\line(0,-1){.059028}}
%\end
%\dashline{1}(64.75,111.25)(67.75,106.5)
\multiput(64.68,111.18)(.032967,-.0521978){13}{\line(0,-1){.0521978}}
\multiput(65.537,109.823)(.032967,-.0521978){13}{\line(0,-1){.0521978}}
\multiput(66.394,108.465)(.032967,-.0521978){13}{\line(0,-1){.0521978}}
\multiput(67.251,107.108)(.032967,-.0521978){13}{\line(0,-1){.0521978}}
%\end
%\dashline{1}(67.5,109)(68.5,107.75)
\multiput(67.43,108.93)(.033333,-.041667){10}{\line(0,-1){.041667}}
\multiput(68.096,108.096)(.033333,-.041667){10}{\line(0,-1){.041667}}
%\end
%\dashline{1}(64,110.5)(67.5,105.75)
\multiput(63.93,110.43)(.0333333,-.0452381){15}{\line(0,-1){.0452381}}
\multiput(64.93,109.073)(.0333333,-.0452381){15}{\line(0,-1){.0452381}}
\multiput(65.93,107.715)(.0333333,-.0452381){15}{\line(0,-1){.0452381}}
\multiput(66.93,106.358)(.0333333,-.0452381){15}{\line(0,-1){.0452381}}
%\end
%\dashline{1}(64.5,109)(66.5,105.5)
\multiput(64.43,108.93)(.033333,-.058333){12}{\line(0,-1){.058333}}
\multiput(65.23,107.53)(.033333,-.058333){12}{\line(0,-1){.058333}}
\multiput(66.03,106.13)(.033333,-.058333){12}{\line(0,-1){.058333}}
%\end
%\dashline{1}(64.25,107.5)(65.5,104.75)
\multiput(64.18,107.43)(.03125,-.06875){10}{\line(0,-1){.06875}}
\multiput(64.805,106.055)(.03125,-.06875){10}{\line(0,-1){.06875}}
%\end
%\dashline{1}(93,118.5)(96.75,113)
\multiput(92.93,118.43)(.0334821,-.0491071){16}{\line(0,-1){.0491071}}
\multiput(94.001,116.858)(.0334821,-.0491071){16}{\line(0,-1){.0491071}}
\multiput(95.073,115.287)(.0334821,-.0491071){16}{\line(0,-1){.0491071}}
\multiput(96.144,113.715)(.0334821,-.0491071){16}{\line(0,-1){.0491071}}
%\end
%\dashline{1}(92.25,117.75)(96.25,112.25)
\multiput(92.18,117.68)(.0333333,-.0458333){15}{\line(0,-1){.0458333}}
\multiput(93.18,116.305)(.0333333,-.0458333){15}{\line(0,-1){.0458333}}
\multiput(94.18,114.93)(.0333333,-.0458333){15}{\line(0,-1){.0458333}}
\multiput(95.18,113.555)(.0333333,-.0458333){15}{\line(0,-1){.0458333}}
%\end
%\dashline{1}(92,116.25)(92.5,115.5)
\multiput(91.93,116.18)(.03125,-.046875){8}{\line(0,-1){.046875}}
%\end
%\dashline{1}(93.5,113.75)(97.5,109.25)
\multiput(93.43,113.68)(.0333333,-.0375){15}{\line(0,-1){.0375}}
\multiput(94.43,112.555)(.0333333,-.0375){15}{\line(0,-1){.0375}}
\multiput(95.43,111.43)(.0333333,-.0375){15}{\line(0,-1){.0375}}
\multiput(96.43,110.305)(.0333333,-.0375){15}{\line(0,-1){.0375}}
%\end
%\dashline{1}(96.5,111.5)(97.5,110.25)
\multiput(96.43,111.43)(.033333,-.041667){10}{\line(0,-1){.041667}}
\multiput(97.096,110.596)(.033333,-.041667){10}{\line(0,-1){.041667}}
%\end
%\dashline{1}(93.75,112)(97.75,107.25)
\multiput(93.68,111.93)(.0333333,-.0395833){15}{\line(0,-1){.0395833}}
\multiput(94.68,110.742)(.0333333,-.0395833){15}{\line(0,-1){.0395833}}
\multiput(95.68,109.555)(.0333333,-.0395833){15}{\line(0,-1){.0395833}}
\multiput(96.68,108.367)(.0333333,-.0395833){15}{\line(0,-1){.0395833}}
%\end
%\dashline{1}(93.75,110)(98,105.25)
\multiput(93.68,109.93)(.0332031,-.0371094){16}{\line(0,-1){.0371094}}
\multiput(94.742,108.742)(.0332031,-.0371094){16}{\line(0,-1){.0371094}}
\multiput(95.805,107.555)(.0332031,-.0371094){16}{\line(0,-1){.0371094}}
\multiput(96.867,106.367)(.0332031,-.0371094){16}{\line(0,-1){.0371094}}
%\end
%\dashline{1}(95.25,106.5)(96.75,105)
\multiput(95.18,106.43)(.0333333,-.0333333){15}{\line(0,-1){.0333333}}
\multiput(96.18,105.43)(.0333333,-.0333333){15}{\line(0,-1){.0333333}}
%\end
\put(79,102.25){$\bf q_1$}
\put(61,115){$\bf p_1$}
\put(77.75,120.75){$q_2$}
\put(98.5,115.5){$\bf p_2$}
\put(75,97.5){Trapezium}
\end{picture}

%\begin{figure}
%\begin{center}
%\includegraphics[width=1.0\textwidth]{Pic5.jpg}
%\end{center}
%\caption{Band and Trapezium}\label{bt}
%\end{figure}

\begin{df}\label{dftrap}
{\rm Let $\Delta$ be a reduced van Kampen  diagram over $M$,
which has  boundary path of the form ${\mathbf p}_1\iv {\mathbf q}_1{\mathbf p}_2{\mathbf q}_2\iv,$ where
${\mathbf p}_1$ and ${\mathbf p}_2$ are sides of $q$-bands, and
${\mathbf q}_1$, ${\mathbf q}_2$ are maximal parts of the sides of
$\theta$-bands such that $\Lab({\mathbf q}_1)$, $\Lab({\mathbf q}_2)$ start and end
with $q$-letters.

%\begin{figure}[ht]
%\begin{center}
%\includegraphics[width=1.0\textwidth]{Pic5.jpg}
%\end{center}
%\caption{Band and Trapezium}\label{bt}
%\end{figure}

Then $\Delta$ is called a \index[g]{trapezium}{\em trapezium}. The path ${\mathbf q}_1$ is
called the \index[g]{trapezium!bottom}{\em bottom}, the path ${\mathbf q}_2$ is called the \index[g]{trapezium!top}{\em top} of
the trapezium, the paths ${\mathbf p}_1$ and ${\mathbf p}_2$ are called the \index[g]{trapezium!left and right sides}{\em left
and right sides} of the trapezium. The history \index[g]{trapezium!step history} \index[g]{trapezium!history}(step history) of the $q$-band
whose side is ${\mathbf p}_2$ is called the {\em history} (resp., step history) of the trapezium;
the length of the history is called the \index[g]{trapezium!height}{\em height}  of the
trapezium. The base of $\Lab ({\mathbf q}_1)$ is called the \index[g]{trapezium!base}{\em base} of the
trapezium.} (Fig. \ref{bt})
\end{df}

\begin{rk} Notice that the top (bottom) side of a
$\theta$-band $\ttt$ does not necessarily coincides with the top
(bottom) side ${\mathbf q}_2$ (side ${\mathbf q}_1$) of the corresponding trapezium of height $1$, and ${\mathbf q}_2$
(${\mathbf q}_1$) is
obtained from $\topp(\ttt)$ (resp. $\bott(\ttt)$) by trimming the
first and the last $Y$-edges
if these paths start and/or end with $Y$-edges.
\end{rk}

By Lemma \ref{NoAnnul}, any trapezium $\Delta$ of height $h\ge 1$
can be decomposed into $\theta$-bands $\ttt_1,...,\ttt_h$ connecting
the left and the right sides of the trapezium.

\begin{lemma}\label{simul} (\cite{OS20}, Lemma 5.12)  (1) Let $\Delta$ be a trapezium
with history $H\equiv\theta(1)\dots\theta(d)$ ($d\ge 1$).
Assume that $\Delta$
has consecutive maximal $\theta$-bands  ${\mathcal T}_1,\dots
{\mathcal T}_d$, and the words
$U_j$
and $V_j$
are the  trimmed bottom and the
trimmed top labels of ${\mathcal T}_j,$ ($j=1,\dots,d$).
Then $H$ is an eligible word, $U_j$, $V_j$ are
permissible
words,
$$V_1^{\emptyset}\equiv U_1^{\emptyset}\cdot \theta(1),\;\; U_2\equiv V_1,\;\;\dots,\;\; U_d \equiv V_{d-1},\;\; V_d^{\emptyset}\equiv U_d^{\emptyset}\cdot \theta(d)$$

Furthermore, if the first and the last $q$-letters of the word $U_j$ or of the word $ V_j$ have some superscripts $(i)$ and $(i')$, then  $i'-i$ (modulo $L$) does not depend on the choice of $U_j$ or $V_j$.

(2) For every eligible computation $U\to\dots\to U\cdot H \equiv V$ of $\mmm$
with $||H||=d\ge 1$
there exists
a trapezium $\Delta$ with bottom label $U_1$ (given by Remark \ref{perad}) such that $U_1^{\emptyset}\equiv U$, top label $V_d$ such that $V_d^{\emptyset}\equiv V$, and with history $H.$
\end{lemma} $\Box$

Using Lemma \ref{simul}, one can immediately derive properties of trapezia from the properties of computations obtained earlier.

If $H'\equiv \theta(i)\dots\theta(j)$ is a subword of the history $H$
from Lemma \ref{simul} (1), then the bands ${\mathcal T}_i,\dots, {\mathcal T}_j$ form a subtrapezium $\Delta'$ of the trapezium $\Delta$ with the same base. A subword of the base of $\Delta$ also defines a subtrapezium with the same history.

\begin{df}\label{d:standard}
We say that a trapezium $\Delta$ is \index[g]{trapezium!standard} {\it standard} if the base of $\Delta$
  is the standard base $\bf B$, and the
  history of $\Delta$ (or the inverse one) contains one of the words (a) $\chi(i-1,i)H'\chi(i,i+1)$ (i.e.,the $S$-machine
works as $\mtt_4$) or (b) $\zeta^{i-1,i} H'\zeta^{i,i+1}$ (i.e.,it works as $\mtt_2$).
\end{df}

\begin{df} \label{dw} A permissible word $V$ is called a \index[g]{disk word} {\it disk word} if $V^{\emptyset}\equiv W^L$ for
some accessible word $W$. (In particular, hub words are disk words.)
\end{df}

\begin{lemma}\label{trivial} (\cite{OS20}, Lemma 7.2). Every disk word $V$ is equal to $1$ in the group $G$.
\end{lemma} $\Box$

\begin{lemma}\label{dr} (\cite{OS20}, Remark 7.3) For a disk word $V$, there
is a reduced {\it disk diagram} $\Delta$ over the presentation (\ref{rel1}) - (\ref{rel3}) with boundary label  $V$ built of one hub and $L$ trapezia corresponding to an accessible computation
for the word $W$, where $V^{\emptyset}=W^L$.
\end{lemma} $\Box$

We will increase the set of relations of $G$ by adding
the (infinite) set of {\it disk relations} $V=1$ , one for
every disk word $V$. So we will consider diagrams
with \label{diskr} {\it disks}, where every disk cell (or just {\it disk}) is labeled by
such a word  $V$. (In particular, a hub is a disk.)

\begin{df}\label{minimald}
We will call a reduced van Kampen or annular diagram $\Delta$ {\it minimal} if
\index[g]{minimal diagram over $G$}

(1) the number of disks is minimal for all diagrams with the same boundary label(s) as $\Delta$ and

(2) $\Delta$ has minimal number of $(\theta,\tt)$-cells among
the diagrams with the same boundary label(s) and with minimal number of  disks.

Clearly, a subdiagram of a minimal diagram is minimal itself.

\end{df}

The following is explained in \cite{OS20}, Subsection 7.1.2.

\begin{lemma}\label{2dis} If two disks of a van Kampen diagram $\Delta$
are connected by at least two $\tt$-bands, then there is a diagram $\Delta'$ with the same boundary label and fewer disks in it. In particular, two disks cannot be connected by two $\tt$-bands in a minimal van Kampen diagram or by three $\tt$-bands in a minimal annular diagram.
\end{lemma} $\Box$

Lemma \ref{2dis} implies the following properties. (Part (1) is Lemma 7.5\cite{OS20}, the proof of part (2) is similar.)

\begin{lemma} \label{extdisk} (1) If a van Kampen diagram  contains at least one disk and has no pairs of disks connected
by at least two $\tt$-bands,
then there is a disk $\Pi$ in $\Delta$ such that $L-3$ consecutive maximal $\tt$-bands ${\mathcal B}_1,\dots
{\mathcal B}_{L-3} $ start on $\partial\Pi$ , end on the boundary $\partial\Delta$, and for any $i\in [1,L-4]$,
there are no disks in the subdiagram $\Gamma_i$ bounded by ${\mathcal B}_i$, ${\mathcal B}_{i+1},$ $\partial\Pi,$ and $\partial\Delta$. (Fig. \ref{extd})
%\begin{figure}
%\begin{center}
%\includegraphics[width=0.7\textwidth]{Pic11.jpg}
%\end{center}
%\caption{Lemma \ref{extdisk}}\label{extd}
%\end{figure}
% This is a LaTeX picture output by TeXCAD.
% File name: [extdis.pic].
% Version of TeXCAD: 4.3
% Reference / build: 30-Jun-2012 (rev. 105)
% For new versions, check: http://texcad.sf.net/
% Options on the following lines.
%\grade{\on}
%\emlines{\off}
%\epic{\off}
%\beziermacro{\on}
%\reduce{\on}
%\snapping{\off}
%\pvinsert{% Your \input, \def, etc. here}
%\quality{8.000}
%\graddiff{0.005}
%\snapasp{1}
%\zoom{4.0000}
\unitlength 1mm % = 2.845pt
\linethickness{0.4pt}
\ifx\plotpoint\undefined\newsavebox{\plotpoint}\fi % GNUPLOT compatibility
\begin{picture}(91,69.5)(25,60)
\thicklines
%\emline(30,99.25)(35.25,109.25)
\multiput(30,99.25)(.033653846,.064102564){156}{\line(0,1){.064102564}}
%\end
%\emline(35.25,109.25)(51.25,116)
\multiput(35.25,109.25)(.07960199,.03358209){201}{\line(1,0){.07960199}}
%\end
\put(51.25,116){\line(1,0){39.75}}
%\emline(30.25,99.5)(33.25,93)
\multiput(30.25,99.5)(.03370787,-.07303371){89}{\line(0,-1){.07303371}}
%\end
%\emline(33.25,93)(45.75,82.25)
\multiput(33.25,93)(.039184953,-.0336990596){319}{\line(1,0){.039184953}}
%\end
%\emline(45.75,82.25)(89.25,92.75)
\multiput(45.75,82.25)(.1394230769,.0336538462){312}{\line(1,0){.1394230769}}
%\end
%\circle*(45.75,98.25){5.385}
\put(43.7898,96.2898){\rule{3.9204\unitlength}{3.9204\unitlength}}
\multiput(44.6634,100.0977)(0,-4.3916){2}{\rule{2.1733\unitlength}{.6961\unitlength}}
\multiput(45.1685,100.6814)(0,-5.1285){2}{\rule{1.1631\unitlength}{.2657\unitlength}}
\multiput(46.2191,100.6814)(-1.3069,0){2}{\multiput(0,0)(0,-5.0643){2}{\rule{.3688\unitlength}{.2015\unitlength}}}
\multiput(46.7242,100.0977)(-2.5263,0){2}{\multiput(0,0)(0,-4.1427){2}{\rule{.578\unitlength}{.4473\unitlength}}}
\multiput(46.7242,100.4326)(-2.2997,0){2}{\multiput(0,0)(0,-4.6134){2}{\rule{.3513\unitlength}{.2483\unitlength}}}
\multiput(46.7242,100.5684)(-2.1816,0){2}{\multiput(0,0)(0,-4.8087){2}{\rule{.2333\unitlength}{.1719\unitlength}}}
\multiput(46.963,100.4326)(-2.6535,0){2}{\multiput(0,0)(0,-4.5483){2}{\rule{.2275\unitlength}{.1832\unitlength}}}
\multiput(47.1897,100.0977)(-3.2041,0){2}{\multiput(0,0)(0,-3.9853){2}{\rule{.3247\unitlength}{.2899\unitlength}}}
\multiput(47.1897,100.2752)(-3.0999,0){2}{\multiput(0,0)(0,-4.2441){2}{\rule{.2205\unitlength}{.1938\unitlength}}}
\multiput(47.4019,100.0977)(-3.5164,0){2}{\multiput(0,0)(0,-3.899){2}{\rule{.2125\unitlength}{.2036\unitlength}}}
\multiput(47.5977,97.1634)(-4.3916,0){2}{\rule{.6961\unitlength}{2.1733\unitlength}}
\multiput(47.5977,99.2242)(-4.1427,0){2}{\multiput(0,0)(0,-2.5263){2}{\rule{.4473\unitlength}{.578\unitlength}}}
\multiput(47.5977,99.6897)(-3.9853,0){2}{\multiput(0,0)(0,-3.2041){2}{\rule{.2899\unitlength}{.3247\unitlength}}}
\multiput(47.5977,99.9019)(-3.899,0){2}{\multiput(0,0)(0,-3.5164){2}{\rule{.2036\unitlength}{.2125\unitlength}}}
\multiput(47.7752,99.6897)(-4.2441,0){2}{\multiput(0,0)(0,-3.0999){2}{\rule{.1938\unitlength}{.2205\unitlength}}}
\multiput(47.9326,99.2242)(-4.6134,0){2}{\multiput(0,0)(0,-2.2997){2}{\rule{.2483\unitlength}{.3513\unitlength}}}
\multiput(47.9326,99.463)(-4.5483,0){2}{\multiput(0,0)(0,-2.6535){2}{\rule{.1832\unitlength}{.2275\unitlength}}}
\multiput(48.0684,99.2242)(-4.8087,0){2}{\multiput(0,0)(0,-2.1816){2}{\rule{.1719\unitlength}{.2333\unitlength}}}
\multiput(48.1814,97.6685)(-5.1285,0){2}{\rule{.2657\unitlength}{1.1631\unitlength}}
\multiput(48.1814,98.7191)(-5.0643,0){2}{\multiput(0,0)(0,-1.3069){2}{\rule{.2015\unitlength}{.3688\unitlength}}}
\put(48.443,98.25){\line(0,1){.3329}}
\put(48.422,98.583){\line(0,1){.1648}}
\put(48.396,98.748){\line(0,1){.1629}}
\put(48.36,98.911){\line(0,1){.1604}}
\put(48.314,99.071){\line(0,1){.1572}}
\put(48.259,99.228){\line(0,1){.1535}}
\put(48.193,99.382){\line(0,1){.1491}}
\put(48.118,99.531){\line(0,1){.1442}}
\put(48.035,99.675){\line(0,1){.1388}}
\put(47.942,99.814){\line(0,1){.1327}}
\put(47.841,99.947){\line(0,1){.1262}}
\put(47.732,100.073){\line(0,1){.1192}}
\put(47.615,100.192){\line(-1,0){.1239}}
\put(47.491,100.304){\line(-1,0){.1305}}
\put(47.361,100.408){\line(-1,0){.1367}}
\put(47.224,100.503){\line(-1,0){.1424}}
\put(47.081,100.59){\line(-1,0){.1475}}
\put(46.934,100.668){\line(-1,0){.152}}
\put(46.782,100.737){\line(-1,0){.156}}
\put(46.626,100.796){\line(-1,0){.1594}}
\put(46.467,100.845){\line(-1,0){.1621}}
\put(46.304,100.885){\line(-1,0){.1643}}
\put(46.14,100.914){\line(-1,0){.1658}}
\put(45.974,100.933){\line(-1,0){.4999}}
\put(45.475,100.928){\line(-1,0){.1654}}
\put(45.309,100.906){\line(-1,0){.1637}}
\put(45.146,100.874){\line(-1,0){.1613}}
\put(44.984,100.831){\line(-1,0){.1584}}
\put(44.826,100.779){\line(-1,0){.1549}}
\put(44.671,100.717){\line(-1,0){.1507}}
\put(44.52,100.645){\line(-1,0){.146}}
\put(44.374,100.565){\line(-1,0){.1407}}
\put(44.233,100.475){\line(-1,0){.1349}}
\put(44.099,100.377){\line(-1,0){.1285}}
\put(43.97,100.27){\line(-1,0){.1217}}
\put(43.848,100.156){\line(0,-1){.1214}}
\put(43.734,100.035){\line(0,-1){.1283}}
\put(43.627,99.907){\line(0,-1){.1346}}
\put(43.529,99.772){\line(0,-1){.1405}}
\put(43.439,99.631){\line(0,-1){.1458}}
\put(43.358,99.486){\line(0,-1){.1505}}
\put(43.286,99.335){\line(0,-1){.1547}}
\put(43.223,99.18){\line(0,-1){.1583}}
\put(43.17,99.022){\line(0,-1){.1612}}
\put(43.128,98.861){\line(0,-1){.1636}}
\put(43.095,98.697){\line(0,-1){.1653}}
\put(43.072,98.532){\line(0,-1){.6656}}
\put(43.085,97.866){\line(0,-1){.1643}}
\put(43.114,97.702){\line(0,-1){.1622}}
\put(43.153,97.54){\line(0,-1){.1595}}
\put(43.202,97.38){\line(0,-1){.1562}}
\put(43.261,97.224){\line(0,-1){.1522}}
\put(43.329,97.072){\line(0,-1){.1477}}
\put(43.406,96.924){\line(0,-1){.1426}}
\put(43.493,96.782){\line(0,-1){.137}}
\put(43.588,96.645){\line(0,-1){.1308}}
\put(43.692,96.514){\line(0,-1){.1241}}
\put(43.803,96.39){\line(1,0){.119}}
\put(43.922,96.273){\line(1,0){.126}}
\put(44.048,96.163){\line(1,0){.1325}}
\put(44.181,96.062){\line(1,0){.1385}}
\put(44.319,95.969){\line(1,0){.144}}
\put(44.463,95.885){\line(1,0){.149}}
\put(44.612,95.81){\line(1,0){.1533}}
\put(44.766,95.744){\line(1,0){.1571}}
\put(44.923,95.688){\line(1,0){.1603}}
\put(45.083,95.641){\line(1,0){.1628}}
\put(45.246,95.605){\line(1,0){.1648}}
\put(45.411,95.579){\line(1,0){.1661}}
\put(45.577,95.563){\line(1,0){.3335}}
\put(45.91,95.562){\line(1,0){.1662}}
\put(46.076,95.577){\line(1,0){.1649}}
\put(46.241,95.603){\line(1,0){.163}}
\put(46.404,95.638){\line(1,0){.1605}}
\put(46.565,95.684){\line(1,0){.1574}}
\put(46.722,95.739){\line(1,0){.1536}}
\put(46.876,95.804){\line(1,0){.1493}}
\put(47.025,95.879){\line(1,0){.1444}}
\put(47.17,95.962){\line(1,0){.139}}
\put(47.309,96.054){\line(1,0){.133}}
\put(47.442,96.155){\line(1,0){.1265}}
\put(47.568,96.264){\line(1,0){.1195}}
\put(47.688,96.38){\line(0,1){.1236}}
\put(47.8,96.504){\line(0,1){.1303}}
\put(47.904,96.634){\line(0,1){.1365}}
\put(48,96.771){\line(0,1){.1422}}
\put(48.087,96.913){\line(0,1){.1473}}
\put(48.165,97.06){\line(0,1){.1519}}
\put(48.235,97.212){\line(0,1){.1559}}
\put(48.294,97.368){\line(0,1){.1593}}
\put(48.344,97.527){\line(0,1){.162}}
\put(48.384,97.689){\line(0,1){.1642}}
\put(48.413,97.853){\line(0,1){.3966}}
%\end
%\emline(47.75,97.75)(61.25,86)
\multiput(47.75,97.75)(.0386819484,-.0336676218){349}{\line(1,0){.0386819484}}
%\end
%\emline(48.25,97.25)(49,96.75)
\multiput(48.25,97.25)(.05,-.0333333){15}{\line(1,0){.05}}
%\end
%\emline(47.75,96.75)(60,85.25)
\multiput(47.75,96.75)(.0359237537,-.0337243402){341}{\line(1,0){.0359237537}}
%\end
%\emline(46.75,96)(52,83.5)
\multiput(46.75,96)(.033653846,-.080128205){156}{\line(0,-1){.080128205}}
%\end
%\emline(46.25,95.75)(50.25,83.25)
\multiput(46.25,95.75)(.033613445,-.105042017){119}{\line(0,-1){.105042017}}
%\end
\put(50.25,83.25){\line(-1,0){.25}}
%\emline(45.25,96)(39,88)
\multiput(45.25,96)(-.033602151,-.043010753){186}{\line(0,-1){.043010753}}
%\end
%\emline(44,96)(38.25,88.75)
\multiput(44,96)(-.033625731,-.042397661){171}{\line(0,-1){.042397661}}
%\end
%\emline(51.75,116.25)(47,100.75)
\multiput(51.75,116.25)(-.033687943,-.109929078){141}{\line(0,-1){.109929078}}
%\end
%\emline(50.75,115.5)(50,116)
\multiput(50.75,115.5)(-.05,.0333333){15}{\line(-1,0){.05}}
%\end
%\emline(50.5,115.75)(45.75,100.25)
\multiput(50.5,115.75)(-.033687943,-.109929078){141}{\line(0,-1){.109929078}}
%\end
%\emline(40.5,111.5)(45.75,99)
\multiput(40.5,111.5)(.033653846,-.080128205){156}{\line(0,-1){.080128205}}
%\end
%\emline(39.5,111.25)(44,100.25)
\multiput(39.5,111.25)(.03358209,-.082089552){134}{\line(0,-1){.082089552}}
%\end
%\emline(33,104.25)(43.75,99.25)
\multiput(33,104.25)(.072147651,-.033557047){149}{\line(1,0){.072147651}}
%\end
%\emline(32.5,103.75)(42.75,98.5)
\multiput(32.5,103.75)(.065705128,-.033653846){156}{\line(1,0){.065705128}}
%\end
%\circle*(65,105){5.385}
\put(63.0398,103.0398){\rule{3.9204\unitlength}{3.9204\unitlength}}
\multiput(63.9134,106.8477)(0,-4.3916){2}{\rule{2.1733\unitlength}{.6961\unitlength}}
\multiput(64.4185,107.4314)(0,-5.1285){2}{\rule{1.1631\unitlength}{.2657\unitlength}}
\multiput(65.4691,107.4314)(-1.3069,0){2}{\multiput(0,0)(0,-5.0643){2}{\rule{.3688\unitlength}{.2015\unitlength}}}
\multiput(65.9742,106.8477)(-2.5263,0){2}{\multiput(0,0)(0,-4.1427){2}{\rule{.578\unitlength}{.4473\unitlength}}}
\multiput(65.9742,107.1826)(-2.2997,0){2}{\multiput(0,0)(0,-4.6134){2}{\rule{.3513\unitlength}{.2483\unitlength}}}
\multiput(65.9742,107.3184)(-2.1816,0){2}{\multiput(0,0)(0,-4.8087){2}{\rule{.2333\unitlength}{.1719\unitlength}}}
\multiput(66.213,107.1826)(-2.6535,0){2}{\multiput(0,0)(0,-4.5483){2}{\rule{.2275\unitlength}{.1832\unitlength}}}
\multiput(66.4397,106.8477)(-3.2041,0){2}{\multiput(0,0)(0,-3.9853){2}{\rule{.3247\unitlength}{.2899\unitlength}}}
\multiput(66.4397,107.0252)(-3.0999,0){2}{\multiput(0,0)(0,-4.2441){2}{\rule{.2205\unitlength}{.1938\unitlength}}}
\multiput(66.6519,106.8477)(-3.5164,0){2}{\multiput(0,0)(0,-3.899){2}{\rule{.2125\unitlength}{.2036\unitlength}}}
\multiput(66.8477,103.9134)(-4.3916,0){2}{\rule{.6961\unitlength}{2.1733\unitlength}}
\multiput(66.8477,105.9742)(-4.1427,0){2}{\multiput(0,0)(0,-2.5263){2}{\rule{.4473\unitlength}{.578\unitlength}}}
\multiput(66.8477,106.4397)(-3.9853,0){2}{\multiput(0,0)(0,-3.2041){2}{\rule{.2899\unitlength}{.3247\unitlength}}}
\multiput(66.8477,106.6519)(-3.899,0){2}{\multiput(0,0)(0,-3.5164){2}{\rule{.2036\unitlength}{.2125\unitlength}}}
\multiput(67.0252,106.4397)(-4.2441,0){2}{\multiput(0,0)(0,-3.0999){2}{\rule{.1938\unitlength}{.2205\unitlength}}}
\multiput(67.1826,105.9742)(-4.6134,0){2}{\multiput(0,0)(0,-2.2997){2}{\rule{.2483\unitlength}{.3513\unitlength}}}
\multiput(67.1826,106.213)(-4.5483,0){2}{\multiput(0,0)(0,-2.6535){2}{\rule{.1832\unitlength}{.2275\unitlength}}}
\multiput(67.3184,105.9742)(-4.8087,0){2}{\multiput(0,0)(0,-2.1816){2}{\rule{.1719\unitlength}{.2333\unitlength}}}
\multiput(67.4314,104.4185)(-5.1285,0){2}{\rule{.2657\unitlength}{1.1631\unitlength}}
\multiput(67.4314,105.4691)(-5.0643,0){2}{\multiput(0,0)(0,-1.3069){2}{\rule{.2015\unitlength}{.3688\unitlength}}}
\put(67.693,105){\line(0,1){.3329}}
\put(67.672,105.333){\line(0,1){.1648}}
\put(67.646,105.498){\line(0,1){.1629}}
\put(67.61,105.661){\line(0,1){.1604}}
\put(67.564,105.821){\line(0,1){.1572}}
\put(67.509,105.978){\line(0,1){.1535}}
\put(67.443,106.132){\line(0,1){.1491}}
\put(67.368,106.281){\line(0,1){.1442}}
\put(67.285,106.425){\line(0,1){.1388}}
\put(67.192,106.564){\line(0,1){.1327}}
\put(67.091,106.697){\line(0,1){.1262}}
\put(66.982,106.823){\line(0,1){.1192}}
\put(66.865,106.942){\line(-1,0){.1239}}
\put(66.741,107.054){\line(-1,0){.1305}}
\put(66.611,107.158){\line(-1,0){.1367}}
\put(66.474,107.253){\line(-1,0){.1424}}
\put(66.331,107.34){\line(-1,0){.1475}}
\put(66.184,107.418){\line(-1,0){.152}}
\put(66.032,107.487){\line(-1,0){.156}}
\put(65.876,107.546){\line(-1,0){.1594}}
\put(65.717,107.595){\line(-1,0){.1621}}
\put(65.554,107.635){\line(-1,0){.1643}}
\put(65.39,107.664){\line(-1,0){.1658}}
\put(65.224,107.683){\line(-1,0){.4999}}
\put(64.725,107.678){\line(-1,0){.1654}}
\put(64.559,107.656){\line(-1,0){.1637}}
\put(64.396,107.624){\line(-1,0){.1613}}
\put(64.234,107.581){\line(-1,0){.1584}}
\put(64.076,107.529){\line(-1,0){.1549}}
\put(63.921,107.467){\line(-1,0){.1507}}
\put(63.77,107.395){\line(-1,0){.146}}
\put(63.624,107.315){\line(-1,0){.1407}}
\put(63.483,107.225){\line(-1,0){.1349}}
\put(63.349,107.127){\line(-1,0){.1285}}
\put(63.22,107.02){\line(-1,0){.1217}}
\put(63.098,106.906){\line(0,-1){.1214}}
\put(62.984,106.785){\line(0,-1){.1283}}
\put(62.877,106.657){\line(0,-1){.1346}}
\put(62.779,106.522){\line(0,-1){.1405}}
\put(62.689,106.381){\line(0,-1){.1458}}
\put(62.608,106.236){\line(0,-1){.1505}}
\put(62.536,106.085){\line(0,-1){.1547}}
\put(62.473,105.93){\line(0,-1){.1583}}
\put(62.42,105.772){\line(0,-1){.1612}}
\put(62.378,105.611){\line(0,-1){.1636}}
\put(62.345,105.447){\line(0,-1){.1653}}
\put(62.322,105.282){\line(0,-1){.6656}}
\put(62.335,104.616){\line(0,-1){.1643}}
\put(62.364,104.452){\line(0,-1){.1622}}
\put(62.403,104.29){\line(0,-1){.1595}}
\put(62.452,104.13){\line(0,-1){.1562}}
\put(62.511,103.974){\line(0,-1){.1522}}
\put(62.579,103.822){\line(0,-1){.1477}}
\put(62.656,103.674){\line(0,-1){.1426}}
\put(62.743,103.532){\line(0,-1){.137}}
\put(62.838,103.395){\line(0,-1){.1308}}
\put(62.942,103.264){\line(0,-1){.1241}}
\put(63.053,103.14){\line(1,0){.119}}
\put(63.172,103.023){\line(1,0){.126}}
\put(63.298,102.913){\line(1,0){.1325}}
\put(63.431,102.812){\line(1,0){.1385}}
\put(63.569,102.719){\line(1,0){.144}}
\put(63.713,102.635){\line(1,0){.149}}
\put(63.862,102.56){\line(1,0){.1533}}
\put(64.016,102.494){\line(1,0){.1571}}
\put(64.173,102.438){\line(1,0){.1603}}
\put(64.333,102.391){\line(1,0){.1628}}
\put(64.496,102.355){\line(1,0){.1648}}
\put(64.661,102.329){\line(1,0){.1661}}
\put(64.827,102.313){\line(1,0){.3335}}
\put(65.16,102.312){\line(1,0){.1662}}
\put(65.326,102.327){\line(1,0){.1649}}
\put(65.491,102.353){\line(1,0){.163}}
\put(65.654,102.388){\line(1,0){.1605}}
\put(65.815,102.434){\line(1,0){.1574}}
\put(65.972,102.489){\line(1,0){.1536}}
\put(66.126,102.554){\line(1,0){.1493}}
\put(66.275,102.629){\line(1,0){.1444}}
\put(66.42,102.712){\line(1,0){.139}}
\put(66.559,102.804){\line(1,0){.133}}
\put(66.692,102.905){\line(1,0){.1265}}
\put(66.818,103.014){\line(1,0){.1195}}
\put(66.938,103.13){\line(0,1){.1236}}
\put(67.05,103.254){\line(0,1){.1303}}
\put(67.154,103.384){\line(0,1){.1365}}
\put(67.25,103.521){\line(0,1){.1422}}
\put(67.337,103.663){\line(0,1){.1473}}
\put(67.415,103.81){\line(0,1){.1519}}
\put(67.485,103.962){\line(0,1){.1559}}
\put(67.544,104.118){\line(0,1){.1593}}
\put(67.594,104.277){\line(0,1){.162}}
\put(67.634,104.439){\line(0,1){.1642}}
\put(67.663,104.603){\line(0,1){.3966}}
%\end
%\emline(47.75,100.75)(46,101.75)
\multiput(47.75,100.75)(-.0583333,.0333333){30}{\line(-1,0){.0583333}}
%\end
%\emline(47.75,99.75)(62.25,105)
\multiput(47.75,99.75)(.092948718,.033653846){156}{\line(1,0){.092948718}}
%\end
%\emline(48,99)(62.75,104)
\multiput(48,99)(.098993289,.033557047){149}{\line(1,0){.098993289}}
%\end
\put(65.75,115.75){\line(0,-1){10.5}}
\put(64.5,116.5){\line(0,-1){10}}
%\emline(67.25,106.75)(72.5,108.75)
\multiput(67.25,106.75)(.0875,.03333333){60}{\line(1,0){.0875}}
%\end
%\emline(67.25,105.75)(73,108)
\multiput(67.25,105.75)(.0858209,.03358209){67}{\line(1,0){.0858209}}
%\end
\put(67.25,104.25){\line(6,-5){4.5}}
%\emline(67.25,103.5)(65.25,103.75)
\multiput(67.25,103.5)(-.25,.03125){8}{\line(-1,0){.25}}
%\end
%\emline(67,103.5)(70.5,100.25)
\multiput(67,103.5)(.036082474,-.033505155){97}{\line(1,0){.036082474}}
%\end
\put(81.5,101.5){$\Delta$}
\put(51,119.5){${\cal B}_1$}
\put(37.75,115){${\cal B}_2$}
\put(63,82.75){${\cal B}_{L-3}$}
\put(51,80){${\cal B}_{L-4}$}
\put(64.25,74.5){Lemma \ref{extdisk} (1)}
\end{picture}

(2) If an annular diagram  contains a least one disk and has no van Kampen subdiagrams with two disks connected
by at least two $\tt$-bands,
then there is a disk $\Pi$ in $\Delta$ and two nonnegative integers
$L', L''$ with $L'+L''\ge L-3$, such that $L'$ (resp. $L''$) consecutive maximal $\tt$-bands ${\mathcal B}_1,\dots
{\mathcal B}_{L'} $ (${\mathcal C}_1,\dots,
{\mathcal C}_{L''}$) start on $\partial\Pi$ , end on the inner (resp., outer) boundary component $\bf p'$ ($\bf p''$) of $\Delta$, and for
any $i\in [1,L'-1]$ (any $i\in [1,L''-1]$)
there are no disks in the van Kampen subdiagram $\Gamma_i$ bounded by ${\mathcal B}_i$, ${\mathcal B}_{i+1},$ $\partial\Pi,$ and $\bf p'$ (by ${\mathcal C}_i$, ${\mathcal C}_{i+1},$ $\partial\Pi,$ and $\bf p''$)
\end{lemma} $\Box$

A maximal $q$-band starting on a disk of  a diagram is called a \index[g]{spoke}
{\it spoke}.
Lemma \ref{extdisk} implies by induction on the number of disks:

\begin{lemma} \label{mnogospits} (see \cite{O12}, Lemma 5.19) If a minimal van Kampen diagram $\Delta$ has $r\ge 1$ disks, then the number of $\tt$-spokes of $\Delta$ ending on the boundary $\partial\Delta$, and therefore
the number of $\tt$-edges in the boundary path of $\Delta$, is greater than  $rL/2$ .
\end{lemma} $\Box$

\begin{lemma} \label{withd} (\cite{OS20}. Lemma 7.7) Let $\Delta$ be a minimal van Kampen diagram.

(1) Assume that a $\theta$-band ${\mathcal T}_0$ crosses $k$
$\tt$-spokes ${\mathcal B}_1,\dots, {\mathcal B}_k$ starting on a disk $\Pi$, and there are no
disks in the subdiagram $\Delta_0$, bounded by these spokes, by ${\mathcal T}_0$ and by $\Pi$. Then $k\le L/2$.

(2) $\Delta$ contains no $\theta$-annuli.
\end{lemma} $\Box$

The proof of the following lemma is given in subsection 7.1.3 of \cite{OS20}.

\begin{lemma} \label{moving} (1) Let $E$ be a van Kampen diagram with the
boundary $\bf x_1y_1x_2y_2$ built of a disk $\Pi$ with boundary ${\bf y_2z^{-1}}$ and  a rim $\theta$-band $\cal T$ with boundary $\bf x_1y_1x_2z$, where $\bf y_1$ and $\bf z$ are the sides of $\cal T$.
Assume that the first and the last cells of $\cal T$ are different $(\theta,\tt)$-cells.
Then there is a diagram
$E'$ with
boundary $\bf x'_1y'_1x'_2y'_2$,
built of a disk $\Pi'$ with boundary ${\bf y'_1(z')^{-1}}$ and  a rim $\theta$-band $\cal T'$, with boundary $\bf x'_1z'x'_2y'_2$,
where $\bf z'$ and $\bf y'_2$ are the sides of $\cal T'$ and
$Lab ({\bf x'_1})\equiv Lab ({\bf x_1})$, $Lab ({\bf x'_2})\equiv Lab ({\bf x_2})$,
$Lab ({\bf y'_1})\equiv Lab ({\bf y_1})$, $Lab ({\bf y'_2})\equiv Lab ({\bf y_2})$.

% This is a LaTeX picture output by TeXCAD.
% File name: [moving.pic].
% Version of TeXCAD: 4.3
% Reference / build: 30-Jun-2012 (rev. 105)
% For new versions, check: http://texcad.sf.net/
% Options on the following lines.
%\grade{\on}
%\emlines{\off}
%\epic{\off}
%\beziermacro{\on}
%\reduce{\on}
%\snapping{\off}
%\pvinsert{% Your \input, \def, etc. here}
%\quality{8.000}
%\graddiff{0.005}
%\snapasp{1}
%\zoom{4.0000}
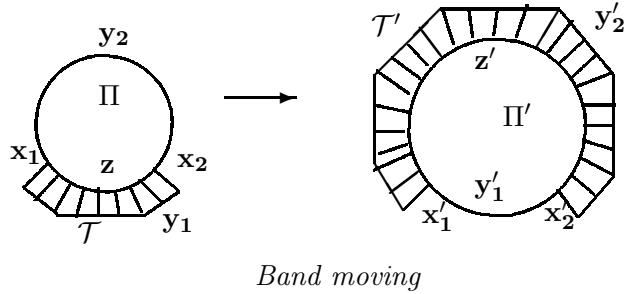
\unitlength 1mm % = 2.845pt
\linethickness{0.4pt}
\ifx\plotpoint\undefined\newsavebox{\plotpoint}\fi % GNUPLOT compatibility
\begin{picture}(174.5,61.5)(0,90)
\thicklines
%\circle(32.75,121.75){18.118}
\put(41.809,121.75){\line(0,1){.4985}}
\put(41.795,122.249){\line(0,1){.497}}
\put(41.754,122.746){\line(0,1){.494}}
\put(41.686,123.24){\line(0,1){.4895}}
\multiput(41.59,123.729)(-.03058,.12087){4}{\line(0,1){.12087}}
\multiput(41.468,124.212)(-.029744,.095203){5}{\line(0,1){.095203}}
\multiput(41.319,124.688)(-.029115,.077851){6}{\line(0,1){.077851}}
\multiput(41.144,125.156)(-.033356,.076131){6}{\line(0,1){.076131}}
\multiput(40.944,125.612)(-.032138,.063583){7}{\line(0,1){.063583}}
\multiput(40.719,126.057)(-.03114,.054003){8}{\line(0,1){.054003}}
\multiput(40.47,126.489)(-.03028,.046407){9}{\line(0,1){.046407}}
\multiput(40.198,126.907)(-.032788,.04467){9}{\line(0,1){.04467}}
\multiput(39.902,127.309)(-.031677,.038518){10}{\line(0,1){.038518}}
\multiput(39.586,127.694)(-.03068,.033379){11}{\line(0,1){.033379}}
\multiput(39.248,128.062)(-.032471,.03164){11}{\line(-1,0){.032471}}
\multiput(38.891,128.41)(-.037579,.032786){10}{\line(-1,0){.037579}}
\multiput(38.515,128.737)(-.039326,.030668){10}{\line(-1,0){.039326}}
\multiput(38.122,129.044)(-.045505,.031619){9}{\line(-1,0){.045505}}
\multiput(37.712,129.329)(-.053073,.032701){8}{\line(-1,0){.053073}}
\multiput(37.288,129.59)(-.054792,.02973){8}{\line(-1,0){.054792}}
\multiput(36.85,129.828)(-.064394,.03048){7}{\line(-1,0){.064394}}
\multiput(36.399,130.041)(-.07697,.031372){6}{\line(-1,0){.07697}}
\multiput(35.937,130.23)(-.094295,.032507){5}{\line(-1,0){.094295}}
\multiput(35.465,130.392)(-.095941,.027268){5}{\line(-1,0){.095941}}
\put(34.986,130.529){\line(-1,0){.4865}}
\put(34.499,130.638){\line(-1,0){.4918}}
\put(34.008,130.721){\line(-1,0){.4956}}
\put(33.512,130.777){\line(-1,0){.4979}}
\put(33.014,130.805){\line(-1,0){.4987}}
\put(32.515,130.806){\line(-1,0){.498}}
\put(32.017,130.779){\line(-1,0){.4958}}
\put(31.522,130.725){\line(-1,0){.4921}}
\put(31.029,130.644){\line(-1,0){.4868}}
\multiput(30.543,130.536)(-.12004,-.0337){4}{\line(-1,0){.12004}}
\multiput(30.063,130.401)(-.0944,-.032201){5}{\line(-1,0){.0944}}
\multiput(29.591,130.24)(-.077071,-.031123){6}{\line(-1,0){.077071}}
\multiput(29.128,130.053)(-.064493,-.030271){7}{\line(-1,0){.064493}}
\multiput(28.677,129.841)(-.054888,-.029553){8}{\line(-1,0){.054888}}
\multiput(28.238,129.605)(-.053178,-.032529){8}{\line(-1,0){.053178}}
\multiput(27.812,129.345)(-.045607,-.031472){9}{\line(-1,0){.045607}}
\multiput(27.402,129.061)(-.039425,-.03054){10}{\line(-1,0){.039425}}
\multiput(27.007,128.756)(-.037685,-.032664){10}{\line(-1,0){.037685}}
\multiput(26.631,128.429)(-.032573,-.031535){11}{\line(-1,0){.032573}}
\multiput(26.272,128.083)(-.030788,-.033279){11}{\line(0,-1){.033279}}
\multiput(25.934,127.716)(-.031801,-.038416){10}{\line(0,-1){.038416}}
\multiput(25.616,127.332)(-.032932,-.044564){9}{\line(0,-1){.044564}}
\multiput(25.319,126.931)(-.03043,-.046309){9}{\line(0,-1){.046309}}
\multiput(25.045,126.514)(-.031315,-.053902){8}{\line(0,-1){.053902}}
\multiput(24.795,126.083)(-.032344,-.063479){7}{\line(0,-1){.063479}}
\multiput(24.568,125.639)(-.033602,-.076023){6}{\line(0,-1){.076023}}
\multiput(24.367,125.183)(-.029367,-.077757){6}{\line(0,-1){.077757}}
\multiput(24.191,124.716)(-.030053,-.095106){5}{\line(0,-1){.095106}}
\multiput(24.04,124.241)(-.03097,-.12077){4}{\line(0,-1){.12077}}
\put(23.916,123.758){\line(0,-1){.4892}}
\put(23.819,123.268){\line(0,-1){.4938}}
\put(23.749,122.775){\line(0,-1){.4969}}
\put(23.707,122.278){\line(0,-1){1.4942}}
\put(23.743,120.784){\line(0,-1){.4942}}
\put(23.81,120.289){\line(0,-1){.4898}}
\multiput(23.904,119.8)(.03018,-.12097){4}{\line(0,-1){.12097}}
\multiput(24.024,119.316)(.029436,-.095299){5}{\line(0,-1){.095299}}
\multiput(24.172,118.839)(.028863,-.077945){6}{\line(0,-1){.077945}}
\multiput(24.345,118.372)(.033109,-.076239){6}{\line(0,-1){.076239}}
\multiput(24.543,117.914)(.031932,-.063687){7}{\line(0,-1){.063687}}
\multiput(24.767,117.468)(.030965,-.054104){8}{\line(0,-1){.054104}}
\multiput(25.015,117.036)(.030129,-.046505){9}{\line(0,-1){.046505}}
\multiput(25.286,116.617)(.032643,-.044776){9}{\line(0,-1){.044776}}
\multiput(25.58,116.214)(.031552,-.038621){10}{\line(0,-1){.038621}}
\multiput(25.895,115.828)(.033629,-.036826){10}{\line(0,-1){.036826}}
\multiput(26.231,115.46)(.032368,-.031745){11}{\line(1,0){.032368}}
\multiput(26.587,115.11)(.037472,-.032907){10}{\line(1,0){.037472}}
\multiput(26.962,114.781)(.039227,-.030795){10}{\line(1,0){.039227}}
\multiput(27.354,114.473)(.045402,-.031767){9}{\line(1,0){.045402}}
\multiput(27.763,114.187)(.052967,-.032872){8}{\line(1,0){.052967}}
\multiput(28.187,113.924)(.054695,-.029908){8}{\line(1,0){.054695}}
\multiput(28.624,113.685)(.064295,-.030689){7}{\line(1,0){.064295}}
\multiput(29.074,113.47)(.076868,-.031621){6}{\line(1,0){.076868}}
\multiput(29.536,113.281)(.09419,-.032812){5}{\line(1,0){.09419}}
\multiput(30.007,113.117)(.095853,-.027579){5}{\line(1,0){.095853}}
\put(30.486,112.979){\line(1,0){.4861}}
\put(30.972,112.867){\line(1,0){.4915}}
\put(31.463,112.783){\line(1,0){.4954}}
\put(31.959,112.726){\line(1,0){.4978}}
\put(32.457,112.696){\line(1,0){.4987}}
\put(32.955,112.693){\line(1,0){.4981}}
\put(33.453,112.719){\line(1,0){.4959}}
\put(33.949,112.771){\line(1,0){.4923}}
\put(34.442,112.851){\line(1,0){.4872}}
\multiput(34.929,112.957)(.12015,.03331){4}{\line(1,0){.12015}}
\multiput(35.409,113.09)(.094504,.031895){5}{\line(1,0){.094504}}
\multiput(35.882,113.25)(.077171,.030873){6}{\line(1,0){.077171}}
\multiput(36.345,113.435)(.06459,.030062){7}{\line(1,0){.06459}}
\multiput(36.797,113.645)(.062838,.033571){7}{\line(1,0){.062838}}
\multiput(37.237,113.88)(.053284,.032356){8}{\line(1,0){.053284}}
\multiput(37.663,114.139)(.045709,.031324){9}{\line(1,0){.045709}}
\multiput(38.075,114.421)(.039524,.030413){10}{\line(1,0){.039524}}
\multiput(38.47,114.725)(.03779,.032542){10}{\line(1,0){.03779}}
\multiput(38.848,115.051)(.032675,.031429){11}{\line(1,0){.032675}}
\multiput(39.207,115.397)(.030896,.03318){11}{\line(0,1){.03318}}
\multiput(39.547,115.761)(.031925,.038312){10}{\line(0,1){.038312}}
\multiput(39.866,116.145)(.033076,.044457){9}{\line(0,1){.044457}}
\multiput(40.164,116.545)(.03058,.04621){9}{\line(0,1){.04621}}
\multiput(40.439,116.961)(.031489,.0538){8}{\line(0,1){.0538}}
\multiput(40.691,117.391)(.032549,.063373){7}{\line(0,1){.063373}}
\multiput(40.919,117.835)(.029013,.065069){7}{\line(0,1){.065069}}
\multiput(41.122,118.29)(.029619,.077661){6}{\line(0,1){.077661}}
\multiput(41.3,118.756)(.030361,.095008){5}{\line(0,1){.095008}}
\multiput(41.452,119.231)(.03136,.12067){4}{\line(0,1){.12067}}
\put(41.577,119.714){\line(0,1){.4888}}
\put(41.676,120.203){\line(0,1){.4935}}
\put(41.747,120.696){\line(0,1){.4967}}
\put(41.792,121.193){\line(0,1){.5571}}
%\end
%\emline(25.25,116.5)(22,113.25)
\multiput(25.25,116.5)(-.033505155,-.033505155){97}{\line(0,-1){.033505155}}
%\end
%\emline(22,113.25)(26.5,109.75)
\multiput(22,113.25)(.043269231,-.033653846){104}{\line(1,0){.043269231}}
%\end
\put(26.5,109.5){\line(1,0){11.5}}
%\emline(38,109.5)(42.75,112.75)
\multiput(38,109.5)(.048969072,.033505155){97}{\line(1,0){.048969072}}
%\end
%\emline(42.75,112.75)(39.25,115.75)
\multiput(42.75,112.75)(-.03932584,.03370787){89}{\line(-1,0){.03932584}}
%\end
%\circle(85,121.25){23.505}
\put(96.753,121.25){\line(0,1){.6176}}
\put(96.736,121.868){\line(0,1){.6159}}
\put(96.688,122.484){\line(0,1){.6125}}
\put(96.607,123.096){\line(0,1){.6074}}
\multiput(96.494,123.703)(-.028964,.120126){5}{\line(0,1){.120126}}
\multiput(96.349,124.304)(-.029364,.098699){6}{\line(0,1){.098699}}
\multiput(96.173,124.896)(-.02958,.083159){7}{\line(0,1){.083159}}
\multiput(95.966,125.478)(-.029671,.071304){8}{\line(0,1){.071304}}
\multiput(95.728,126.049)(-.033377,.069646){8}{\line(0,1){.069646}}
\multiput(95.461,126.606)(-.032881,.060263){9}{\line(0,1){.060263}}
\multiput(95.165,127.148)(-.032402,.052606){10}{\line(0,1){.052606}}
\multiput(94.841,127.674)(-.031929,.04621){11}{\line(0,1){.04621}}
\multiput(94.49,128.183)(-.031454,.040762){12}{\line(0,1){.040762}}
\multiput(94.113,128.672)(-.033553,.039053){12}{\line(0,1){.039053}}
\multiput(93.71,129.141)(-.0328232,.0343715){13}{\line(0,1){.0343715}}
\multiput(93.283,129.587)(-.0345842,.0325991){13}{\line(-1,0){.0345842}}
\multiput(92.834,130.011)(-.03927,.033298){12}{\line(-1,0){.03927}}
\multiput(92.362,130.411)(-.040966,.031188){12}{\line(-1,0){.040966}}
\multiput(91.871,130.785)(-.046416,.031628){11}{\line(-1,0){.046416}}
\multiput(91.36,131.133)(-.052816,.032059){10}{\line(-1,0){.052816}}
\multiput(90.832,131.453)(-.060475,.032488){9}{\line(-1,0){.060475}}
\multiput(90.288,131.746)(-.069861,.032923){8}{\line(-1,0){.069861}}
\multiput(89.729,132.009)(-.081709,.033379){7}{\line(-1,0){.081709}}
\multiput(89.157,132.243)(-.08335,.029039){7}{\line(-1,0){.08335}}
\multiput(88.574,132.446)(-.098887,.028721){6}{\line(-1,0){.098887}}
\multiput(87.98,132.619)(-.120312,.028182){5}{\line(-1,0){.120312}}
\put(87.379,132.759){\line(-1,0){.6081}}
\put(86.771,132.869){\line(-1,0){.613}}
\put(86.158,132.946){\line(-1,0){.6162}}
\put(85.541,132.99){\line(-1,0){1.2352}}
\put(84.306,132.982){\line(-1,0){.6156}}
\put(83.69,132.929){\line(-1,0){.612}}
\put(83.079,132.845){\line(-1,0){.6067}}
\multiput(82.472,132.728)(-.119936,-.029744){5}{\line(-1,0){.119936}}
\multiput(81.872,132.579)(-.098506,-.030005){6}{\line(-1,0){.098506}}
\multiput(81.281,132.399)(-.082965,-.03012){7}{\line(-1,0){.082965}}
\multiput(80.7,132.188)(-.071109,-.030134){8}{\line(-1,0){.071109}}
\multiput(80.131,131.947)(-.061713,-.03007){9}{\line(-1,0){.061713}}
\multiput(79.576,131.676)(-.060048,-.033272){9}{\line(-1,0){.060048}}
\multiput(79.036,131.377)(-.052394,-.032743){10}{\line(-1,0){.052394}}
\multiput(78.512,131.049)(-.046001,-.032229){11}{\line(-1,0){.046001}}
\multiput(78.006,130.695)(-.040557,-.031718){12}{\line(-1,0){.040557}}
\multiput(77.519,130.314)(-.0358468,-.0312053){13}{\line(-1,0){.0358468}}
\multiput(77.053,129.909)(-.0341573,-.033046){13}{\line(-1,0){.0341573}}
\multiput(76.609,129.479)(-.0323735,-.0347954){13}{\line(0,-1){.0347954}}
\multiput(76.188,129.027)(-.033042,-.039486){12}{\line(0,-1){.039486}}
\multiput(75.792,128.553)(-.033732,-.04491){11}{\line(0,-1){.04491}}
\multiput(75.421,128.059)(-.031325,-.046621){11}{\line(0,-1){.046621}}
\multiput(75.076,127.546)(-.031715,-.053023){10}{\line(0,-1){.053023}}
\multiput(74.759,127.016)(-.032094,-.060685){9}{\line(0,-1){.060685}}
\multiput(74.47,126.47)(-.032468,-.070074){8}{\line(0,-1){.070074}}
\multiput(74.21,125.909)(-.032847,-.081924){7}{\line(0,-1){.081924}}
\multiput(73.98,125.335)(-.033246,-.09746){6}{\line(0,-1){.09746}}
\multiput(73.781,124.751)(-.033693,-.118886){5}{\line(0,-1){.118886}}
\multiput(73.612,124.156)(-.027399,-.120493){5}{\line(0,-1){.120493}}
\put(73.475,123.554){\line(0,-1){.6088}}
\put(73.37,122.945){\line(0,-1){.6135}}
\put(73.297,122.331){\line(0,-1){1.2343}}
\put(73.248,121.097){\line(0,-1){.6174}}
\put(73.273,120.48){\line(0,-1){.6152}}
\put(73.329,119.865){\line(0,-1){.6114}}
\multiput(73.418,119.253)(.03024,-.15147){4}{\line(0,-1){.15147}}
\multiput(73.539,118.647)(.030523,-.11974){5}{\line(0,-1){.11974}}
\multiput(73.692,118.049)(.030645,-.098308){6}{\line(0,-1){.098308}}
\multiput(73.876,117.459)(.030659,-.082767){7}{\line(0,-1){.082767}}
\multiput(74.09,116.879)(.030595,-.070912){8}{\line(0,-1){.070912}}
\multiput(74.335,116.312)(.030471,-.061516){9}{\line(0,-1){.061516}}
\multiput(74.609,115.758)(.033662,-.05983){9}{\line(0,-1){.05983}}
\multiput(74.912,115.22)(.033083,-.05218){10}{\line(0,-1){.05218}}
\multiput(75.243,114.698)(.032527,-.045791){11}{\line(0,-1){.045791}}
\multiput(75.601,114.194)(.031981,-.04035){12}{\line(0,-1){.04035}}
\multiput(75.985,113.71)(.0314377,-.0356431){13}{\line(0,-1){.0356431}}
\multiput(76.393,113.247)(.0332674,-.0339417){13}{\line(0,-1){.0339417}}
\multiput(76.826,112.806)(.0350052,-.0321466){13}{\line(1,0){.0350052}}
\multiput(77.281,112.388)(.0397,-.032784){12}{\line(1,0){.0397}}
\multiput(77.757,111.994)(.045129,-.033439){11}{\line(1,0){.045129}}
\multiput(78.254,111.626)(.046824,-.031022){11}{\line(1,0){.046824}}
\multiput(78.769,111.285)(.053228,-.03137){10}{\line(1,0){.053228}}
\multiput(79.301,110.972)(.060893,-.031699){9}{\line(1,0){.060893}}
\multiput(79.849,110.686)(.070284,-.032012){8}{\line(1,0){.070284}}
\multiput(80.411,110.43)(.082136,-.032314){7}{\line(1,0){.082136}}
\multiput(80.986,110.204)(.097674,-.032611){6}{\line(1,0){.097674}}
\multiput(81.572,110.008)(.119103,-.03292){5}{\line(1,0){.119103}}
\multiput(82.168,109.844)(.15084,-.03327){4}{\line(1,0){.15084}}
\put(82.771,109.711){\line(1,0){.6095}}
\put(83.381,109.609){\line(1,0){.614}}
\put(83.995,109.54){\line(1,0){.6168}}
\put(84.611,109.504){\line(1,0){.6178}}
\put(85.229,109.5){\line(1,0){.6172}}
\put(85.846,109.528){\line(1,0){.6149}}
\put(86.461,109.589){\line(1,0){.6108}}
\multiput(87.072,109.681)(.15127,.03122){4}{\line(1,0){.15127}}
\multiput(87.677,109.806)(.119539,.031301){5}{\line(1,0){.119539}}
\multiput(88.275,109.963)(.098107,.031283){6}{\line(1,0){.098107}}
\multiput(88.864,110.151)(.082566,.031197){7}{\line(1,0){.082566}}
\multiput(89.442,110.369)(.070711,.031056){8}{\line(1,0){.070711}}
\multiput(90.007,110.617)(.061317,.03087){9}{\line(1,0){.061317}}
\multiput(90.559,110.895)(.053649,.030645){10}{\line(1,0){.053649}}
\multiput(91.096,111.202)(.051964,.033422){10}{\line(1,0){.051964}}
\multiput(91.615,111.536)(.045578,.032824){11}{\line(1,0){.045578}}
\multiput(92.117,111.897)(.040141,.032243){12}{\line(1,0){.040141}}
\multiput(92.598,112.284)(.035438,.0316688){13}{\line(1,0){.035438}}
\multiput(93.059,112.696)(.0337247,.0334874){13}{\line(1,0){.0337247}}
\multiput(93.497,113.131)(.0319183,.0352135){13}{\line(0,1){.0352135}}
\multiput(93.912,113.589)(.032526,.039912){12}{\line(0,1){.039912}}
\multiput(94.303,114.068)(.033145,.045345){11}{\line(0,1){.045345}}
\multiput(94.667,114.566)(.030716,.047024){11}{\line(0,1){.047024}}
\multiput(95.005,115.084)(.031023,.053431){10}{\line(0,1){.053431}}
\multiput(95.315,115.618)(.031302,.061097){9}{\line(0,1){.061097}}
\multiput(95.597,116.168)(.031554,.07049){8}{\line(0,1){.07049}}
\multiput(95.849,116.732)(.031779,.082344){7}{\line(0,1){.082344}}
\multiput(96.072,117.308)(.031975,.097884){6}{\line(0,1){.097884}}
\multiput(96.264,117.895)(.032145,.119315){5}{\line(0,1){.119315}}
\multiput(96.424,118.492)(.03229,.15105){4}{\line(0,1){.15105}}
\put(96.554,119.096){\line(0,1){.6101}}
\put(96.651,119.706){\line(0,1){.6144}}
\put(96.716,120.321){\line(0,1){.9292}}
%\end
%\emline(75.75,113.5)(72.5,110.25)
\multiput(75.75,113.5)(-.033505155,-.033505155){97}{\line(0,-1){.033505155}}
%\end
\put(72.5,110.25){\line(-3,5){3.75}}
\put(68.75,116.5){\line(0,1){11.75}}
\put(68.75,128.25){\line(1,1){8.75}}
\put(77.5,137){\line(1,0){15.75}}
%\emline(93.25,137)(100.5,128.25)
\multiput(93.25,137)(.03372093,-.040697674){215}{\line(0,-1){.040697674}}
%\end
\put(100.5,128.25){\line(0,-1){13.5}}
%\emline(100.5,114.75)(95.75,109.25)
\multiput(100.5,114.75)(-.033687943,-.039007092){141}{\line(0,-1){.039007092}}
%\end
%\emline(95.75,109.25)(93,113)
\multiput(95.75,109.25)(-.03353659,.04573171){82}{\line(0,1){.04573171}}
%\end
%\emline(26.75,115)(23.5,111.5)
\multiput(26.75,115)(-.033505155,-.036082474){97}{\line(0,-1){.036082474}}
%\end
%\emline(28.25,113.75)(26.25,110.25)
\multiput(28.25,113.75)(-.03333333,-.05833333){60}{\line(0,-1){.05833333}}
%\end
%\emline(36.25,113.25)(37.75,110)
\multiput(36.25,113.25)(.03333333,-.07222222){45}{\line(0,-1){.07222222}}
%\end
%\emline(37.5,114.25)(40.5,111.25)
\multiput(37.5,114.25)(.03370787,-.03370787){89}{\line(0,-1){.03370787}}
%\end
\put(32,113.25){\line(0,-1){3.5}}
%\emline(34.25,112.75)(34.75,109.75)
\multiput(34.25,112.75)(.0333333,-.2){15}{\line(0,-1){.2}}
%\end
%\emline(30,113)(29,109.25)
\multiput(30,113)(-.0333333,-.125){30}{\line(0,-1){.125}}
%\end
\put(74.75,115.5){\line(-4,-3){4}}
%\emline(73.5,117.75)(68.5,115.5)
\multiput(73.5,117.75)(-.07462687,-.03358209){67}{\line(-1,0){.07462687}}
%\end
%\emline(69,128.75)(73.25,125.75)
\multiput(69,128.75)(.04775281,-.03370787){89}{\line(1,0){.04775281}}
%\end
%\emline(174.5,141.5)(173.75,141)
\multiput(174.5,141.5)(-.05,-.0333333){15}{\line(-1,0){.05}}
%\end
\put(68.75,120.75){\line(1,0){4}}
%\emline(68.5,124.25)(73.25,123.25)
\multiput(68.5,124.25)(.1583333,-.0333333){30}{\line(1,0){.1583333}}
%\end
%\emline(69.25,117.75)(73,119)
\multiput(69.25,117.75)(.09868421,.03289474){38}{\line(1,0){.09868421}}
%\end
%\emline(77.75,137)(79.25,132.25)
\multiput(77.75,137)(.03333333,-.10555556){45}{\line(0,-1){.10555556}}
%\end
%\emline(73.25,132.25)(75.75,129.5)
\multiput(73.25,132.25)(.03333333,-.03666667){75}{\line(0,-1){.03666667}}
%\end
\put(74.75,134.5){\line(3,-4){3}}
%\emline(71,130.5)(74.75,127.75)
\multiput(71,130.5)(.04573171,-.03353659){82}{\line(1,0){.04573171}}
%\end
\put(93.25,137){\line(-3,-5){3}}
%\emline(85,136.75)(84.5,133.5)
\multiput(85,136.75)(-.0333333,-.2166667){15}{\line(0,-1){.2166667}}
%\end
%\emline(89,137.5)(87.5,133.25)
\multiput(89,137.5)(-.03333333,-.09444444){45}{\line(0,-1){.09444444}}
%\end
%\emline(81.25,137)(81.75,133.25)
\multiput(81.25,137)(.0333333,-.25){15}{\line(0,-1){.25}}
%\end
%\emline(100.25,128.5)(95.75,126.75)
\multiput(100.25,128.5)(-.08653846,-.03365385){52}{\line(-1,0){.08653846}}
%\end
%\emline(96.75,133)(93.25,129.75)
\multiput(96.75,133)(-.036082474,-.033505155){97}{\line(-1,0){.036082474}}
%\end
%\emline(95,135.5)(91.75,131.25)
\multiput(95,135.5)(-.033505155,-.043814433){97}{\line(0,-1){.043814433}}
%\end
%\emline(98,131.25)(98.5,130.75)
\multiput(98,131.25)(.0333333,-.0333333){15}{\line(0,-1){.0333333}}
%\end
%\emline(94.75,127.75)(98.25,131)
\multiput(94.75,127.75)(.036082474,.033505155){97}{\line(1,0){.036082474}}
%\end
%\emline(96,117.25)(100.25,115)
\multiput(96,117.25)(.06343284,-.03358209){67}{\line(1,0){.06343284}}
%\end
\put(96.5,121.5){\line(1,0){4}}
\put(96.5,124.25){\line(1,0){4.25}}
%\emline(96.5,119.5)(100.5,118.25)
\multiput(96.5,119.5)(.10526316,-.03289474){38}{\line(1,0){.10526316}}
%\end
%\emline(94.75,115)(98,112.25)
\multiput(94.75,115)(.03963415,-.03353659){82}{\line(1,0){.03963415}}
%\end
\put(32,132.75){$\bf y_2$}
\put(20.25,117){$\bf x_1$}
\put(42.5,116){$\bf x_2$}
\put(40.5,108){$\bf y_1$}
\put(32.25,115.25){$\bf z$}
\put(32,124.25){$\Pi$}
\put(29.25,106){$\cal T$}
\put(75,108.5){$\bf x'_1$}
\put(91.25,108.75){$\bf x'_2$}
\put(82,112.5){$\bf y'_1$}
\put(82,129){$\bf z'$}
\put(85.75,121.5){$\Pi'$}
\put(98,135){$\bf y'_2$}
\put(68.25,133.5){$\cal T'$}
\put(48.75,125.25){\vector(1,0){9.5}}
\put(52.75,100.25){Band moving}
\end{picture}

2) Let $\Delta$ be a van Kampen diagram with  boundary $\bf pq$ and $\Delta$
a union of a minimal diagram
$\Gamma$ with $r>0$ disks and a rim $\theta$-band $\cal T$ with  side $\bf p$.
Assume that there are two $\tt$-spokes
in $\Delta$ starting on a disk $\Pi$ and ending on $\bf p$. Then
there exists a van Kampen diagram $\Delta'$
with boundary $\bf p'q'$, and $\Delta'$ is a union
of a minimal diagram $\Gamma'$ with $r'<r$ disks and a rim $\theta$-band $\cal T'$ having side $\bf p'$, where  $Lab({\bf p'})=Lab({\bf p})$ in the group $G$, and $Lab({\bf q'})\equiv Lab({\bf q})$.
\end{lemma} $\Box$

\begin{lemma} \label{quadr} (\cite{OS20}, Lemmas 8.2). The group $G$ has quadratic Dehn function.
\end{lemma} $\Box$

\section{Isomorphism problem for groups with quadratic Dehn function}\label{iso}

In this section, we assume that the construction of the machine $\bf M$ is based on the S-machine $M_1$ provided by Lemma \ref{S}>

\begin{lemma} \label{norep} Let a reduced computation $\ccc\colon W\to\dots\to W'$ of $\mmm$ with standard base
have no rules of Sets $\Theta_1$ and $\Theta_2$. If $W'\equiv W$,
then the computation $\ccc$ is empty.
\end{lemma}

\proof Proving by contradiction, one may assume that the history $H$ of $\ccc$ is a non-empty
cyclically reduced word, because otherwise one could replace $\ccc$ with a shorter computation
with equal the first and the last configurations.

Assume first that $\ccc$ is a one-step computation. This step cannot be $\Theta_3$ or $\Theta_5$, since the computations with rules from these sets multiplies the words in history sectors by a copy of $H^{\pm 1}$. So $\ccc$ is of type $\Theta_4$.
This assumption reformulates our problem as the same problem for the $S$-machine $\mmm_5$. If
$\ccc$ has a $\chi$-rule of $\mmm_5$, then one may consider the computation $W\to\dots\to W$
with history $H^2$, where this rule occurs at least twice, which contradicts  Lemma \ref{M31}.
Therefore there are no $\chi$-rules in $H$, and so $\ccc$ is just a computation of either
$\bf RL$ or $\bf LR$, or $\mmm_2$. The first and second cases contradicts Lemma \ref{Hprim}, because
the length of powers $H^s$ are unbounded. In the later case we restrict the computation
of $\mmm_2$ to a history sector:  $V\to\dots\to V$, where a computation with any history $H'$ multiplies
the tape word from the left and from the right by copies of $H^{\pm 1}$ in disjoint alphabets.
Clearly one cannot obtain a repetition, provided the word $H$ is non-empty.

If $H$ has at least two steps, then its step history (or a cyclic permutation of it) is a power of $(3)(4)(5)(4)$
by Lemma \ref{212}. So $H^{\pm 2}$  has to contain a subword $H_1H_2H_3$, where
$H_2$ has step history $(43)(3)(34)$, $H_3$ and $H_1^{-1}$ are of type $(4)(45)(5)(54)$.
Since $H_2$ does change history sectors but does not change the input ones, the
computations with histories $H_3$ and $H_1^{-1}$ start working with configurations
having  equal input sectors but different history sectors. Considering the
subcomputations in $(34)(4)(45)$ corresponding to the work of
$\mmm_5$, $\mmm_3$, and $\mmm_2$ (as we did in the beginning of this proof),
we see that the $S$-machine $\mmm_2$ can connect both $I_2(\alpha^k,H')$ and $I_2(\alpha^k, H'')$
with $A_2(H')$ and $A_2(H'')$, resp., where $H'\ne H''$. It follows that there are two different
reduced computations of $\mmm_1$ accepting the same input word $\alpha^k$, contrary
to Lemma \ref{S}.  The lemma is proved.
\endproof

Recall that the rule $\theta(23)$ locks all sectors of the standard base of $\mmm$
except for the input sector ${\tilde R}_0{\tilde P}_1$ and its mirror copy. Hence every
$\theta(23)^{-1}$-admissible configuration has the form $W(k,k')\equiv w_1\alpha^kw_2(\alpha')^{-k'}w_3$,
where $k$ and $k'$ are integers and $w_1, w_2, w_3$ are fixed word in state letters;
$w_1$ starts with $\tt$.

\begin{lemma} \label{conj} (\cite{OS20}, Lemma 8.3). A word $W(k,k)$ is a conjugate of the word $W_{ac}$
in the group $G$ (and in the group $M$) if and only if the subword $\alpha^k$
is accepted by the Turing machine $\mmm_0$.
\end{lemma} $\Box$

\begin{lemma} \label{order} For arbitrary integer $k$, the word $W(k,k)$ has order $L$ in $G$.
\end{lemma}

\proof Starting with the word $W_{st}$, a computation of Set $\Theta_1$ can insert the words $\alpha^k$ and $(\alpha')^{-k}$ in the  input sectors. So, after the  application of the connecting rule $\theta(12)$,
the rules of Set $\Theta_2$ can successfully check the content of the input sectors, and
the rule $\theta(23)$ gives us the word $W(k,k)$, the last configuration of this computation. Therefore the word $W(k,k)$ is
accessible. Thus, the power $W(k,k)^L$ is a disk word equal to
$1$ in $G$ by Lemma \ref{trivial}.

Assume that $W(k,k)^l=1$ for a positive $l\le L/2$. Then on the one hand, the minimal diagram $\Delta$ for this equality has $l\le L/2$ $\tt$-edges in the boundary, and so it has no disks by Lemma \ref{mnogospits}. On the other hand, since
all $\tt$-letters of the boundary label occur with the same sign, a maximal
$\tt$-band of $\Delta$ cannot start and end on the boundary, and therefore the word $W(k,k)$ has no $\tt$-letters, a contradiction. The lemma is proved.
\endproof

\begin{lemma}\label{finor} Every element of finite order is conjugate of a power of some
word $W(k,k)$.
\end{lemma}

\proof Consider a minimal diagram $\Delta$ for an equality $U^s=1$, $s>0$, assuming that $U$
has minimal number of $\theta$-letters in the conjugacy class and, under this assumption,
$U$ has  minimal number of $q$-letters. If $\Delta$ has no disk, then it has
a rim $\theta$-band $\mathcal T$. The exterior side $\bf y$ of $\mathcal T$ cannot have length $\ge ||U||$,
since then the whole boundary of $\Delta$ has to have no $\theta$-letters, contrary to
the existence of the rim $\theta$-band. If $||{\mathbf y}|| = ||U||-1$ then the ends
of $\mathcal T$ must have the same label, since the boundary label of $\Delta$ has period $U$, but
this is not possible since one of these $\theta$-letters is positive and the other one
is negative for the boundary label of a band. If  $||{\mathbf y}|| \le ||U||-2$, then one can replace the common boundary
path of $\mathcal T$ and $\partial\Delta$ with a path separating $\mathcal T$ from $\Delta$. Hence a cyclic permutation of $U$ can be replaced with an equal in $G$ word having less $\theta$-letters, a contradiction.

Therefore $\Delta$ contains disks. If it has $\theta$-edges in the boundary, then there is
a maximal $\theta$-band $\mathcal T$ such that the van Kampen diagram bounded by $\mathcal T$ and a subpath
$\bf y$ of $\partial\Delta$ has no $\theta$-cells. Therefore $\bf y$ has no $\theta$-edges, and one comes to
a contradiction, as in the previous paragraph. Hence $\partial\Delta$ has no $\theta$-edges
and therefore $\Delta$ has no $\theta$-edges by Lemma \ref{withd} (2).

Let us consider a disk $\Pi$ provided by Lemma \ref{extdisk}. Since $\Delta$ has no $\theta$-cells,
there is a common subpath $\bf p$ of $\partial\Pi$ and $\partial\Delta$ containing $L-3$ $\tt$-letters.
The word $U$ has at most $L/2$ $\tt$-letters since otherwise there is a cyclic permutation
of $U$ containing a subword of $\Lab({\mathbf p})$ with $>L/2$ $\tt$-letters. So the disk relation
makes $U$ conjugate in $G$ to a word $U'$ with $|U'|_{\theta}=|U|_{\theta}=0$ and with $|U'|_q<|U|_q$,
a contradiction.

Since $L/2<L-3$ and $\Lab({\mathbf p})$ is a subword of a power of $U$, there is a $\tt$-letter
occurring in $\Lab({\mathbf p})$ at least twice. Therefore the disk word on $\partial\Pi$
has no letters with superscripts. Hence it is a power $V^L$, where the letters of $V$ have
no superscripts and a cyclic permutation of $U$ is a power $V^l$, where $|l|<L$.
It remains to show that the word $V$ is a conjugate of some $W(k,k)$.

By Definition \ref{dw}, the word $V$ is accessible. Hence it can be connected by a computation $\ccc$
either with $W_{st}$ or with $W_{ac}$. In the former case $\ccc$ has a rule $\theta(23)^{-1}$
since the letters of $V$ have no superscripts. The maximal prefix $\mathcal D$ of $\ccc$ containing no
rules $\theta(23)^{-1}$ connects $V$ with some $\theta(23)^{-1}$-admissible word $W(k,k)$, and
the configurations of $\mathcal D$ have no letters with superscripts. Hence the corresponding
to $\mathcal D$ trapezium (see Lemma \ref{simul} (2)) has equal side labels, and so $V$ is conjugate of $W(k,k)$, as required.

In the later case, we may assume  that the computation $\ccc$ connecting $V$ and $W_{ac}$ has
no rules $\theta(23)^{-1}$ (otherwise one could argue as above), and so it has no rules of Sets 1 and 2. Therefore the
word $V$ is a conjugate of $W_{ac}$ by Lemma \ref{simul} (2), since
the side labels of the corresponding trapezium have no superscripts and
therefore are equal. In turn,  by Lemma \ref{conj}, the word $W_{ac}$ is a conjugate of any $W(k_0,k_0)$
if the word $\alpha^{k_0}$ is accepted by the machine $M_0$. This completes the proof of the lemma.
\endproof

 \begin{lemma} \label{cent} (1) For every $k$, the cyclic subgroup $\langle W(k,k)\rangle $ is malnormal in $G$, that is
$I=\langle W(k,k)\rangle \cap Z\langle W(k,k)\rangle Z^{-1} =\{1\}$ if $Z\notin \langle W(k,k)\rangle $. The centralizer of an element $g\in G$ of order $L$ is equal to
the cyclic subgroup $\langle g \rangle$.

(2) The subgroup $\langle W(k,k)\rangle $ has trivial intersection with every conjugate subgroup
of $\langle W_{ac}\rangle $ provided the word $\alpha^k$ is not accepted by the machine $\mmm_0$.
 \end{lemma}

\proof (1) To prove the statement about the centralizers, one may assume  by Lemmas \ref{order} and \ref{finor}, that the element $g$ of order $L$ is represented by some word $W\equiv W(k,k)$. Therefore it suffices to prove the first claim of the lemma.

Assuming that
the intersection $I$ is nontrivial, we can find two exponents $s$ and $r$ such that
$ W^s = Z W^r Z^{-1}$ and $0<s\le L/3$, $|r|\le L/2$ if the order of $I$ is
odd, or $s=r=L/2$ otherwise.

We should prove that the word $Z$ is equal to a power of $W$ in $G$. For this goal, we consider a minimal van Kampen diagram $\Delta$ for this equality $ W^s = Z W^r Z^{-1}$
and identifying the subpaths of the boundary labeled by $Z$, we obtain an
annular diagram $\Gamma$ whose two (clockwise) boundary labels read from some vertices $o$ and $o'$ are $W^s$ and $W^r$, and there is a simple path $\bf z$ connecting $o$ with $o'$ and labeled by $Z$. (To obtain $\Gamma$ homeomorphic to topological annulus and to make the path $\bf z$ simple, one
can use 0-cells corresponding to trivial relations as in Section 11 of \cite{book}.)

One can cancel out the pairs of mirror cells if $\Gamma$ is not reduced. Also
if $\Gamma$ contains a pair of disks $\Pi_1$, $\Pi_2$ connected by two $\tt$-bands, and
these disks and the bands do not surround the hole of $\Gamma$, one can replace a van Kampen subdiagram
having two disks with a diagram without disks by Lemma \ref{2dis}
%as it is explained in \cite{OS20}, 7.1.2.
The obtained reduced diagram
$\Gamma'$ has a simple path ${\mathbf z}'$ connecting $o$ and $o'$, whose label is equal in $G$ to
$\Lab({\mathbf z})\equiv Z$. (See Section 13.6 of \cite{book}.)

The reduced diagram $\Gamma'$ has no disks. Indeed,  the two boundary components of $\Gamma'$ have at most
$L/3+L/2$ $\tt$-edges if $s\le L/3$,
but by Lemma \ref{extdisk} (2), an annular diagram with disks has to have at least $L-3>\frac 56 L+1$ $\tt$-edges on the boundary.
If $s=r=L/2$
we have a contradiction again since all the $\tt$-letters of $W^s$ and $W^r$ are positive, but
the disk has to have spokes ending on both boundaries of $\Gamma'$ since $L-3>L/2$.
Hence there
are only  maximal $\tt$-bands,
but there are no disks in $\Gamma'$.

If $\Gamma'$ has no cells, then it is a diagram over the free group, and so $s=r$ and the word $\Lab({\mathbf z}')$ commuting with $W^s$ in the free group, equal to a power of $W$, as required. Arguing  by contradiction, assume that $\Gamma'$ has $\theta$-cells. Then
all maximal $\theta$-bands of $\Gamma'$ are $\theta$-annuli surrounding the hole of $\Gamma'$ by Lemma \ref{NoAnnul}. Cutting along a side of a maximal $\tt$-band $\mathcal T$ we obtain a reduced van Kampen diagram $\Gamma''$ over $M$, which is a trapezium
of height $h\ge 1$ with equal side labels. Therefore the maximal $\theta$-bands of $\Gamma'$
have no superscripts in the labels of their  cells, because it follows from Lemma \ref{simul} (1), that the label of the
right side of a trapezium with $s$ maximal $\tt$-bands must be the $\pm s$-shift of the label
of the left side of it, but $s<L$.

By Lemma \ref{simul} (1), a  subtrapezium  of $\Gamma''$ with the same history gives us a non-empty computation $W\to\dots\to W$
without rules of Sets $\Theta_1$ and $\Theta_2$. This contradicts Lemma \ref{norep}.

(2) We have the same proof as in item (1) considering now a hypothetical conjugation of $W(k,k)^s$ and $W_{ac}^r$. Clearly these cyclically reduced words are not conjugate in the free group since they involve different $q$-letters.
Then as above, one obtain a computation $W(k,k)\to\dots\to W_{ac}$ without rules from Sets $\Theta_1$ and $\Theta_2$, contrary to Lemma \ref{I6A6}.
The lemma is proved.
\endproof

\begin{rk}\label{Wkk} If the machine $\bf M_0$ is chosen with nonrecursive
language of input words $\alpha^k$, then by Lemmas \ref{conj} and \ref{cent} (2),
there exists no algorithm deciding whether some nontrivial powers of the words $W(k,k)$ and $W_{ac}$ are conjugate in $G$ or not.
\end{rk}

\medskip

Basing on Lemmas \ref{conj} and \ref{order}, we can introduce the HNN-extension $G_k$ ($k=1,2,\dots$)
of the group $G$ by adding a stable letter $x$ to the set of generators and the relation
$xW(k,k)x^{-1}=W_{ac}$ to the set of defining relations of $G$.

We need a property of HNN-extensions similar to the property of amalgamated products
obtained in \cite{BC}. Both the formulation and the proof of the first part of the following lemma were left
to the reader in \cite{BC}.

\begin{lemma} \label{same} Let a function $f$ bound from above the Dehn function of a finitely presented group
$A$ and $f$ is super-additive, i.e.,$ f(n_1)+f(n_2)\le f(n_1+n_2)$ for every integers $n_1,n_2\ge 0$,
and $f(1)\ge 1$.
Let $B $ be an HNN-extension of $A$ with finite associated subgroups: $xCx^{-1}=D$. Then
the Dehn function $g(n)$ of $B$ is bounded from above by a function equivalent to $f(n)$.
In particular, every group $G_k$ has quadratic Dehn function.
\end{lemma}

\proof  We have a finite presentation of $B$ with one extra letter $x$ and finitely
many relations $xU_ix^{-1}=V_i$, where all elements of the subgroups $C$ and $D$ are
presented by some words $U_i$-s and $V_i$-s. Let $c-1$ be the maximum of the lengths
of all $U_i$-s and $V_i$-s.

Assume that a word $W\equiv W_1x^{\pm 1}W_2x^{\pm 1}\dots $ is equal to 1 in $B$,
where the words $W_j$-s have no $x$-letters. We will induct on $s=s(W)=||W||_B+cr$, where
$r$ is the number of $x$-letters in $W$, with trivial base $s=0$, to show that the area $\area(W)$ in $B$
does not exceed  $f(s)$. If $W$ has no $x$-letters, then $W=1$ in $G$, and therefore $\area(W)\le f(s)$.

If the word $W$ has $x$-letters, then the word $W$ has a pinch subword by Britton's lemma (see Section IV.2 in \cite{LS}), i.e., a subword $xW_jx^{-1}$ (a subword  $x^{-1}W_jx$ ), where $W_j$ is equal
in $A$ to some word $U_i$ (resp. to some $V_j$). Therefore one can replace $W_j$ with $U_i$
using an auxiliary diagram of area $\le f(||W_j||+||U_i||)\le f(n_j+c-1)$, where $n_j=||W_j||$. Then
the application of one conjugacy by $x$ replaces the subword $xU_ix^{-1}$ with $V_i$.
We obtain a word $W'$ with $s(W')< s(W)-n_j-c$ since the number of $x$-letters is decreased
by $2$. Therefore $$\area(W)< f(s-n_j-c)+ f(n_j+c-1) +1 \le f(s-n_j-c)+ f(n_j+c))\le f(s)\le f((c+1)||W||),$$
and so $g(n) \le f((c+1)n)$ for every $n\ge 0$,  which proves the first statement of the lemma.
It implies the second one by Lemmas \ref{quadr} and \ref{order}.
\endproof

\begin{rk} It is unknown if there is a finitely presented group whose Dehn function
is not equivalent to a super-additive function. This problem was raised by V.S. Guba
and M.V. Sapir in \cite{GS}.
\end{rk}

\begin{lemma}\label{maln} Let $B$ be an HNN-extension of a group $A$ with associate
malnormal subgroups $C$ and $D$:  $xCx^{-1}= D$. Assume also that $gCg^{-1}\cap D=\{1\}$
for every element $g\in A$. Then the centralizer of any nontrivial element $h\in A$
in $B$ is equal to the centralizer of $h$ in $A$.
\end{lemma}

\proof  Let an element $z$ commute with $h$ in $B$. Assume first that
 it has only one stable letter $x$ in the normal form: $z=g_1xg_2$, where $g_1, g_2\in A$.
Then the equality $g_1xg_2hg_2^{-1}x^{-1}g_1^{-1}=zhz^{-1}=h \in C$ implies
that the subword $x(g_2hg_2^{-1})x^{-1}$ is a pinch, and so  $g_2hg_2^{-1}\in C\backslash\{1\}$.
Then $xg_2hg_2^{-1}x^{-1}=d\in D\backslash\{1\}$, but the conjugate in $A$ element
$g_1^{-1}dg_1=h$ belongs to $C$, contrary to the assumption of the lemma.

Now assume that the normal form of $z$ (without pinches) has at least two $x$-letters:
$z=g_1x^{\epsilon}g_2x^{\eta}g_3\dots$, where $\epsilon, \eta\in \{1,-1\}$. Then  the only pinch in the product
$z^{-1}hz$ is $x^{-\epsilon}g_1^{-1}hg_1x^{\epsilon}$.

If $g_1^{-1}hg_1\in C$, then
 $\epsilon =-1$ and $d=xg_1^{-1}hg_1x^{-1}\in  D\backslash\{1\}$. Then we have
the pinch $x^{-\eta}g_2^{-1}dg_2x^{\eta}$, where the product $g_2^{-1}dg_2$ cannot belong
to $C$ by the assumption of the lemma. Thus, it is in $D$, $\eta=1$, and since $D$ is a malnormal subgroup of $A$, we should
have $g_2\in D$. But this gives the pinch $x^{-1}g_2x$ in the normal form of $z$, a contradiction.

If we had $g_1^{-1}hg_1\in D$, then the same argument would give a pinch $xg_2x^{-1}$ in the normal form of $z$. Therefore $z\in A$, and the lemma is proved.
\endproof

\begin{lemma} \label{noniso} (1) The group $G_k$ has an element of order $L$ with infinite centralizer
if the word $\alpha^k$ is accepted by the Turing machine $\mmm_0$. For all accepted
$\alpha^k$ the groups $G_k$ are isomorphic with the HNN extension $\overline G$ of $G$ with stable
letter $y$ and the additional relation $yW_{ac}y^{-1}=W_{ac}$.

(2) If the word $\alpha^k$ is not accepted by $\mmm_0$, then the centralizers of elements with
order $L$ in $G_k$ have order $L$.
\end{lemma}

\proof (1) By Lemma \ref{conj}, there is an element $g\in G$ such that $gW(k,k)g^{-1}=W_{ac}$
in $G$. So in $G_k$, we obtain the relation $x^{-1}gW(k,k)g^{-1}x= W(k,k)$, i.e.
$z^{-1}W(k,k)z=W(k,k)$ for $z=g^{-1}x$ and $yW_{ac}y^{-1}=W_{ac}$ for $y=xg^{-1}$. Here $W(k,k)$ has order $L$ by Lemma \ref{order} and $y, z$ have infinite order.
Furthermore, one can replace the generator $x$ by $y$ in the presentation of $G_k$ and
obtain the presentation of $\overline G$.

(2) Let an element $g$ has order $L$ in $G_k$, then it is a conjugate of an element of order
$L$ from $G$ (Theorem IV.2.4\cite{LS}). Therefore one may assume that $g\in G$, and its centralizer
$C_G(g)$ has order $L$ by Lemma \ref{cent} (1).
The same order
has the centralizer of $g$ in $G(k)$ by Lemma \ref{maln}, because the assumptions of that lemma
are guaranteed by Lemma \ref{cent} (1, 2).
\endproof

\medskip

{\bf Proof of Theorem \ref{isom}}.It follows from Lemma \ref{noniso} that the group $G_k$ is isomorphic to the group $\overline G$
(which is isomorphic to every $G_i$ with $\alpha^i$ accepted by the machine $\mmm_0$)
if and only if the word $\alpha^k$ is accepted by the Turing machine $M_0$. So the isomorphism
problem is not decidable in the set  of finitely presented groups $\{G_k\}_{k=1}^{\infty}$
if the language of accepted words of $\mmm_0$ is not recursive. Hence by Lemmas \ref{quadr}
and \ref{same},
the isomorphism problem is algorithmically undecidable in the class of finitely presented
groups with quadratic Dehn function. Theorem \ref{isom} is proved.

\section{Dehn functions of subgroups}\label{sub}

In this section, we assume that the language of accepted words of the
Turing machine ${\bf M}_0$ consists of all non-negative powers $\alpha^k$ in the one-letter alphabet $\{\alpha\}$, but we select ${\bf M}_0$ so that the time function $T_{{\bf M}_0}(n)$ of ${\bf M}_0$ grows fast. Given a recursive function $f(n)$, we can define
a symmetric Turing machine ${\bf M}_0$ with time function satisfying the inequalities
\begin{equation}\label{Tn}
T_{{\bf M}_0}(n)>f(n) \;\;
for\;\; n\ge 0.
\end{equation}
Here $T_{{\bf M}_0}(n)$ is the
time of the shortest ${\bf M}_0$-computation accepting the word $\alpha^n$.

 The  next S-machine ${\bf M}_1^+$ depends on ${\bf M}_0$ only, and we will assume that it  coincides with the machine ${\bf M_1}$ from \cite{OS06} provided by Lemma \ref{negat}.

The language
of accepted words for ${\bf M}_1^+$ is $\{\alpha^k\}_{k=0}^{\infty}$ by Lemma \ref{negat}, and we have the inequalities $T_{{\bf M}_0}(n)\le T_{{\bf M}_1^+}(n)$ for the time functions since a computation
of ${\bf M}_1^+$ accepting an input configuration has to simulate an accepting computation of ${\bf M}_0$ by Lemmas \ref{negat} and 4.1\cite{OS06}.

By definition, the machine ${\bf M}_1^-$ is a copy of
${\bf M}_1^+$, but the language of the accepted words
of ${\bf M}_1^-$ is $\{\alpha^k\}_{k=0}^{-\infty}$. We will assume that these two machines
have disjoint sets of rules, and the common state letters of them are only the letters of the start and the accept configurations.  The machine ${\bf M}_1$ is defined now as
the union of ${\bf M}_1^+$ and ${\bf M}_1^-$, where
every admissible words is admissible either for ${\bf M}_1^+$ or for ${\bf M}_1^-$. So an input (the accept) configuration of ${\bf M}_1$ is an input (resp., the accept) configuration of ${\bf M}_1^+$ and of ${\bf M}_1^-$.

\begin{lemma}\label{M1} The language of accepted words
of ${\bf M}_1$ is $\{\alpha^k\}_{k=-\infty}^{\infty}$,
and for every $n>0$, we have the inequality $T_{{\bf M_1}}(n)\ge f(n)- C(0)$, where $C(0)=T_{{\bf M_1}}(0)$.
\end{lemma}
\proof The first statement is obvious since ${\bf M}_1$
is the union of ${\bf M}_1^+$ and ${\bf M}_1^-$.

Let ${\cal C}: C_0\to\dots\to C_t$ be a shortest accepting computation of ${\bf M}_1$ starting
with an input configuration $C_0$ with a non-empty input word $\alpha^n$.
The computation $\cal C$ has an alternating factorization
${\cal C}={\cal C}_1\dots{\cal C}_s$, where every factor
belongs to either ${\bf M}_1^+$ or ${\bf M}_1^-$. Without loss of generality we assume that ${\cal C}_1: C_0\to\dots\to C_r$ is a
computation of ${\bf M}_1^+$.

If $n<0$, then the computation ${\cal C}_1$ cannot accept or end
with an input configuration by Lemma \ref{negat}. Therefore it cannot be followed by a computation ${\cal C}_2$ of ${\bf M}_1^-$, a contradiction.

If $n>0$ and $s=1$, then  ${\cal C}_1$ is an accepting computation
of ${\bf M}_1^+$, and so $$t\ge T_{{\bf M}_1^+}(n)\ge T_{{\bf M}_0}(n) \ge f(n)$$

If $s>1$, then ${\cal C}_1$ does not accept,
and so $C_r$ is a start configuration for both ${\bf M}_1^+$ and ${\bf M}_1^-$. Let $\alpha^k$ be the input word in $C_r$. Then $k\ge 0$ by Lemma \ref{negat} applied to ${\cal C}_1^{-1}$. If $k>0$, then ${\cal C}_2$
could not end working by the dual lemma applied to ${\cal C}_2$; hence $k=0$.

Note that since $k=0$, one can construct a computation $\cal D$ as computation ${\cal C}_1$ followed by the ${\bf M}_1^+$-computation of length $C(0)$ accepting the empty word, and $\cal D$ accepts $\alpha^n$.
Hence $T_{{\bf M_0}}(n)\le T_{{\bf M_1}}(n)\le r +C(0)$,
and so by (\ref{Tn}, we have
$$t\ge r \ge T_{{\bf M_0}}(n)-C(0)\ge f(n)-C(0),$$
which proves the lemma.
\endproof

The group $G$ is defined by $\bf M$ in Subsection \ref{md}.
In this section, we also consider a ``trimmed'' version of the machine $\bf M$. The set of rules of this machine $\bf\bar M$ is ${\bf \Theta}_3\cup {\bf \Theta}_4\cup{\bf \Theta}_5$, i.e. we now remove the sets ${\bf \Theta}_1$ and ${\bf \Theta}_2$, as well as the transition rule $\theta(23)$, from the definition of $\bf M$. The state
letters occurring in the removed rules only are removed
too. The words of the form $W(k,k')$ for arbitrary integers $k$ and $k'$
become the start configurations of the machine $\bf\bar M$.

The definitions of the machines ${\bf M}_2-{\bf M}_5, {\bf M}$
depend on ${\bf M}_1$ only, and ${\bf\bar M}$
accepts the same language $\{\alpha^k\}_{k=-\infty}^{\infty}$. Furthermore, every computation of each of these machines accepting a word $\alpha^k$ must simulate a work of the previous machine
accepting the same word, and so by Lemma \ref{M1}, for every non-negative $n$, we have inequalities for the time functions:
$$f(n)-C(0)\le T_{{\bf M}_1}(n)\le T_{{\bf M}_2}(n)\le  T_{{\bf M}_3}(n)$$

By the definition of the set ${\bf \Theta}_4$, an accepting computation for a word $W(k,k)$
is longer than the computation of ${\bf M}_3$ accepting the input $\alpha^n$. It follows that for every $n\ge 0$, we have
\begin{equation} \label{time}
f(n)-C(0)\le
T_{\bf\bar M}(n),
\end{equation}

Now we define the group $\bar M$  given by the generators and relations occurring in
formulas (\ref{rel111}) only (which correspond to the rules from ${\bf \Theta}_3\cup {\bf \Theta}_4\cup{\bf \Theta}_5$). The group $\bar G$ is obtained from $\bar M$ by imposing only one hub relation $W_{ac}^L=1$ from \ref{rel3}. In particular, the generators of the groups
$\bar M$ and $\bar G$ have no superscripts $(i)$, $i=1,\dots,L$.

\begin{lemma}\label{inj} The canonical homomorphisms
$\bar M\to M$ and $\bar G\to G$ are injective. So one may
identify $\bar M$ and $\bar G$ with the subgroups of the groups $M$ and $G$, respectively.
\end{lemma}

\proof We should prove that if a word $w$ in the generators of $\bar G$ is equal to 1 in the group $M$ (in $G$), then it represents $1$ in $\bar M$ (in $\bar G$, resp.)

Let $\Delta$ be a minimal diagram over $G$ with boundary label $w$. If $\Delta$ has no disks, then by Lemma \ref{NoAnnul}, every maximal $\theta$-band of $\Delta$
ends on the boundary $\partial\Delta$, and so the one-letter history of it is a history of $\bf\bar M$ since $w$
is a word in the generators of $\bar M$. It follows that
$\Delta$ is a diagram over $\bar M$ and $w=1$ in $\bar M$, as required. Thus, one may assume that $\Delta$ has
at least one disk and induct on the number of disks $l$
with base $l=0$.

If $l>0$, then Lemma \ref{extdisk} provides us with a disk $\Pi$ and a $\tt$-band $\cal B$ connecting this disk
with the boundary $\partial\Delta$. Since by Lemma \ref{withd} (2), $\Delta$ contains no $\theta$-annuli, every $\theta$-band crossing $\cal B$ ends on $\partial\Delta$
and it is a diagram over $\bar M$. It follows that the
$\tt$-letter labelling the common edge of $\cal B$ and $\partial\Pi$
has no superscripts. Hence no letter of the accessible
boundary label of $\Pi$ has a superscript; this label has the form $W^L$.

By Lemma \ref{dr}, we have a computation $\cal C$ of $\bf M$,
connecting the word $W$ either with $W_{st}$ or with $W_{ac}$. If the history of this computation has no rules $\theta(23)^{-1}$, then it is a computation of $\bf\bar M$,
and by Lemmas \ref{I6A6} (1) and \ref{dr}, the disk $\Pi$ can be filled in with the cells corresponding to the relations of the group $\bar G$
(including the hub relation $W_{ac}^L$).

If $\cal C$ has a rule $\theta(23)^{-1}$, then this
rule is applied to a configuration $W(k,k')$. Here
$k=k'$ since the word $W$, and so the word $W(k,k')$,
is accessible.  However for every integer $k$, the configuration $W(k,k)$
is accepted by the machine $\bf\bar M$ since the language of
$\bf\bar M$-accepted input words is $\{\alpha^k\}_{k=-\infty}^{\infty}$. Therefore by Lemmas \ref{I6A6} (1) and  \ref{dr}, the
subdisk with boundary label $W(k,k)^L$ can be filled in
with cells corresponding to the presentation of $\bar G$. The same is true for the whole disk $\Pi$.
 Then $\Pi$ can be removed from $\Delta$ along with the
 $\tt$-band $\cal B$, and the boundary of the remaining
 part $\Delta'$ of $\Delta$ is again labelled over $\bar G$. Since the number of disks in $\Delta'$ is $l-1$,
 the lemma is proved by induction.
 \endproof

 Since every word $W(k,k)$ is accepted by the machine $\bf\bar M$, Lemma \ref{dr} gives us a {\it disk diagram} $\Delta$ with boundary label  $W(k,k)^L$ built of one hub and $L$ trapezia corresponding to a reduced accepting computation for $W(k,k)$. The boundary of the hub in $\Delta$ is labeled by $W_{ac}^L$, since $\bar G$ has only one hub relation.

One can prove that the disk diagram with boundary label $W(k,k)^L$ is minimal, but we need an estimate from below for the area of the word $W(k,k)^L$ with respect to
the {\it finite} presentation of $\bar G$ (which contains a hub, but no other disks).
The following statement is Lemma 10.2 of \cite{O18}. (Although the machine is different in \cite{O18}, the proof of Lemma 10.2 works
for $\bar G$ without any changes.)

\begin{lemma} \label{2ar} The area of a disk diagram $\Delta$
with boundary label $W(k,k)^L$
does not exceed two areas of the disk word $W(k,k)^L$ with respect to
the finite presentation of $\bar G$.
\end{lemma} $\Box$

 \begin{lemma}\label{Wkk} The area of the word $W(k,k)^L$ with respect
 to the finite presentation of $\bar G$
 is at least $L(f(k)-C(0))$.
 \end{lemma}

 \proof Consider the disk diagram $\Delta$
 with boundary label $W(k,k)^L$ provided by Lemma \ref{dr}.
 $\Delta$ contains $L$ trapezia corresponding to an
 accepting computation of the machine $\bf \bar M$
 starting with the configuration $W(k,k)$. By Lemma
 \ref{simul}, the height of each trapezium $\Gamma$ is at least $T_{\bf\bar M}(k)$, which is greater than $f(k)-C(0)$ by inequality (\ref{time}). Hence $\Gamma$
 contains at least $N(f(k)-C(0))$ $(\theta,q)$-cells, and
 therefore $\Delta$ has at least $NL(f(k)-C(0))$ cells. By Lemma \ref{2ar}, the area of the word $W(k,k)^L$ is at least $NL(f(k)-C(0))/2$
 which proves the lemma.
 \endproof

 {\bf Proof of Theorem \ref{subgr}}. Note that the length of the word $W(k,k)$ is a linear function of $k$. Therefore to bound the Dehn function
 of the subgroup $\bar G$ from below by $f(n)$ (up to equivalence), it suffices to obtain the inequalities
 $\area_{\bar G}(W(k,k))> f(k)- C(0)$ for every $k\ge 1$.
 Indeed, these inequalities follow from Lemmas  \ref{Wkk}.

 Since the group $H=\bar G$ embeds in $G$ by Lemma \ref{inj}, Theorem \ref{subgr} is proved, because by Lemma \ref{quadr}, the group $G$ defined by $\bf M$ in Subsection \ref{md}
has quadratic Dehn function.
 $\Box$

 \begin{rk} Since the  word $W(k,k)$ is accepted by the machine $\bf\bar M$ for every integer $k$, it follows from Lemmas \ref{conj},
 \ref{order}, \ref{finor}, \ref{enough}, and \ref{coninf} that the
 conjugacy problem is decidable for both groups $G$ and $\bar G=H$
 constructed in this section.
 \end{rk}

 \section{Conjugacy in the group G}\label{Co}

 In this section and in the next one, the construction of the machine $\bf M$
 can be based on the machine ${\bf M}_1$ provided either by Lemma \ref{negat} or by Lemma \ref{S}.

\begin{lemma}\label{malo} Let $\Gamma$ be a reduced diagram over $M$
with boundary $\bf xy$. Suppose there are  neither $\theta$-bands nor
$q$-bands starting and ending on $\bf y$. Then

(1) $|{\bf y}|_{\theta}\le |{\bf x}|_{\theta}$ and $|{\bf y}|_{q}\le |{
\bf x}|_{q}$;

(2) If $\Gamma$ is a subdiagram of a diagram $\Delta$ and $\bf y =y_1y_2y_3$, where $\bf y_1$ and $\bf y_3$ are sides of $q$-bands,
$\bf y_2$ is a side of a $\theta$-band or a subpath of the boundary of a disk $\Pi$, then $||{\bf y_2}||$  is bounded by a quadratic function of $||\bf x||$. The perimeter $||\partial\Pi||$ is also bounded by a quadratic function of $||\bf x||$  if $\bf y_2$
contains at least two $\tt$-letters.
\end{lemma}

\proof (1) Since every maximal $\theta$- or $q$-band of $\Gamma$ starting on $\bf y$ has to end on $\bf x$, the inequalities follow.

(2) It follows from (1) that the numbers of maximal $\theta$- and
 $q$-bands of $\Gamma$ are bounded by a linear function of $||{\bf x}||$. Therefore the number of $(\theta,q)$-cells of $\Gamma$ is bounded by a quadratic function by Lemma \ref{NoAnnul}. Note that the $Y$-lengths of the sides $\bf y_1$, $\bf y_2$ of $q$-bands are linearly bounded by their $\theta$-lengths. Since
 every maximal $Y$-band starting on ${\bf y}_2$ has to end on
 a $(\theta, q)$-cell or on $\bf y_1$, or on $\bf y_2$, or on $\bf x$,
 we have $|{\bf y_2}|_Y$ and $||{\bf y}_2||$ bounded by a quadratic function of $||{\bf x}||$. If $\bf y_2$ is a subpath of $\Pi$ having at least two $\tt$-edges, we have $||\partial\Pi||< L||{\bf y}_2||$, and the statement (2) follows.
\endproof

 \begin{lemma} \label{simpl} There is an algorithm replacing given word $W$
with a conjugate in $G$ word $W'$ having the property that
there is no minimal van Kampen diagram $\Gamma$ with boundary  $\bf pq$, where
$Lab (\bf p)$ is a subword of a cyclic permutation of the word $W'$
and

(a) ${\bf q}= {\bf y_1y_2y_3}$, where $\bf y_1$ and $\bf y_3$ are
sides of $\tt$-spokes connecting the subpath $\bf y_2$ of the boundary of a disk $\Pi$ with $\bf p$,
provided that there are no other disks in $\Gamma$  and $\Pi$ is connected with $\bf p$ by $s>L/2 $ $\tt$-bands or

(b) ${\bf q}$ is a side of a rim $q$-band $\cal C$ of $\Gamma$ starting and ending
on $\bf p$, or

(c) ${\bf q}$ is a side of a rim $\theta$-band $\cal T$ of $\Gamma$ starting and ending on $\bf p$.

The algorithms  providing properties (a) and (b) do not increase
the $q$-length of the word $W$, and the algorithm providing property (c)
 increases neither $q$-length nor $\theta$-length.
\end{lemma}

% This is a LaTeX picture output by TeXCAD.
% File name: [lemma62.pic].
% Version of TeXCAD: 4.3
% Reference / build: 30-Jun-2012 (rev. 105)
% For new versions, check: http://texcad.sf.net/
% Options on the following lines.
%\grade{\on}
%\emlines{\off}
%\epic{\off}
%\beziermacro{\on}
%\reduce{\on}
%\snapping{\off}
%\pvinsert{% Your \input, \def, etc. here}
%\quality{8.000}
%\graddiff{0.005}
%\snapasp{1}
%\zoom{4.0000}
\unitlength 1mm % = 2.845pt
\linethickness{0.4pt}
\ifx\plotpoint\undefined\newsavebox{\plotpoint}\fi % GNUPLOT compatibility
\begin{picture}(127.25,53.75)(10,70)
\thicklines
\put(19,110){\line(1,0){42}}
\put(40.75,92.625){\oval(17,8.75)[]}
%\emline(23.5,109.75)(34,95.5)
\multiput(23.5,109.75)(.0336538462,-.0456730769){312}{\line(0,-1){.0456730769}}
%\end
%\emline(25.5,109.75)(34.5,96.75)
\multiput(25.5,109.75)(.0337078652,-.0486891386){267}{\line(0,-1){.0486891386}}
%\end
%\emline(33,110)(36.5,97)
\multiput(33,110)(.033653846,-.125){104}{\line(0,-1){.125}}
%\end
%\emline(34.5,109.75)(37.5,97.25)
\multiput(34.5,109.75)(.03370787,-.14044944){89}{\line(0,-1){.14044944}}
%\end
%\emline(47.5,96.5)(57,110.25)
\multiput(47.5,96.5)(.0336879433,.0487588652){282}{\line(0,1){.0487588652}}
%\end
%\emline(58.25,109.75)(48.25,95.25)
\multiput(58.25,109.75)(-.0336700337,-.0488215488){297}{\line(0,-1){.0488215488}}
%\end
\put(31.5,113.75){$Lab({\bf p})\equiv W'$}
\put(39.75,106){$\Gamma$}
\put(23.5,103){$\bf y_1$}
\put(55.75,102.5){$\bf y_2$}
\put(40.5,99.25){$\bf y_3$}
\put(39.75,92.25){$\Pi$}
%\emline(91.5,99.5)(100.25,105.5)
\multiput(91.5,99.5)(.049157303,.033707865){178}{\line(1,0){.049157303}}
%\end
\put(100.25,105.5){\line(1,0){20.75}}
\put(121,105.5){\line(1,-1){6.25}}
\put(94.75,101){\line(1,-1){5.5}}
\put(100.25,95.5){\line(1,0){18.5}}
%\emline(118.75,95.5)(125.25,101.25)
\multiput(118.75,95.5)(.038011696,.033625731){171}{\line(1,0){.038011696}}
%\end
\put(92.75,99.5){\line(1,-1){6.5}}
\put(99.25,93){\line(1,0){21.25}}
%\emline(120.5,93)(127,99)
\multiput(120.5,93)(.036516854,.033707865){178}{\line(1,0){.036516854}}
%\end
%\emline(96,99.25)(95.25,97)
\multiput(96,99.25)(-.0326087,-.0978261){23}{\line(0,-1){.0978261}}
%\end
%\emline(96.5,95.5)(98.5,97)
\multiput(96.5,95.5)(.04444444,.03333333){45}{\line(1,0){.04444444}}
%\end
%\emline(98.5,97)(98.25,97.25)
\multiput(98.5,97)(-.03125,.03125){8}{\line(0,1){.03125}}
%\end
%\emline(100.25,95.75)(99.25,93)
\multiput(100.25,95.75)(-.0333333,-.0916667){30}{\line(0,-1){.0916667}}
%\end
%\emline(102.5,95.5)(103.5,93.5)
\multiput(102.5,95.5)(.0333333,-.0666667){30}{\line(0,-1){.0666667}}
%\end
%\emline(103.75,93.5)(103.5,94)
\multiput(103.75,93.5)(-.03125,.0625){8}{\line(0,1){.0625}}
%\end
\put(106.5,95){\line(0,-1){2.25}}
%\emline(110,95.25)(111.25,93.25)
\multiput(110,95.25)(.03289474,-.05263158){38}{\line(0,-1){.05263158}}
%\end
%\emline(114.5,95.5)(114,93)
\multiput(114.5,95.5)(-.0333333,-.1666667){15}{\line(0,-1){.1666667}}
%\end
%\emline(118.5,95.5)(116.75,93)
\multiput(118.5,95.5)(-.03365385,-.04807692){52}{\line(0,-1){.04807692}}
%\end
\put(116.75,93){\line(0,1){0}}
\put(120.25,96.5){\line(0,-1){3.25}}
%\emline(122.75,98.25)(123.75,96.5)
\multiput(122.75,98.25)(.0333333,-.0583333){30}{\line(0,-1){.0583333}}
%\end
\put(104.5,89.5){$\bf q$}
\put(110.75,107.75){$\bf p$}
\put(103,98.25){$\cal C$ or $\cal T$}
\put(21,81.5){Diagram $\Gamma$ in Lemma \ref{simpl} (a)}
\put(85.75,80.75){Diagram $\Gamma$ in Lemma \ref{simpl} (b,c)}
\end{picture}

\proof (a) Assume that $\Pi$ is connected  with $\bf p$ in $\Gamma$ by $r> L/2$ consecutive $\tt$-spokes ${\cal C}_1,\dots, {\cal C}_r$.
By Lemma \ref{withd} (1), $\Gamma$ contains no $\theta$-bands connecting ${\cal C}_1$ and ${\cal C}_r$. Lemma \ref{malo}  gives a linear bound (in terms of $||W||$) for the lengths of the spokes and quadratic upper bound for
 the perimeter of the disk $\Pi$ . So there is a subdiagram $\Gamma'$
 containing $\Pi$ and having the boundary $\bf pq'$ with $|{\bf q'}|_q<|{\bf q}|_q$ and bounded $||{\bf q'}||$.

 Replacing the subword $Lab ({\bf p})$ with $Lab ({\bf q'})^{-1})$, one obtains a conjugate in $G$ word of smaller $q$-length, where the  length
 of modified word is quadratically bounded in terms of $||W||$, and the search of it has effectively bounded time by Lemma \ref{quadr}. This gives the required algorithm.

 (b) Consider such a diagram $\Gamma$ if any exists. We may assume that $\Gamma$ contains no disks. Indeed,  the $q$-band $\cal C$ has no side $q$-edges, and so by Lemma \ref{extdisk} (1), the existence of a disk should imply the existence of a sudiagram already eliminated in item (a), a contradiction.  Then by Lemma \ref{NoAnnul}, every maximal $\theta$-band starting from $\cal C$ in $\Gamma$ ends on $\bf p$. Hence the length
 of $\cal C$ is less than $||\bf p||$. The side label of $\cal C$ is equal in $M$ to $Lab ({\bf p})$ but has $q$-length $0$.
 Hence the subword $Lab ({\bf p})$  can be replaced with an equal in $M$ word
 $Z$ having smaller $q$-length and $||Z||$ is linearly bounded in terms of $||U||$. Since the word problem is decidable in $G$, one can
 efficiently execute such a replacement. This procedure provides us with the desired algorithm.

 (c) Let $\Gamma$ be the diagram from the formulation of item (c). Now we may assume that the word $W$ has property (a). We also may assume that
 $\Gamma$ has no maximal $\theta$-bands except for $\cal T$ since otherwise $\Gamma$ should contain a proper subdiagram with the same property.
 Assume that $\Gamma$ has a disk. Then let $\Pi$ be a disk provided by Lemma \ref{extdisk} (1). By item (a), at least two $\tt$-spokes $\cal C$ and $\cal C'$ start on $\Pi$ and end on $\cal T$ since $L-5> L/2+1$. The uniqueness of the maximal $\theta$-band $\cal T$ in $\Gamma$ implies that $\cal C$ and $\cal C'$ have length $0$. Then by Lemma \ref{moving} (2), one can decrease the number
 of disks in the diagram $\Gamma$, replacing it with a diagram $\Gamma'$
 satisfying the same assumptions and having $Lab({\bf p'})\equiv Lab({\bf p})$. The induction on the number of disks allows us to assume that $\Gamma$ is a diskless diagram, and so it coincides with $\cal T$.
 Replacing $\bf p$ with ${\bf q}^{-1}$, one decreases the $\theta$-length of the boundary label
 preserving the $q$-length of it. Since the length of $\cal T$ does not exceed $||W||$  one can remove $\cal T$ effectively.
 \endproof

 \begin{rk} \label{awa}Let us call the word $W'$ from Lemma \ref{simpl} an {\it adapted} word.

 One can change the formulation of Lemma \ref{simpl} by replacing
 the minimal diagram $\Gamma$ over $G$ with a reduced diagram over
 the group $M$ and removing Property (a). The statement remains true
 (the proof is a simplified proof of Lemma \ref{simpl}).
 We will call the words $W'$ obtained from $W$ according to these
 weaker version of Lemma \ref{simpl} a {weakly} adapted word.
\end{rk}

The following statement will be used in the next section.

 \begin{lemma} \label{coninf} There is an algorithm that decides whether two words
 $U$ and $V$ representing elements of infinite order in the group $G$
 are conjugate in $G$ or not.
 \end{lemma}

 \proof Below the proof is subdivided in a few steps.

 {\bf 0.}  Lemma \ref{simpl} allows us to
 assume that the words $U$ and $V$ are adapted.
 By Schupp Lemma, $U$ and $V$ are conjugate in $G$ if and only if there is an annular diagram $\Delta$ over $G$ with boundary components $\bf p$ and $\bf q$ labelled by $U$ and $V$, respectively. If there is a recursive function $f$ such that one always can choose $\Delta$ so that $\bf p$ and $\bf q$ can be connected
 by a path $\bf x$ of length $\le f(||U||+||V||)$, then the cut along $\bf x$
 replaces $\Delta$ with a van Kampen diagram of bounded perimeter, and so the conjugacy
 problem is reduced to the word problem. Since the word problem is decidable in a group with quadratic Dehn function (see Lemma \ref{quadr}), our goal is to find
 such a ''short cut'' under the assumption that $U$ and $V$ are conjugate
 in $G$.

 {\bf 1.} Let $\Delta$ be a minimal annular diagram whose boundary contours $\bf p$ and $\bf q$ are labeled by $U$ and $V$, resp.
 The number $r$ of disks in $\Delta$ cannot exceed
 $s=|U|_q+|V|_q$ for the following reason. Let $\Pi$ be the disk
 provided by Lemma \ref{extdisk} (2). Since the words $U$ and $V$ are adapted, $\Pi$ is connected by spokes with both $\bf p$ and $\bf q$.
 The complement to the union of $\Pi$ and the diskless subdiagrams between $\Pi$ and the boundary components (bounded by the spokes at $\Pi$) is a van Kampen
 diagram $\Delta'$ with at most  $s - (L-3) + 3< s$ $\tt$-edges on the boundary.
 Now Lemma \ref{mnogospits} bounds the number of disks in $\Delta'$.

 {\bf 2.} If $\Delta$ has a disk, then Lemma \ref{extdisk} (2) gives
 a disk $\Pi$ connected with $\bf p$ (or with $\bf q$) by at least
 two $\tt$-spokes. Assuming that this two $\tt$-bands $\cal C$ and $\cal C'$ are consecutive, we consider the subdiagram $\Gamma$ over the group $M$ bounded  by $\cal C$, $\cal C'$, $\partial\Pi$ and $\bf p$.

 If $\Gamma$ contains no $\theta$-bands connecting $\cal C$ and $\cal C'$, then  Lemma \ref{malo}, gives a linear bound (in terms of $||U||$) for the perimeter of the disk $\Pi$ and the length of $\cal C$.

 Making a cut along the boundary of $\cal C$ and around $\Pi$, one can
 remove the disk $\Pi$ and obtain an annular diagram $\Delta'$ with fewer disks, where the boundary label $U$ is replaced with equal in the group $G$ word $U'$, whose length is quadratically bounded in terms of $||U||$.

 If $\Gamma$ has a $\theta$-band connecting $\cal C$ and $\cal C'$, then
 such a $\theta$-band closest to $\Pi$ has to share a side with
 $\partial\Pi$. This $\theta$-band $\cal T$ and $\Pi$ form a diagram $E$
 satisfying the assumption of Lemma \ref{moving} (1). Therefore $\Pi$ and the subdiagram $E$ can be replaced in $\Delta$ with a disk $\Pi'$ and a diagram $E'$. This surgery removes the $\theta$-band $\cal T$ from $\Gamma$ and shortens the connecting $\tt$-bands $\cal C$ and $\cal C'$
 in the obtained annular diagram $\Delta'$.
 (We do not care of the minimality of the entire $\Delta'$.) Then we can continue moving the disk closer to $\bf p$ until we obtain a subdiagram, where no $\theta$-band connects $\cal C$ and
 $\cal C'$.

 Thus, if $\Delta$ has a disk, then there is a minimal annular  diagram $\Delta'$ with fewer disks and boundary label $U'$ and $V'$ equal to $U$ and $V$ in $G$, where $||U'||+ ||V'||$ is effectively bounded in terms
 of $||U||+||V||$. Since the number of disks in $\Delta$ does not
 exceed $||U||+||V||$, the iteration of this argument, provides us
 with a diskless annular diagram $\bar\Delta$ and with boundary
 labels $\bar U$ and $\bar V$ equal to $U$ and $V$, resp., in $G$.
 Hence an effective exhaustive search gives a finite set $S=S(U,V)$ of pairs $(U_i, V_i)$ such that $U$ and $V$ are conjugate in $G$ if and only if for some $i$, the words $U_i$ and $V_i$ are conjugate in the group $M$. Moreover, by Remark \ref{awa}, the words $U_i$ and $V_i$
 can be assumed weakly adapted. So, keeping the same notation,  we may assume that the annular diagram $\Delta$ contains no disks.

 $\bf 3.$  Since we may assume that   the words $U$ and $V$ are weakly adapted,
 it remains to consider three options: (a) neither $U$ nor $V$ have
 $\theta$-letters; borrowing the term from \cite{OS04}, the corresponding annular diagram $\Delta$ will be called a {\it ring}; (b) there are no $q$-edges
 in the boundary of $\Delta$ and every maximal $\theta$-band  connects $\bf p$ and $\bf q$; such a diagram is a {\it roll}; (c) there
 are $q$-letters and there are  $\theta$-letters in both $U$ and $V$,
 and every maximal $q$- or $\theta$-band of $\Delta$ connects $\bf p$
 and $\bf q$; $\Delta$ is a {\it spiral}.

 {\bf 4.} Assume that $\Delta$ is a ring. Then, by Lemma \ref{NoAnnul},
 the annular diagram $\Delta$ is built of $\theta$-annuli surrounding the hole of $\Delta$. Different $\theta$-annuli cannot copy each other, since otherwise one could remove some $\theta$-annuli from $\Delta$.

 If $\Delta$ has no $(\theta,q)$-cells, then by Lemma
 \ref{NoAnnul}, every maximal $Y$-band connects $\bf p$ and $\bf q$,
 and so all the $\theta$-annuli have the same length $|U|_a$, and the number of different $\theta$-annuli of this length is effectively bounded. Therefore there is a path $\bf x$ of bounded length
 connecting $\bf p$ and $\bf q$, as desired.
 Therefore, we may assume that there are $(\theta,q)$-cells
 in $\Delta$, and  we have $s\ge 1$ maximal $q$-bands ${\cal C}_1$,..., ${\cal C}_s$,
 each of them has the same length $h$ and connects $\bf p$ and $\bf q$. Let the $\theta$-annuli of $\Delta$ have boundary components with lengths $l_0$ and $l_1$, $l_1$ and $l_2$,..., $l_{h-1}$ and $l_h$.

 If $\max_{i=0}^h l_i \le c_4\max (||U||, ||V||)$, then there is an effective upper bound for $h$ since the $\theta$-annuli cannot copy each other. So proving by contradiction, we assume that

 \begin{equation}\label{c5L}
 max_{i=0}^h l_i > c_4\max (||U||, ||V||)
 \end{equation}

 We consider the ''power'' $\Delta^{L}$ of $\Delta$. This reduced annular diagram has boundary labels $U'=U^{L}$ and $V'=V^{L}$, respectively, and maximal $q$-bands ${\cal D}_1,\dots {\cal D}_{sL}$.
 ($\Delta^{L}$ covers $\Delta$ with multiplicity $L$).

 Let us cut $\Delta^{L}$ along a side of ${\cal D}_{sL}$ and attach
 a copy ${\cal D}_0$ of the $q$-band ${\cal D}_{sL}$ obtaining a
 trapezium $\Gamma$ of height $h$. The base of $\Gamma$ has form
 $xv_1xv_2x\dots x v_rx$, where the letter $x$ does not occur in the subwords $b_1,\dots, b_r$ and $r$ is a multiple of $L$. So $\Gamma$ is covered by $r$ trapezia $\Gamma_i$ with bases $xb_ix$.

 By inequality (\ref{c5L}) and Lemma \ref{xvx}, the base of $\Gamma$
 is a power of a cyclic permutation of the standard base.
Replacing $U$ with a cyclic permutation, we may assume that the base of $\Gamma_i$ is standard, it is equal $tv$, where $x\equiv t$, $v\equiv v_1\equiv\dots\equiv v_r$. Moreover, by Lemmas \ref{simul} (1) and \ref{narrow}, for the top/bottom labels $W_i$/$W'_i$ of every $\Gamma_i$, the words $W_i^{\emptyset}$/$(W'_i)^{\emptyset}$ are accessible words.
Furthermore by Remark \ref{les} and Lemma \ref{resto}, all the words $W_1,\dots, W_r$ are equal up to the superscripts. It follows
that the word $U^L$ is permissible, and it is a power of a disk word, because $r$ is divisible by $L$. So $U^L=1$ in the group $G$ by  Lemma \ref{dw}, contrary to the assumption that $U$ has infinite order.

{\bf 5.} Assume now that $\Delta$ is a roll having no $(\theta,q)$-cells.
If it has a $Y$-annulus $\cal A$ surrounding the hole of $\Delta$, then it follows from the form of $(\theta,a)$-relations, that the inner and the outer boundary components have the same boundary label. Then one can just remove $\cal A$ from $\Delta$ identifying these two sides of $\cal A$. Therefore we may assume that there are no such $Y$-annuli in $\Delta$, and every maximal $Y$-band starts or ends on $\bf p$ or $\bf q$ by
Lemma \ref{NoAnnul}. It follows that the number of maximal $Y$-bands cannot exceed
the sum $||U||+||V||$. To obtain an upper bound for the length of a connected path $\bf x$ from item $\bf 1$, it suffices to bound the number of $(\theta,a)$-cells in
$\Delta$, that is to bound from above the length of every maximal $Y$-band.

If an $Y$-band $\cal A$ starts and ends on $\bf p$ (or $\bf q$), then its length does not
exceed $||U||$ since by Lemma \ref{NoAnnul}, every maximal $\theta$-band crosses $\cal A$ at most once and starts or ends on $\bf p$. Let $\cal A$ connects $\bf p$ and $\bf q$. To bound the length of $\cal A$ it suffices to bound the number
of $(\theta,a)$-cells belonging to the intersection of $\cal A$ with arbitrary maximal $\theta$-band
$\cal T$, because the number of maximal $\theta$-bands (connecting $\bf p$ and $\bf q$) is at most $||U||+||V||$.

 By Lemma \ref{NoAnnul}, a  subband of $\cal A$ cannot cross twice a subband of $\cal T$ in a van Kampen subdiagram of $\Delta$. Therefore after the $Y$-band  $\cal A$ crosses $\cal T$ at some cell $\pi$,
 it has to cross every other maximal $\theta$-band of $\Delta$ before  $\cal A$ crosses $\cal T$ again at some cell $\pi'$. So we have a convolution,
 i.e., the subband of $\cal A$ of length at most $||U||+||V||$ between $\pi$ and $\pi'$. Hence it suffices to bound from above the number of such convolutions
 in $\cal A$.

 The subband ${\cal T'}$ of $\cal T$ between $\pi$ and $\pi'$ can be crossed by
 another maximal $Y$-band $\cal A'$ at most once (and  $\cal A'$ has to connect
 $\bf p$ and $\bf q$). Therefore the length of ${\cal T'}$ is at most $||U||+||V||$.
 So a side of the convolution and a side of ${\cal T'}$
  form a loop $\bf z$ of  length $O(||U||+||V||)$ surrounding the hole of $\Delta$.

  If $\bf z$ surrounds another loop $\bf z'$ of this type with $Lab ({\bf z'})\equiv Lab ({\bf z})$, then one can remove the diagram between $\bf z$ and $\bf z'$ and identify these
  two loops decreasing the number of cells in the annular diagram.
  Since the lengths of loops of this type are bounded, the number of
  them is effectively bounded as well. It follows that we have an
  effective upper bound for the lengths of $Y$-bands as desired.

 {\bf 6.} Assume now that $\Delta$ is a roll containing $(\theta,q)$-cells. Then by Lemma \ref{NoAnnul},
 every maximal $q$-band of $\Delta$ is a $q$-annulus surrounding the hole of $\Delta$ since a roll has no $q$-edges in the boundary.
 By the same lemma, every maximal $\theta$-band crossing a $q$-annulus
 $\cal C$ ends on $\bf p$ and $\bf q$ and cannot intersect $\cal C$ twice. Therefore the length of arbitrary $q$-annulus does not
 exceed $\min (|U|_{\theta}, |V|_{\theta})$. This observation effectively
 bounds the number of $q$-annuli in a minimal annular diagram $\Delta$,
 because two different annuli cannot copy each other. (One could remove
 the diagram between them and identify such annuli, contrary to
 the minimality of $\Delta$.)

If ${\bf p}_1,  {\bf q}_1,\dots,  {\bf p}_k,  {\bf q}_k$ are pairs of boundary components of all the $q$-annuli counting from $\bf p$ to $\bf q$, then the annular diagrams between ${\bf q_0}=\bf p$ and ${\bf p}_1$,
 ${\bf q}_1$ and ${\bf p}_2$,\dots, ${\bf q}_k$ and ${\bf p}_{k+1}={\bf q}$ contain no
 $(\theta,q)$-cells. Therefore there is a roll having boundary
 labels $U$ and $V$ if and only if there are at most $k$ different $q$-annuli (where $k$ and the lengths of the $q$- annuli are effectively bounded) and $\le k+1$ rolls without $(\theta, q)$-cells between ${\bf q}_i$ and
 ${\bf p}_{i+1}$ ($i=0,\dots,k$, and the length of all ${\bf p}_j$-s, ${\bf q}_j$-s are efficiently bounded). So the case of rolls is efficiently
 reduced to the special case considered in item $\bf 5$.

 {\bf 7.} It remains to assume that $\Delta$ is a spiral. It follows from Lemma \ref{NoAnnul}, that every $q$-edge of $\Delta$ belongs to one
 of the (clockwise) consecutive maximal $q$-bands ${\cal C}_1,\dots, {\cal C}_s$ connecting $\bf p$ and $\bf q$. Clearly we have $s\le \min(|U|_q, |V|_q)$. Similarly, all maximal $\theta$-bands ${\cal T}_1,\dots, {\cal T}_r$ connect
 $\bf p$ and $\bf q$, and  $r\le \min(|U|_{\theta}, |V|_{\theta})$.

 Assume that going from $\bf p$ to $\bf q$ a $\theta$-band ${\cal T = T}_j$ crosses a $q$-band ${\cal C = C}_i$ clockwise. Then Lemma \ref{NoAnnul} implies, that if after this intersection, $\cal T$
 crosses a $q$-band again, then this next intersection is with ${\cal C}_{i+1}$ (the indices taken modulo $s$) and ${\cal C}_{i+1}$
 is crossed clockwise too. So if the number of intersections
 of $\cal T$ with $q$-bands is greater than $s$, it has to cross
 one of the band at least twice.

 Our nearest goal is to bound from above the lengths of the $q$-bands. And so we will assume that every maximal $\theta$-band
 contains greater than $s$ different $(\theta,q)$-cells. (If there is $\cal T$ crossing
 $q$-bands at most $s$ times, then the same has to be true for every
 $q$-band.)

 The spiral structure of $\Delta$ implies the dual property:
 If a maximal $q$-band $\cal C$ (directed from $\bf p$ to $\bf q$) crosses some
 ${\cal T}_j$, then the next intersection of $\cal C$ (if any) is the intersection with ${\cal T}_{j+1}$ (the indices taken modulo $r$).
 It follows that the history of any $q$-band is periodic with
 a period $H$ of length $r$.

 Consider now a van Kampen subdiagram $\Gamma$ bounded by two
 $q$-bands $\cal C$ and $\cal C'$ and parts of $\bf p$ and $\bf q$,
 such that $\Gamma$ has no other $q$-bands between $\cal C$ and $\cal C'$.  Let $E$ be the maximal trapezium (possibly empty) of $\Gamma$, bounded by
 subbands $\cal D$ and $\cal D'$ of $\cal C$ and $\cal C'$, and by $\theta$-bands $\cal S$ and $\cal S'$ connecting $\cal C$ and $\cal C'$
 in $\Gamma$. Thus, the compliments of $E$ in $\Gamma$ have no
 $\theta$-bands connecting $\cal C$ and $\cal C'$, and Lemma \ref{malo} (2) gives
 a quadratic bound for the lengths of $\cal S$ and $\cal S'$

  By Lemma \ref{simul} (1), the trapezium $E$ corresponds to an eligible computation with periodic history $H$. Therefore Lemma \ref{wi} gives a linear upper bound of the lengths
  of arbitrary $\theta$-bands of $E$ in terms of the lengths of the top and the bottom of $E$, and the length of the period $H$ of the history. So the lengths of such $\theta$-bands are quadratically bounded in terms of $||U||+||V||$.

  % This is a LaTeX picture output by TeXCAD.
% File name: [Ez.pic].
% Version of TeXCAD: 4.3
% Reference / build: 30-Jun-2012 (rev. 105)
% For new versions, check: http://texcad.sf.net/
% Options on the following lines.
%\grade{\on}
%\emlines{\off}
%\epic{\off}
%\beziermacro{\on}
%\reduce{\on}
%\snapping{\off}
%\pvinsert{% Your \input, \def, etc. here}
%\quality{8.000}
%\graddiff{0.005}
%\snapasp{1}
%\zoom{4.0000}
\unitlength 1mm % = 2.845pt
\linethickness{0.4pt}
\ifx\plotpoint\undefined\newsavebox{\plotpoint}\fi % GNUPLOT compatibility
\begin{picture}(145,80)(0,70)
\thicklines
%\emline(144.75,147)(145,146.75)
\multiput(144.75,147)(.03125,-.03125){8}{\line(0,-1){.03125}}
%\end
%\emline(12,136.25)(21.75,140.25)
\multiput(12,136.25)(.081932773,.033613445){119}{\line(1,0){.081932773}}
%\end
\put(21.75,140.25){\line(1,0){29}}
%\emline(50.75,140.25)(58,136.25)
\multiput(50.75,140.25)(.06092437,-.033613445){119}{\line(1,0){.06092437}}
%\end
\put(16.75,104){\line(5,3){6.25}}
\put(23,107.75){\line(1,0){25.75}}
%\emline(48.75,107.75)(53.25,104.5)
\multiput(48.75,107.75)(.046391753,-.033505155){97}{\line(1,0){.046391753}}
%\end
\put(25.75,140.25){\line(0,-1){32.5}}
\put(27.5,140.5){\line(0,-1){32.5}}
\put(44.25,140.5){\line(0,-1){32.75}}
\put(45.75,140){\line(0,-1){33}}
\put(26.5,131.25){\line(1,0){19.25}}
\put(26,129.75){\line(1,0){19.75}}
\put(26.5,117.25){\line(1,0){18.75}}
\put(26,115.75){\line(1,0){19.25}}
\put(14.75,134.5){$\bf p$}
\put(14.25,103.25){$\bf q$}
\put(25.75,143.5){$\cal C$}
\put(44.25,144){$\cal C'$}
\put(35,133.5){$\cal S$}
\put(34.5,112.75){$\cal S'$}
\put(33.75,124.25){$E$}
\put(36.75,100.25){$\Gamma$}
\thinlines
%\dashline{1}(29,130.75)(27.5,129.5)
\multiput(28.93,130.68)(-.0384615,-.0320513){13}{\line(-1,0){.0384615}}
\multiput(27.93,129.846)(-.0384615,-.0320513){13}{\line(-1,0){.0384615}}
%\end
%\dashline{1}(25.75,136.25)(31.25,140.25)
\multiput(25.68,136.18)(.0458333,.0333333){15}{\line(1,0){.0458333}}
\multiput(27.055,137.18)(.0458333,.0333333){15}{\line(1,0){.0458333}}
\multiput(28.43,138.18)(.0458333,.0333333){15}{\line(1,0){.0458333}}
\multiput(29.805,139.18)(.0458333,.0333333){15}{\line(1,0){.0458333}}
%\end
%\dashline{1}(26.5,133.75)(34.75,140.25)
\multiput(26.43,133.68)(.0404412,.0318627){17}{\line(1,0){.0404412}}
\multiput(27.805,134.763)(.0404412,.0318627){17}{\line(1,0){.0404412}}
\multiput(29.18,135.846)(.0404412,.0318627){17}{\line(1,0){.0404412}}
\multiput(30.555,136.93)(.0404412,.0318627){17}{\line(1,0){.0404412}}
\multiput(31.93,138.013)(.0404412,.0318627){17}{\line(1,0){.0404412}}
\multiput(33.305,139.096)(.0404412,.0318627){17}{\line(1,0){.0404412}}
%\end
%\dashline{1}(26,131.25)(39,140.25)
\multiput(25.93,131.18)(.0481481,.0333333){15}{\line(1,0){.0481481}}
\multiput(27.374,132.18)(.0481481,.0333333){15}{\line(1,0){.0481481}}
\multiput(28.819,133.18)(.0481481,.0333333){15}{\line(1,0){.0481481}}
\multiput(30.263,134.18)(.0481481,.0333333){15}{\line(1,0){.0481481}}
\multiput(31.707,135.18)(.0481481,.0333333){15}{\line(1,0){.0481481}}
\multiput(33.152,136.18)(.0481481,.0333333){15}{\line(1,0){.0481481}}
\multiput(34.596,137.18)(.0481481,.0333333){15}{\line(1,0){.0481481}}
\multiput(36.041,138.18)(.0481481,.0333333){15}{\line(1,0){.0481481}}
\multiput(37.485,139.18)(.0481481,.0333333){15}{\line(1,0){.0481481}}
%\end
%\dashline{1}(26.25,128)(26.25,128.5)
\put(26.18,127.93){\line(0,1){.25}}
%\end
%\dashline{1}(26,129)(33.5,134.5)
\multiput(25.93,128.93)(.0441176,.0323529){17}{\line(1,0){.0441176}}
\multiput(27.43,130.03)(.0441176,.0323529){17}{\line(1,0){.0441176}}
\multiput(28.93,131.13)(.0441176,.0323529){17}{\line(1,0){.0441176}}
\multiput(30.43,132.23)(.0441176,.0323529){17}{\line(1,0){.0441176}}
\multiput(31.93,133.33)(.0441176,.0323529){17}{\line(1,0){.0441176}}
%\end
%\dashline{1}(33.5,134.5)(41.5,140.25)
\multiput(33.43,134.43)(.0454545,.0326705){16}{\line(1,0){.0454545}}
\multiput(34.884,135.475)(.0454545,.0326705){16}{\line(1,0){.0454545}}
\multiput(36.339,136.521)(.0454545,.0326705){16}{\line(1,0){.0454545}}
\multiput(37.793,137.566)(.0454545,.0326705){16}{\line(1,0){.0454545}}
\multiput(39.248,138.612)(.0454545,.0326705){16}{\line(1,0){.0454545}}
\multiput(40.702,139.657)(.0454545,.0326705){16}{\line(1,0){.0454545}}
%\end
%\dashline{1}(26,126.75)(43.75,140.25)
\multiput(25.93,126.68)(.0428744,.0326087){18}{\line(1,0){.0428744}}
\multiput(27.473,127.854)(.0428744,.0326087){18}{\line(1,0){.0428744}}
\multiput(29.017,129.028)(.0428744,.0326087){18}{\line(1,0){.0428744}}
\multiput(30.56,130.201)(.0428744,.0326087){18}{\line(1,0){.0428744}}
\multiput(32.104,131.375)(.0428744,.0326087){18}{\line(1,0){.0428744}}
\multiput(33.647,132.549)(.0428744,.0326087){18}{\line(1,0){.0428744}}
\multiput(35.191,133.723)(.0428744,.0326087){18}{\line(1,0){.0428744}}
\multiput(36.734,134.897)(.0428744,.0326087){18}{\line(1,0){.0428744}}
\multiput(38.278,136.071)(.0428744,.0326087){18}{\line(1,0){.0428744}}
\multiput(39.821,137.245)(.0428744,.0326087){18}{\line(1,0){.0428744}}
\multiput(41.364,138.419)(.0428744,.0326087){18}{\line(1,0){.0428744}}
\multiput(42.908,139.593)(.0428744,.0326087){18}{\line(1,0){.0428744}}
%\end
%\dashline{1}(25.75,124.75)(45.75,139)
\multiput(25.68,124.68)(.0452489,.0322398){17}{\line(1,0){.0452489}}
\multiput(27.218,125.776)(.0452489,.0322398){17}{\line(1,0){.0452489}}
\multiput(28.757,126.872)(.0452489,.0322398){17}{\line(1,0){.0452489}}
\multiput(30.295,127.968)(.0452489,.0322398){17}{\line(1,0){.0452489}}
\multiput(31.834,129.064)(.0452489,.0322398){17}{\line(1,0){.0452489}}
\multiput(33.372,130.16)(.0452489,.0322398){17}{\line(1,0){.0452489}}
\multiput(34.91,131.257)(.0452489,.0322398){17}{\line(1,0){.0452489}}
\multiput(36.449,132.353)(.0452489,.0322398){17}{\line(1,0){.0452489}}
\multiput(37.987,133.449)(.0452489,.0322398){17}{\line(1,0){.0452489}}
\multiput(39.526,134.545)(.0452489,.0322398){17}{\line(1,0){.0452489}}
\multiput(41.064,135.641)(.0452489,.0322398){17}{\line(1,0){.0452489}}
\multiput(42.603,136.737)(.0452489,.0322398){17}{\line(1,0){.0452489}}
\multiput(44.141,137.834)(.0452489,.0322398){17}{\line(1,0){.0452489}}
%\end
%\dashline{1}(25.25,122.25)(45.75,136.25)
\multiput(25.18,122.18)(.0492788,.0336538){16}{\line(1,0){.0492788}}
\multiput(26.757,123.257)(.0492788,.0336538){16}{\line(1,0){.0492788}}
\multiput(28.334,124.334)(.0492788,.0336538){16}{\line(1,0){.0492788}}
\multiput(29.91,125.41)(.0492788,.0336538){16}{\line(1,0){.0492788}}
\multiput(31.487,126.487)(.0492788,.0336538){16}{\line(1,0){.0492788}}
\multiput(33.064,127.564)(.0492788,.0336538){16}{\line(1,0){.0492788}}
\multiput(34.641,128.641)(.0492788,.0336538){16}{\line(1,0){.0492788}}
\multiput(36.218,129.718)(.0492788,.0336538){16}{\line(1,0){.0492788}}
\multiput(37.795,130.795)(.0492788,.0336538){16}{\line(1,0){.0492788}}
\multiput(39.372,131.872)(.0492788,.0336538){16}{\line(1,0){.0492788}}
\multiput(40.949,132.949)(.0492788,.0336538){16}{\line(1,0){.0492788}}
\multiput(42.526,134.026)(.0492788,.0336538){16}{\line(1,0){.0492788}}
\multiput(44.103,135.103)(.0492788,.0336538){16}{\line(1,0){.0492788}}
%\end
%\dashline{1}(26,120.25)(32.5,124.5)
\multiput(25.93,120.18)(.0515873,.0337302){14}{\line(1,0){.0515873}}
\multiput(27.374,121.124)(.0515873,.0337302){14}{\line(1,0){.0515873}}
\multiput(28.819,122.069)(.0515873,.0337302){14}{\line(1,0){.0515873}}
\multiput(30.263,123.013)(.0515873,.0337302){14}{\line(1,0){.0515873}}
\multiput(31.707,123.957)(.0515873,.0337302){14}{\line(1,0){.0515873}}
%\end
%\dashline{1}(32.5,124.5)(45.5,134)
\multiput(32.43,124.43)(.0451389,.0329861){16}{\line(1,0){.0451389}}
\multiput(33.874,125.485)(.0451389,.0329861){16}{\line(1,0){.0451389}}
\multiput(35.319,126.541)(.0451389,.0329861){16}{\line(1,0){.0451389}}
\multiput(36.763,127.596)(.0451389,.0329861){16}{\line(1,0){.0451389}}
\multiput(38.207,128.652)(.0451389,.0329861){16}{\line(1,0){.0451389}}
\multiput(39.652,129.707)(.0451389,.0329861){16}{\line(1,0){.0451389}}
\multiput(41.096,130.763)(.0451389,.0329861){16}{\line(1,0){.0451389}}
\multiput(42.541,131.819)(.0451389,.0329861){16}{\line(1,0){.0451389}}
\multiput(43.985,132.874)(.0451389,.0329861){16}{\line(1,0){.0451389}}
%\end
%\dashline{1}(25.75,118.5)(24.75,119)
\multiput(25.68,118.43)(-.0625,.03125){8}{\line(-1,0){.0625}}
%\end
%\dashline{1}(26.25,117.75)(45.25,131.25)
\multiput(26.18,117.68)(.0447059,.0317647){17}{\line(1,0){.0447059}}
\multiput(27.7,118.76)(.0447059,.0317647){17}{\line(1,0){.0447059}}
\multiput(29.22,119.84)(.0447059,.0317647){17}{\line(1,0){.0447059}}
\multiput(30.74,120.92)(.0447059,.0317647){17}{\line(1,0){.0447059}}
\multiput(32.26,122)(.0447059,.0317647){17}{\line(1,0){.0447059}}
\multiput(33.78,123.08)(.0447059,.0317647){17}{\line(1,0){.0447059}}
\multiput(35.3,124.16)(.0447059,.0317647){17}{\line(1,0){.0447059}}
\multiput(36.82,125.24)(.0447059,.0317647){17}{\line(1,0){.0447059}}
\multiput(38.34,126.32)(.0447059,.0317647){17}{\line(1,0){.0447059}}
\multiput(39.86,127.4)(.0447059,.0317647){17}{\line(1,0){.0447059}}
\multiput(41.38,128.48)(.0447059,.0317647){17}{\line(1,0){.0447059}}
\multiput(42.9,129.56)(.0447059,.0317647){17}{\line(1,0){.0447059}}
\multiput(44.42,130.64)(.0447059,.0317647){17}{\line(1,0){.0447059}}
%\end
%\dashline{1}(25.75,115.5)(45,129.5)
\multiput(25.68,115.43)(.0452941,.0329412){17}{\line(1,0){.0452941}}
\multiput(27.22,116.55)(.0452941,.0329412){17}{\line(1,0){.0452941}}
\multiput(28.76,117.67)(.0452941,.0329412){17}{\line(1,0){.0452941}}
\multiput(30.3,118.79)(.0452941,.0329412){17}{\line(1,0){.0452941}}
\multiput(31.84,119.91)(.0452941,.0329412){17}{\line(1,0){.0452941}}
\multiput(33.38,121.03)(.0452941,.0329412){17}{\line(1,0){.0452941}}
\multiput(34.92,122.15)(.0452941,.0329412){17}{\line(1,0){.0452941}}
\multiput(36.46,123.27)(.0452941,.0329412){17}{\line(1,0){.0452941}}
\multiput(38,124.39)(.0452941,.0329412){17}{\line(1,0){.0452941}}
\multiput(39.54,125.51)(.0452941,.0329412){17}{\line(1,0){.0452941}}
\multiput(41.08,126.63)(.0452941,.0329412){17}{\line(1,0){.0452941}}
\multiput(42.62,127.75)(.0452941,.0329412){17}{\line(1,0){.0452941}}
\multiput(44.16,128.87)(.0452941,.0329412){17}{\line(1,0){.0452941}}
%\end
%\dashline{1}(26,113.5)(45.5,127.5)
\multiput(25.93,113.43)(.0458824,.0329412){17}{\line(1,0){.0458824}}
\multiput(27.49,114.55)(.0458824,.0329412){17}{\line(1,0){.0458824}}
\multiput(29.05,115.67)(.0458824,.0329412){17}{\line(1,0){.0458824}}
\multiput(30.61,116.79)(.0458824,.0329412){17}{\line(1,0){.0458824}}
\multiput(32.17,117.91)(.0458824,.0329412){17}{\line(1,0){.0458824}}
\multiput(33.73,119.03)(.0458824,.0329412){17}{\line(1,0){.0458824}}
\multiput(35.29,120.15)(.0458824,.0329412){17}{\line(1,0){.0458824}}
\multiput(36.85,121.27)(.0458824,.0329412){17}{\line(1,0){.0458824}}
\multiput(38.41,122.39)(.0458824,.0329412){17}{\line(1,0){.0458824}}
\multiput(39.97,123.51)(.0458824,.0329412){17}{\line(1,0){.0458824}}
\multiput(41.53,124.63)(.0458824,.0329412){17}{\line(1,0){.0458824}}
\multiput(43.09,125.75)(.0458824,.0329412){17}{\line(1,0){.0458824}}
\multiput(44.65,126.87)(.0458824,.0329412){17}{\line(1,0){.0458824}}
%\end
%\dashline{1}(25.75,111.5)(45.75,125.5)
\multiput(25.68,111.43)(.0480769,.0336538){16}{\line(1,0){.0480769}}
\multiput(27.218,112.507)(.0480769,.0336538){16}{\line(1,0){.0480769}}
\multiput(28.757,113.584)(.0480769,.0336538){16}{\line(1,0){.0480769}}
\multiput(30.295,114.66)(.0480769,.0336538){16}{\line(1,0){.0480769}}
\multiput(31.834,115.737)(.0480769,.0336538){16}{\line(1,0){.0480769}}
\multiput(33.372,116.814)(.0480769,.0336538){16}{\line(1,0){.0480769}}
\multiput(34.91,117.891)(.0480769,.0336538){16}{\line(1,0){.0480769}}
\multiput(36.449,118.968)(.0480769,.0336538){16}{\line(1,0){.0480769}}
\multiput(37.987,120.045)(.0480769,.0336538){16}{\line(1,0){.0480769}}
\multiput(39.526,121.122)(.0480769,.0336538){16}{\line(1,0){.0480769}}
\multiput(41.064,122.199)(.0480769,.0336538){16}{\line(1,0){.0480769}}
\multiput(42.603,123.276)(.0480769,.0336538){16}{\line(1,0){.0480769}}
\multiput(44.141,124.353)(.0480769,.0336538){16}{\line(1,0){.0480769}}
%\end
%\dashline{1}(25.75,109.5)(32.25,113.75)
\multiput(25.68,109.43)(.0515873,.0337302){14}{\line(1,0){.0515873}}
\multiput(27.124,110.374)(.0515873,.0337302){14}{\line(1,0){.0515873}}
\multiput(28.569,111.319)(.0515873,.0337302){14}{\line(1,0){.0515873}}
\multiput(30.013,112.263)(.0515873,.0337302){14}{\line(1,0){.0515873}}
\multiput(31.457,113.207)(.0515873,.0337302){14}{\line(1,0){.0515873}}
%\end
%\dashline{1}(32.25,113.75)(45.25,123.25)
\multiput(32.18,113.68)(.0451389,.0329861){16}{\line(1,0){.0451389}}
\multiput(33.624,114.735)(.0451389,.0329861){16}{\line(1,0){.0451389}}
\multiput(35.069,115.791)(.0451389,.0329861){16}{\line(1,0){.0451389}}
\multiput(36.513,116.846)(.0451389,.0329861){16}{\line(1,0){.0451389}}
\multiput(37.957,117.902)(.0451389,.0329861){16}{\line(1,0){.0451389}}
\multiput(39.402,118.957)(.0451389,.0329861){16}{\line(1,0){.0451389}}
\multiput(40.846,120.013)(.0451389,.0329861){16}{\line(1,0){.0451389}}
\multiput(42.291,121.069)(.0451389,.0329861){16}{\line(1,0){.0451389}}
\multiput(43.735,122.124)(.0451389,.0329861){16}{\line(1,0){.0451389}}
%\end
%\dashline{1}(27.25,109)(31,112.25)
\multiput(27.18,108.93)(.0367647,.0318627){17}{\line(1,0){.0367647}}
\multiput(28.43,110.013)(.0367647,.0318627){17}{\line(1,0){.0367647}}
\multiput(29.68,111.096)(.0367647,.0318627){17}{\line(1,0){.0367647}}
%\end
%\dashline{1}(28.5,108.5)(45.25,121)
\multiput(28.43,108.43)(.0447861,.0334225){17}{\line(1,0){.0447861}}
\multiput(29.952,109.566)(.0447861,.0334225){17}{\line(1,0){.0447861}}
\multiput(31.475,110.702)(.0447861,.0334225){17}{\line(1,0){.0447861}}
\multiput(32.998,111.839)(.0447861,.0334225){17}{\line(1,0){.0447861}}
\multiput(34.521,112.975)(.0447861,.0334225){17}{\line(1,0){.0447861}}
\multiput(36.043,114.112)(.0447861,.0334225){17}{\line(1,0){.0447861}}
\multiput(37.566,115.248)(.0447861,.0334225){17}{\line(1,0){.0447861}}
\multiput(39.089,116.384)(.0447861,.0334225){17}{\line(1,0){.0447861}}
\multiput(40.612,117.521)(.0447861,.0334225){17}{\line(1,0){.0447861}}
\multiput(42.134,118.657)(.0447861,.0334225){17}{\line(1,0){.0447861}}
\multiput(43.657,119.793)(.0447861,.0334225){17}{\line(1,0){.0447861}}
%\end
%\dashline{1}(31,107.5)(30,108)
\multiput(30.93,107.43)(-.0625,.03125){8}{\line(-1,0){.0625}}
%\end
%\dashline{1}(30,108)(45.5,118.5)
\multiput(29.93,107.93)(.0484375,.0328125){16}{\line(1,0){.0484375}}
\multiput(31.48,108.98)(.0484375,.0328125){16}{\line(1,0){.0484375}}
\multiput(33.03,110.03)(.0484375,.0328125){16}{\line(1,0){.0484375}}
\multiput(34.58,111.08)(.0484375,.0328125){16}{\line(1,0){.0484375}}
\multiput(36.13,112.13)(.0484375,.0328125){16}{\line(1,0){.0484375}}
\multiput(37.68,113.18)(.0484375,.0328125){16}{\line(1,0){.0484375}}
\multiput(39.23,114.23)(.0484375,.0328125){16}{\line(1,0){.0484375}}
\multiput(40.78,115.28)(.0484375,.0328125){16}{\line(1,0){.0484375}}
\multiput(42.33,116.33)(.0484375,.0328125){16}{\line(1,0){.0484375}}
\multiput(43.88,117.38)(.0484375,.0328125){16}{\line(1,0){.0484375}}
%\end
%\dashline{1}(33.5,108)(45.25,116.75)
\multiput(33.43,107.93)(.0435185,.0324074){18}{\line(1,0){.0435185}}
\multiput(34.996,109.096)(.0435185,.0324074){18}{\line(1,0){.0435185}}
\multiput(36.563,110.263)(.0435185,.0324074){18}{\line(1,0){.0435185}}
\multiput(38.13,111.43)(.0435185,.0324074){18}{\line(1,0){.0435185}}
\multiput(39.696,112.596)(.0435185,.0324074){18}{\line(1,0){.0435185}}
\multiput(41.263,113.763)(.0435185,.0324074){18}{\line(1,0){.0435185}}
\multiput(42.83,114.93)(.0435185,.0324074){18}{\line(1,0){.0435185}}
\multiput(44.396,116.096)(.0435185,.0324074){18}{\line(1,0){.0435185}}
%\end
%\dashline{1}(35.75,108)(45.75,114.75)
\multiput(35.68,107.93)(.0480769,.0324519){16}{\line(1,0){.0480769}}
\multiput(37.218,108.968)(.0480769,.0324519){16}{\line(1,0){.0480769}}
\multiput(38.757,110.007)(.0480769,.0324519){16}{\line(1,0){.0480769}}
\multiput(40.295,111.045)(.0480769,.0324519){16}{\line(1,0){.0480769}}
\multiput(41.834,112.084)(.0480769,.0324519){16}{\line(1,0){.0480769}}
\multiput(43.372,113.122)(.0480769,.0324519){16}{\line(1,0){.0480769}}
\multiput(44.91,114.16)(.0480769,.0324519){16}{\line(1,0){.0480769}}
%\end
%\dashline{1}(38.75,107.5)(38,107.75)
\multiput(38.68,107.43)(-.09375,.03125){4}{\line(-1,0){.09375}}
%\end
%\dashline{1}(38,107.75)(45.25,112.75)
\multiput(37.93,107.68)(.0483333,.0333333){15}{\line(1,0){.0483333}}
\multiput(39.38,108.68)(.0483333,.0333333){15}{\line(1,0){.0483333}}
\multiput(40.83,109.68)(.0483333,.0333333){15}{\line(1,0){.0483333}}
\multiput(42.28,110.68)(.0483333,.0333333){15}{\line(1,0){.0483333}}
\multiput(43.73,111.68)(.0483333,.0333333){15}{\line(1,0){.0483333}}
%\end
%\dashline{1}(41.5,108)(45.25,111)
\multiput(41.43,107.93)(.0416667,.0333333){15}{\line(1,0){.0416667}}
\multiput(42.68,108.93)(.0416667,.0333333){15}{\line(1,0){.0416667}}
\multiput(43.93,109.93)(.0416667,.0333333){15}{\line(1,0){.0416667}}
%\end
\thicklines
%\circle(103,121.75){39.756}
\put(122.878,121.75){\line(0,1){.9197}}
\put(122.856,122.67){\line(0,1){.9177}}
\put(122.793,123.587){\line(0,1){.9138}}
\multiput(122.686,124.501)(-.029674,.181576){5}{\line(0,1){.181576}}
\multiput(122.538,125.409)(-.031702,.150007){6}{\line(0,1){.150007}}
\multiput(122.348,126.309)(-.033093,.127182){7}{\line(0,1){.127182}}
\multiput(122.116,127.199)(-.030288,.097623){9}{\line(0,1){.097623}}
\multiput(121.844,128.078)(-.031295,.086505){10}{\line(0,1){.086505}}
\multiput(121.531,128.943)(-.032058,.077241){11}{\line(0,1){.077241}}
\multiput(121.178,129.793)(-.032631,.069369){12}{\line(0,1){.069369}}
\multiput(120.786,130.625)(-.0330513,.0625703){13}{\line(0,1){.0625703}}
\multiput(120.357,131.438)(-.0333457,.0566189){14}{\line(0,1){.0566189}}
\multiput(119.89,132.231)(-.0335343,.0513477){15}{\line(0,1){.0513477}}
\multiput(119.387,133.001)(-.0336319,.0466324){16}{\line(0,1){.0466324}}
\multiput(118.849,133.747)(-.0336503,.0423778){17}{\line(0,1){.0423778}}
\multiput(118.277,134.468)(-.0335985,.0385103){18}{\line(0,1){.0385103}}
\multiput(117.672,135.161)(-.0334841,.0349717){19}{\line(0,1){.0349717}}
\multiput(117.036,135.826)(-.0350662,.033385){19}{\line(-1,0){.0350662}}
\multiput(116.37,136.46)(-.0386051,.0334895){18}{\line(-1,0){.0386051}}
\multiput(115.675,137.063)(-.0424729,.0335303){17}{\line(-1,0){.0424729}}
\multiput(114.953,137.633)(-.0467274,.0334999){16}{\line(-1,0){.0467274}}
\multiput(114.205,138.169)(-.0514424,.0333889){15}{\line(-1,0){.0514424}}
\multiput(113.433,138.67)(-.056713,.0331854){14}{\line(-1,0){.056713}}
\multiput(112.639,139.134)(-.0626636,.0328741){13}{\line(-1,0){.0626636}}
\multiput(111.825,139.561)(-.069461,.032435){12}{\line(-1,0){.069461}}
\multiput(110.991,139.951)(-.077331,.03184){11}{\line(-1,0){.077331}}
\multiput(110.141,140.301)(-.086594,.03105){10}{\line(-1,0){.086594}}
\multiput(109.275,140.611)(-.097708,.030012){9}{\line(-1,0){.097708}}
\multiput(108.395,140.882)(-.127276,.032733){7}{\line(-1,0){.127276}}
\multiput(107.504,141.111)(-.150096,.031278){6}{\line(-1,0){.150096}}
\multiput(106.604,141.298)(-.181659,.02916){5}{\line(-1,0){.181659}}
\put(105.695,141.444){\line(-1,0){.9141}}
\put(104.781,141.548){\line(-1,0){.9179}}
\put(103.863,141.609){\line(-1,0){1.8393}}
\put(102.024,141.604){\line(-1,0){.9175}}
\put(101.107,141.537){\line(-1,0){.9135}}
\multiput(100.193,141.429)(-.181491,-.030187){5}{\line(-1,0){.181491}}
\multiput(99.286,141.278)(-.149917,-.032127){6}{\line(-1,0){.149917}}
\multiput(98.386,141.085)(-.127088,-.033453){7}{\line(-1,0){.127088}}
\multiput(97.497,140.851)(-.097537,-.030564){9}{\line(-1,0){.097537}}
\multiput(96.619,140.576)(-.086416,-.03154){10}{\line(-1,0){.086416}}
\multiput(95.755,140.26)(-.07715,-.032277){11}{\line(-1,0){.07715}}
\multiput(94.906,139.905)(-.069276,-.032827){12}{\line(-1,0){.069276}}
\multiput(94.075,139.511)(-.0624766,-.0332282){13}{\line(-1,0){.0624766}}
\multiput(93.262,139.079)(-.0565243,-.0335058){14}{\line(-1,0){.0565243}}
\multiput(92.471,138.61)(-.0512527,-.0336794){15}{\line(-1,0){.0512527}}
\multiput(91.702,138.105)(-.0437996,-.0317776){17}{\line(-1,0){.0437996}}
\multiput(90.958,137.565)(-.0399334,-.0318939){18}{\line(-1,0){.0399334}}
\multiput(90.239,136.991)(-.038415,-.0337074){18}{\line(-1,0){.038415}}
\multiput(89.547,136.384)(-.0348768,-.0335829){19}{\line(-1,0){.0348768}}
\multiput(88.885,135.746)(-.0332857,-.0351605){19}{\line(0,-1){.0351605}}
\multiput(88.252,135.078)(-.0333801,-.0386997){18}{\line(0,-1){.0386997}}
\multiput(87.652,134.381)(-.03341,-.0425675){17}{\line(0,-1){.0425675}}
\multiput(87.084,133.658)(-.0333675,-.046822){16}{\line(0,-1){.046822}}
\multiput(86.55,132.908)(-.0332432,-.0515367){15}{\line(0,-1){.0515367}}
\multiput(86.051,132.135)(-.0330249,-.0568066){14}{\line(0,-1){.0568066}}
\multiput(85.589,131.34)(-.0326967,-.0627563){13}{\line(0,-1){.0627563}}
\multiput(85.164,130.524)(-.032238,-.069552){12}{\line(0,-1){.069552}}
\multiput(84.777,129.69)(-.031621,-.077421){11}{\line(0,-1){.077421}}
\multiput(84.429,128.838)(-.030805,-.086681){10}{\line(0,-1){.086681}}
\multiput(84.121,127.971)(-.033452,-.110017){8}{\line(0,-1){.110017}}
\multiput(83.853,127.091)(-.032373,-.127368){7}{\line(0,-1){.127368}}
\multiput(83.627,126.2)(-.030853,-.150184){6}{\line(0,-1){.150184}}
\multiput(83.442,125.298)(-.028646,-.181741){5}{\line(0,-1){.181741}}
\put(83.298,124.39){\line(0,-1){.9144}}
\put(83.197,123.475){\line(0,-1){.9181}}
\put(83.139,122.557){\line(0,-1){.9198}}
\put(83.123,121.638){\line(0,-1){.9195}}
\put(83.149,120.718){\line(0,-1){.9173}}
\put(83.218,119.801){\line(0,-1){.9132}}
\multiput(83.329,118.888)(.0307,-.181405){5}{\line(0,-1){.181405}}
\multiput(83.483,117.98)(.03255,-.149825){6}{\line(0,-1){.149825}}
\multiput(83.678,117.082)(.029586,-.111119){8}{\line(0,-1){.111119}}
\multiput(83.915,116.193)(.03084,-.09745){9}{\line(0,-1){.09745}}
\multiput(84.192,115.316)(.031784,-.086327){10}{\line(0,-1){.086327}}
\multiput(84.51,114.452)(.032495,-.077058){11}{\line(0,-1){.077058}}
\multiput(84.868,113.605)(.033023,-.069183){12}{\line(0,-1){.069183}}
\multiput(85.264,112.774)(.0334048,-.0623823){13}{\line(0,-1){.0623823}}
\multiput(85.698,111.963)(.0336656,-.0564293){14}{\line(0,-1){.0564293}}
\multiput(86.17,111.173)(.0317103,-.0479599){16}{\line(0,-1){.0479599}}
\multiput(86.677,110.406)(.0319014,-.0437095){17}{\line(0,-1){.0437095}}
\multiput(87.219,109.663)(.0320068,-.039843){18}{\line(0,-1){.039843}}
\multiput(87.795,108.946)(.0320361,-.0363027){19}{\line(0,-1){.0363027}}
\multiput(88.404,108.256)(.0336814,-.0347816){19}{\line(0,-1){.0347816}}
\multiput(89.044,107.595)(.0352546,-.0331861){19}{\line(1,0){.0352546}}
\multiput(89.714,106.965)(.038794,-.0332705){18}{\line(1,0){.038794}}
\multiput(90.412,106.366)(.0426619,-.0332894){17}{\line(1,0){.0426619}}
\multiput(91.137,105.8)(.0469162,-.0332349){16}{\line(1,0){.0469162}}
\multiput(91.888,105.268)(.0516305,-.0330973){15}{\line(1,0){.0516305}}
\multiput(92.663,104.772)(.0568998,-.032864){14}{\line(1,0){.0568998}}
\multiput(93.459,104.312)(.0628486,-.032519){13}{\line(1,0){.0628486}}
\multiput(94.276,103.889)(.069643,-.032041){12}{\line(1,0){.069643}}
\multiput(95.112,103.504)(.07751,-.031401){11}{\line(1,0){.07751}}
\multiput(95.965,103.159)(.086768,-.03056){10}{\line(1,0){.086768}}
\multiput(96.832,102.853)(.110111,-.033141){8}{\line(1,0){.110111}}
\multiput(97.713,102.588)(.127459,-.032013){7}{\line(1,0){.127459}}
\multiput(98.605,102.364)(.150271,-.030428){6}{\line(1,0){.150271}}
\multiput(99.507,102.182)(.181821,-.028132){5}{\line(1,0){.181821}}
\put(100.416,102.041){\line(1,0){.9146}}
\put(101.331,101.942){\line(1,0){.9182}}
\put(102.249,101.886){\line(1,0){.9198}}
\put(103.169,101.873){\line(1,0){.9195}}
\put(104.088,101.902){\line(1,0){.9171}}
\put(105.005,101.974){\line(1,0){.9128}}
\multiput(105.918,102.088)(.181317,.031214){5}{\line(1,0){.181317}}
\multiput(106.825,102.244)(.149732,.032974){6}{\line(1,0){.149732}}
\multiput(107.723,102.442)(.111035,.0299){8}{\line(1,0){.111035}}
\multiput(108.611,102.681)(.097362,.031116){9}{\line(1,0){.097362}}
\multiput(109.488,102.961)(.086237,.032028){10}{\line(1,0){.086237}}
\multiput(110.35,103.281)(.076966,.032713){11}{\line(1,0){.076966}}
\multiput(111.197,103.641)(.069089,.033219){12}{\line(1,0){.069089}}
\multiput(112.026,104.04)(.0622876,.0335811){13}{\line(1,0){.0622876}}
\multiput(112.835,104.476)(.0525782,.0315701){15}{\line(1,0){.0525782}}
\multiput(113.624,104.95)(.04787,.0318458){16}{\line(1,0){.04787}}
\multiput(114.39,105.459)(.0436191,.0320249){17}{\line(1,0){.0436191}}
\multiput(115.132,106.004)(.0397523,.0321194){18}{\line(1,0){.0397523}}
\multiput(115.847,106.582)(.0362119,.0321387){19}{\line(1,0){.0362119}}
\multiput(116.535,107.192)(.0329519,.0320907){20}{\line(1,0){.0329519}}
\multiput(117.194,107.834)(.0330862,.0353483){19}{\line(0,1){.0353483}}
\multiput(117.823,108.506)(.0331606,.038888){18}{\line(0,1){.038888}}
\multiput(118.42,109.206)(.0331686,.0427559){17}{\line(0,1){.0427559}}
\multiput(118.984,109.933)(.0331021,.04701){16}{\line(0,1){.04701}}
\multiput(119.513,110.685)(.0329511,.0517239){15}{\line(0,1){.0517239}}
\multiput(120.007,111.461)(.0327029,.0569926){14}{\line(0,1){.0569926}}
\multiput(120.465,112.259)(.0323411,.0629403){13}{\line(0,1){.0629403}}
\multiput(120.886,113.077)(.031844,.069733){12}{\line(0,1){.069733}}
\multiput(121.268,113.914)(.031182,.077599){11}{\line(0,1){.077599}}
\multiput(121.611,114.767)(.033683,.096504){9}{\line(0,1){.096504}}
\multiput(121.914,115.636)(.032829,.110204){8}{\line(0,1){.110204}}
\multiput(122.177,116.517)(.031652,.127549){7}{\line(0,1){.127549}}
\multiput(122.398,117.41)(.030003,.150356){6}{\line(0,1){.150356}}
\multiput(122.578,118.312)(.027617,.1819){5}{\line(0,1){.1819}}
\put(122.716,119.222){\line(0,1){.9149}}
\put(122.812,120.137){\line(0,1){1.6133}}
%\end
\put(101.5,142){\line(0,-1){15}}
\put(103.5,141.25){\line(0,-1){14.25}}
%\circle*(102.5,121.75){10.012}
\put(98.9038,118.1538){\rule{7.1924\unitlength}{7.1924\unitlength}}
\multiput(100.528,125.2337)(0,-8.1651){2}{\rule{3.9441\unitlength}{1.1977\unitlength}}
\multiput(101.4671,126.3189)(0,-9.5352){2}{\rule{2.0658\unitlength}{.3974\unitlength}}
\multiput(101.9531,126.6038)(0,-9.8922){2}{\rule{1.0939\unitlength}{.1846\unitlength}}
\multiput(103.4204,126.3189)(-2.4299,0){2}{\multiput(0,0)(0,-9.4158){2}{\rule{.589\unitlength}{.278\unitlength}}}
\multiput(103.4204,126.4844)(-2.1931,0){2}{\multiput(0,0)(0,-9.6469){2}{\rule{.3522\unitlength}{.178\unitlength}}}
\multiput(103.897,126.3189)(-3.1398,0){2}{\multiput(0,0)(0,-9.3388){2}{\rule{.3458\unitlength}{.2009\unitlength}}}
\multiput(104.3596,125.2337)(-4.6971,0){2}{\multiput(0,0)(0,-7.7025){2}{\rule{.978\unitlength}{.7351\unitlength}}}
\multiput(104.3596,125.8563)(-4.2757,0){2}{\multiput(0,0)(0,-8.5777){2}{\rule{.5566\unitlength}{.365\unitlength}}}
\multiput(104.3596,126.1089)(-4.0563,0){2}{\multiput(0,0)(0,-8.9407){2}{\rule{.3371\unitlength}{.2229\unitlength}}}
\multiput(104.5842,126.1089)(-4.3913,0){2}{\multiput(0,0)(0,-8.8868){2}{\rule{.2229\unitlength}{.1691\unitlength}}}
\multiput(104.8037,125.8563)(-4.9337,0){2}{\multiput(0,0)(0,-8.4565){2}{\rule{.3263\unitlength}{.2439\unitlength}}}
\multiput(104.8037,125.9878)(-4.8275,0){2}{\multiput(0,0)(0,-8.6499){2}{\rule{.2201\unitlength}{.1743\unitlength}}}
\multiput(105.0175,125.8563)(-5.2521,0){2}{\multiput(0,0)(0,-8.3921){2}{\rule{.2171\unitlength}{.1795\unitlength}}}
\multiput(105.2251,125.2337)(-5.9573,0){2}{\multiput(0,0)(0,-7.4098){2}{\rule{.5071\unitlength}{.4424\unitlength}}}
\multiput(105.2251,125.5636)(-5.7635,0){2}{\multiput(0,0)(0,-7.8909){2}{\rule{.3134\unitlength}{.2636\unitlength}}}
\multiput(105.2251,125.7148)(-5.664,0){2}{\multiput(0,0)(0,-8.1141){2}{\rule{.2138\unitlength}{.1844\unitlength}}}
\multiput(105.426,125.5636)(-6.0622,0){2}{\multiput(0,0)(0,-7.8165){2}{\rule{.2103\unitlength}{.1892\unitlength}}}
\multiput(105.6197,125.2337)(-6.5379,0){2}{\multiput(0,0)(0,-7.2493){2}{\rule{.2985\unitlength}{.2819\unitlength}}}
\multiput(105.6197,125.4032)(-6.4459,0){2}{\multiput(0,0)(0,-7.5002){2}{\rule{.2065\unitlength}{.1939\unitlength}}}
\multiput(105.8058,125.2337)(-6.814,0){2}{\multiput(0,0)(0,-7.1657){2}{\rule{.2025\unitlength}{.1983\unitlength}}}
\multiput(105.9837,119.778)(-8.1651,0){2}{\rule{1.1977\unitlength}{3.9441\unitlength}}
\multiput(105.9837,123.6096)(-7.7025,0){2}{\multiput(0,0)(0,-4.6971){2}{\rule{.7351\unitlength}{.978\unitlength}}}
\multiput(105.9837,124.4751)(-7.4098,0){2}{\multiput(0,0)(0,-5.9573){2}{\rule{.4424\unitlength}{.5071\unitlength}}}
\multiput(105.9837,124.8697)(-7.2493,0){2}{\multiput(0,0)(0,-6.5379){2}{\rule{.2819\unitlength}{.2985\unitlength}}}
\multiput(105.9837,125.0558)(-7.1657,0){2}{\multiput(0,0)(0,-6.814){2}{\rule{.1983\unitlength}{.2025\unitlength}}}
\multiput(106.1532,124.8697)(-7.5002,0){2}{\multiput(0,0)(0,-6.4459){2}{\rule{.1939\unitlength}{.2065\unitlength}}}
\multiput(106.3136,124.4751)(-7.8909,0){2}{\multiput(0,0)(0,-5.7635){2}{\rule{.2636\unitlength}{.3134\unitlength}}}
\multiput(106.3136,124.676)(-7.8165,0){2}{\multiput(0,0)(0,-6.0622){2}{\rule{.1892\unitlength}{.2103\unitlength}}}
\multiput(106.4648,124.4751)(-8.1141,0){2}{\multiput(0,0)(0,-5.664){2}{\rule{.1844\unitlength}{.2138\unitlength}}}
\multiput(106.6063,123.6096)(-8.5777,0){2}{\multiput(0,0)(0,-4.2757){2}{\rule{.365\unitlength}{.5566\unitlength}}}
\multiput(106.6063,124.0537)(-8.4565,0){2}{\multiput(0,0)(0,-4.9337){2}{\rule{.2439\unitlength}{.3263\unitlength}}}
\multiput(106.6063,124.2675)(-8.3921,0){2}{\multiput(0,0)(0,-5.2521){2}{\rule{.1795\unitlength}{.2171\unitlength}}}
\multiput(106.7378,124.0537)(-8.6499,0){2}{\multiput(0,0)(0,-4.8275){2}{\rule{.1743\unitlength}{.2201\unitlength}}}
\multiput(106.8589,123.6096)(-8.9407,0){2}{\multiput(0,0)(0,-4.0563){2}{\rule{.2229\unitlength}{.3371\unitlength}}}
\multiput(106.8589,123.8342)(-8.8868,0){2}{\multiput(0,0)(0,-4.3913){2}{\rule{.1691\unitlength}{.2229\unitlength}}}
\multiput(107.0689,120.7171)(-9.5352,0){2}{\rule{.3974\unitlength}{2.0658\unitlength}}
\multiput(107.0689,122.6704)(-9.4158,0){2}{\multiput(0,0)(0,-2.4299){2}{\rule{.278\unitlength}{.589\unitlength}}}
\multiput(107.0689,123.147)(-9.3388,0){2}{\multiput(0,0)(0,-3.1398){2}{\rule{.2009\unitlength}{.3458\unitlength}}}
\multiput(107.2344,122.6704)(-9.6469,0){2}{\multiput(0,0)(0,-2.1931){2}{\rule{.178\unitlength}{.3522\unitlength}}}
\multiput(107.3538,121.2031)(-9.8922,0){2}{\rule{.1846\unitlength}{1.0939\unitlength}}
\put(107.506,121.75){\line(0,1){.592}}
\put(107.471,122.342){\line(0,1){.2934}}
\put(107.427,122.635){\line(0,1){.2903}}
\put(107.366,122.926){\line(0,1){.2862}}
\put(107.288,123.212){\line(0,1){.281}}
\put(107.193,123.493){\line(0,1){.2749}}
\multiput(107.082,123.768)(-.03189,.06696){4}{\line(0,1){.06696}}
\multiput(106.954,124.036)(-.028643,.051962){5}{\line(0,1){.051962}}
\multiput(106.811,124.296)(-.03167,.050174){5}{\line(0,1){.050174}}
\multiput(106.652,124.546)(-.028822,.040175){6}{\line(0,1){.040175}}
\multiput(106.479,124.787)(-.031151,.038398){6}{\line(0,1){.038398}}
\multiput(106.293,125.018)(-.033371,.036485){6}{\line(0,1){.036485}}
\multiput(106.092,125.237)(-.030406,.029524){7}{\line(-1,0){.030406}}
\multiput(105.879,125.443)(-.037452,.032282){6}{\line(-1,0){.037452}}
\multiput(105.655,125.637)(-.039298,.030007){6}{\line(-1,0){.039298}}
\multiput(105.419,125.817)(-.049208,.033152){5}{\line(-1,0){.049208}}
\multiput(105.173,125.983)(-.051085,.030179){5}{\line(-1,0){.051085}}
\multiput(104.918,126.134)(-.052783,.0271){5}{\line(-1,0){.052783}}
\put(104.654,126.269){\line(-1,0){.2715}}
\put(104.382,126.389){\line(-1,0){.2781}}
\put(104.104,126.492){\line(-1,0){.2837}}
\put(103.82,126.579){\line(-1,0){.2884}}
\put(103.532,126.649){\line(-1,0){.292}}
\put(103.24,126.701){\line(-1,0){.2946}}
\put(102.945,126.736){\line(-1,0){.2961}}
\put(102.649,126.754){\line(-1,0){.2967}}
\put(102.353,126.754){\line(-1,0){.2962}}
\put(102.056,126.737){\line(-1,0){.2946}}
\put(101.762,126.702){\line(-1,0){.292}}
\put(101.47,126.649){\line(-1,0){.2884}}
\put(101.181,126.579){\line(-1,0){.2838}}
\put(100.898,126.493){\line(-1,0){.2781}}
\put(100.62,126.39){\line(-1,0){.2715}}
\multiput(100.348,126.27)(-.052793,-.027081){5}{\line(-1,0){.052793}}
\multiput(100.084,126.135)(-.051096,-.03016){5}{\line(-1,0){.051096}}
\multiput(99.829,125.984)(-.04922,-.033134){5}{\line(-1,0){.04922}}
\multiput(99.583,125.818)(-.039309,-.029993){6}{\line(-1,0){.039309}}
\multiput(99.347,125.638)(-.037464,-.032269){6}{\line(-1,0){.037464}}
\multiput(99.122,125.445)(-.030417,-.029512){7}{\line(-1,0){.030417}}
\multiput(98.909,125.238)(-.033385,-.036473){6}{\line(0,-1){.036473}}
\multiput(98.709,125.019)(-.031166,-.038386){6}{\line(0,-1){.038386}}
\multiput(98.522,124.789)(-.028837,-.040165){6}{\line(0,-1){.040165}}
\multiput(98.349,124.548)(-.031689,-.050163){5}{\line(0,-1){.050163}}
\multiput(98.19,124.297)(-.028662,-.051952){5}{\line(0,-1){.051952}}
\multiput(98.047,124.037)(-.03192,-.06695){4}{\line(0,-1){.06695}}
\put(97.919,123.77){\line(0,-1){.2749}}
\put(97.808,123.495){\line(0,-1){.281}}
\put(97.713,123.214){\line(0,-1){.2862}}
\put(97.634,122.928){\line(0,-1){.2903}}
\put(97.573,122.637){\line(0,-1){.2934}}
\put(97.529,122.344){\line(0,-1){1.4775}}
\put(97.572,120.866){\line(0,-1){.2903}}
\put(97.633,120.576){\line(0,-1){.2862}}
\put(97.711,120.29){\line(0,-1){.2811}}
\put(97.806,120.009){\line(0,-1){.275}}
\multiput(97.918,119.734)(.03187,-.06697){4}{\line(0,-1){.06697}}
\multiput(98.045,119.466)(.028624,-.051973){5}{\line(0,-1){.051973}}
\multiput(98.188,119.206)(.031652,-.050186){5}{\line(0,-1){.050186}}
\multiput(98.347,118.955)(.028807,-.040186){6}{\line(0,-1){.040186}}
\multiput(98.519,118.714)(.031137,-.038409){6}{\line(0,-1){.038409}}
\multiput(98.706,118.484)(.033358,-.036497){6}{\line(0,-1){.036497}}
\multiput(98.906,118.265)(.030395,-.029535){7}{\line(1,0){.030395}}
\multiput(99.119,118.058)(.03744,-.032296){6}{\line(1,0){.03744}}
\multiput(99.344,117.864)(.039287,-.030022){6}{\line(1,0){.039287}}
\multiput(99.58,117.684)(.049196,-.03317){5}{\line(1,0){.049196}}
\multiput(99.825,117.518)(.051074,-.030198){5}{\line(1,0){.051074}}
\multiput(100.081,117.367)(.052773,-.02712){5}{\line(1,0){.052773}}
\put(100.345,117.231){\line(1,0){.2714}}
\put(100.616,117.112){\line(1,0){.278}}
\put(100.894,117.008){\line(1,0){.2837}}
\put(101.178,116.921){\line(1,0){.2883}}
\put(101.466,116.852){\line(1,0){.292}}
\put(101.758,116.799){\line(1,0){.2946}}
\put(102.053,116.764){\line(1,0){.2961}}
\put(102.349,116.746){\line(1,0){.2967}}
\put(102.646,116.746){\line(1,0){.2962}}
\put(102.942,116.763){\line(1,0){.2946}}
\put(103.236,116.798){\line(1,0){.292}}
\put(103.528,116.851){\line(1,0){.2884}}
\put(103.817,116.92){\line(1,0){.2838}}
\put(104.101,117.007){\line(1,0){.2782}}
\put(104.379,117.11){\line(1,0){.2716}}
\multiput(104.65,117.229)(.052803,.027061){5}{\line(1,0){.052803}}
\multiput(104.914,117.364)(.051108,.030141){5}{\line(1,0){.051108}}
\multiput(105.17,117.515)(.049232,.033116){5}{\line(1,0){.049232}}
\multiput(105.416,117.681)(.03932,.029978){6}{\line(1,0){.03932}}
\multiput(105.652,117.861)(.037476,.032255){6}{\line(1,0){.037476}}
\multiput(105.877,118.054)(.030428,.029501){7}{\line(1,0){.030428}}
\multiput(106.09,118.261)(.033398,.03646){6}{\line(0,1){.03646}}
\multiput(106.29,118.479)(.03118,.038374){6}{\line(0,1){.038374}}
\multiput(106.477,118.71)(.028852,.040154){6}{\line(0,1){.040154}}
\multiput(106.65,118.95)(.031708,.050151){5}{\line(0,1){.050151}}
\multiput(106.809,119.201)(.028681,.051941){5}{\line(0,1){.051941}}
\multiput(106.952,119.461)(.03194,.06694){4}{\line(0,1){.06694}}
\put(107.08,119.729){\line(0,1){.2748}}
\put(107.192,120.004){\line(0,1){.281}}
\put(107.287,120.285){\line(0,1){.2861}}
\put(107.365,120.571){\line(0,1){.2903}}
\put(107.427,120.861){\line(0,1){.2934}}
\put(107.471,121.154){\line(0,1){.5957}}
%\end
\put(103.875,136.25){\oval(.25,0)[r]}
\put(103.75,135.375){\oval(.5,.25)[]}
%\emline(103.25,136.5)(111.5,134)
\multiput(103.25,136.5)(.11,-.03333333){75}{\line(1,0){.11}}
%\end
%\emline(111.5,134)(115.75,126)
\multiput(111.5,134)(.033730159,-.063492063){126}{\line(0,-1){.063492063}}
%\end
%\emline(115.75,126)(112.75,116)
\multiput(115.75,126)(-.03370787,-.11235955){89}{\line(0,-1){.11235955}}
%\end
%\emline(112.75,116)(103,109.25)
\multiput(112.75,116)(-.048507463,-.03358209){201}{\line(-1,0){.048507463}}
%\end
%\emline(103,109.25)(92.25,114.75)
\multiput(103,109.25)(-.06554878,.033536585){164}{\line(-1,0){.06554878}}
%\end
\put(92.25,114.75){\line(0,1){11.25}}
%\emline(92.25,126)(103.5,132.5)
\multiput(92.25,126)(.058290155,.033678756){193}{\line(1,0){.058290155}}
%\end
\put(103.5,132.5){\line(0,1){0}}
\put(102.25,145.75){${\cal C}_i$}
\put(22,123){$\cal D$}
\put(47.5,123.5){$\cal D'$}
\put(110.5,111.25){$\bf z$}
\put(90.5,132.5){$\Delta$}
\put(20.5,91.75){The subdiagram $\Gamma$\;\;\;\;\;\;\;\;\;\;\;\;\;\;\;\;\;\;\;\;\;\;\;\;\;\;\;\;\;\;\;\;\;\;\;\;\;\;\;\;The loop $\bf z$ }
\end{picture}

  If ${\cal T}_0$ is a part of some $\theta$-band $\cal T$, such that ${\cal T}_0$  starts and ends on the same $q$-band ${\cal C}_i$ and crosses once every other maximal $q$-band, then the above argument
  provides us with a cubic upper bound for the length of ${\cal T}_0$.   The ends of the side of ${\cal T}_0$ are connected by a part
  of length $O(r)$ of the $q$-band ${\cal C}_i$. So this side and the
  connecting path form a loop $\bf z$ of at most cubic length surrounding the hole of $\Delta$.

  If $\bf z$ surrounds another loop $\bf z'$ of this type with $Lab ({\bf z'})\equiv Lab ({\bf z})$, then one can remove the diagram between $\bf z$ and $\bf z'$ and identify these
  two loops decreasing the number of cells in the annular diagram.
  Since the lengths of loops of this type are bounded, the number of
  them is effectively bounded as well. It follows that we have an
  effective upper bound for the lengths of $q$-bands in $\Delta$, and so, for the length of a path connecting $\bf p$ and $\bf q$.
  The lemma is proved. \endproof

  Lemmas \ref{coninf} and \ref{finor} prove {\bf Theorem \ref{c}}.

  \section{The power conjugacy problem}\label{pow}

  \begin{lemma}\label{enough} (a) There is an algorithm such that given a word $W$ in the generators of the group $G$, it decides whether the order of $W$ in $G$ is finite or infinite.

 (b) To obtain an algorithm solving the power conjugacy problem in $G$ for arbitrary pairs of words $(U,V)$, it suffices to obtain such an algorithm under assumption that the words $U$ and $V$ have infinite order in $G$.
 \end{lemma}

 \proof (a) By Lemmas \ref{order} and \ref{finor}, $W$ has a finite
 order in $G$ if and only if $W^L=1$ in $G$. Since Lemma \ref{quadr}
 implies that the word problem is decidable in $G$, Statement (a) is proved.

 (b) This statement follows from (a) since some positive powers
 of elements having finite orders are trivial, and so are conjugate.
\endproof

  The relations (\ref{rel1} - \ref{rel3}) defining the group $G$ immediately imply that
  there exists a homomorphism $\mu$ from $G$ to the additive group $\Z/L\Z$ which sends all $\tt$-generators to 1 and all other generators to $0$. Since to solve the power conjugacy problem in $G$ for given
pair of words $(U,V)$, one may replace this pair with the pair $(U^L, V^L)$, we may {\bf assume further} that

\begin{equation} \label{0i0}
\mu(U) = 0\;\;\; and \;\;\; \mu(V)=0
\end{equation}

Let $U$ and $V$ be two words representing elements of infinite order from the group $G$.
Under the condition (\ref{0i0}) and the assumption that the words $U$ and $V$ are adapted, we will also assume that some powers $U^k$ and $V^l$
are conjugate in $G$ for $k,l \ne 0$. Without loss of generality we
assume that $k,l >0$, and there is no pair of positive exponents $k',l'$
such that $U^{k'}$ is a conjugate of $V^{l'}$, where $k'<k$, $l'<l$.
Thus, we will study a minimal annular diagram $\Delta$, where the outer
boundary component $\bf p$ has the clockwise label $U^k$ and the  inner
boundary component $\bf q$ is labeled by $V^l$.

If two consecutive $\tt$-spokes $\cal C$  and $\cal C'$ start on a disk $\Pi$ of $\Delta$, end both on $\bf p$ (or both on $\bf q$), and the van Kampen subdiagram $\Gamma$ bounded by a subpath of $\bf p$ (resp. $\bf q$), a subpath
of $\partial\Pi$, by $\cal C$ and $\cal C'$ contains no disks (but contains
the spokes $\cal C$  and $\cal C'$), then we shall call $\Gamma$ a
$tp$-{\it bond} at $\Pi$ (resp., $tq$-bond).

% This is a LaTeX picture output by TeXCAD.
% File name: [tp.pic].
% Version of TeXCAD: 4.3
% Reference / build: 30-Jun-2012 (rev. 105)
% For new versions, check: http://texcad.sf.net/
% Options on the following lines.
%\grade{\on}
%\emlines{\off}
%\epic{\off}
%\beziermacro{\on}
%\reduce{\on}
%\snapping{\off}
%\pvinsert{% Your \input, \def, etc. here}
%\quality{8.000}
%\graddiff{0.005}
%\snapasp{1}
%\zoom{4.0000}
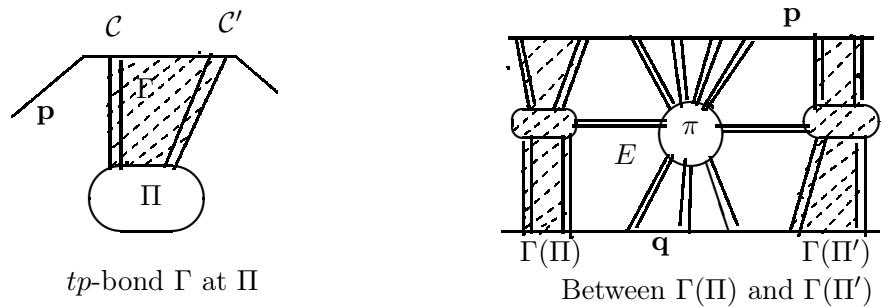
\unitlength 1mm % = 2.845pt
\linethickness{0.4pt}
\ifx\plotpoint\undefined\newsavebox{\plotpoint}\fi % GNUPLOT compatibility
\begin{picture}(129.25,80.75)(0,80)
\thicklines
\put(31.875,122.125){\oval(15.25,8.75)[]}
%\emline(14,133)(23.5,141)
\multiput(14,133)(.039915966,.033613445){238}{\line(1,0){.039915966}}
%\end
\put(23.5,141){\line(1,0){20.25}}
%\emline(43.75,141)(49.25,136.25)
\multiput(43.75,141)(.039007092,-.033687943){141}{\line(1,0){.039007092}}
%\end
\put(27,141){\line(0,-1){15}}
\put(28.5,140.5){\line(0,-1){14.5}}
%\emline(40.5,141.5)(34.25,126.25)
\multiput(40.5,141.5)(-.033602151,-.081989247){186}{\line(0,-1){.081989247}}
%\end
%\emline(42.5,141)(35.5,126)
\multiput(42.5,141)(-.033653846,-.072115385){208}{\line(0,-1){.072115385}}
%\end
\put(31,121.25){$\Pi$}
\put(26.5,143.75){$\cal C$}
\put(41.5,144.25){$\cal C'$}
\put(17.25,132.75){$\bf p$}
\put(30.5,135.5){$\Gamma$}
\thinlines
%\dashline{1}(28.25,139.25)(28.75,139)
\multiput(28.18,139.18)(.0625,-.03125){4}{\line(1,0){.0625}}
%\end
%\dashline{1}(26.75,138.25)(31,140.5)
\multiput(26.68,138.18)(.059028,.03125){12}{\line(1,0){.059028}}
\multiput(28.096,138.93)(.059028,.03125){12}{\line(1,0){.059028}}
\multiput(29.513,139.68)(.059028,.03125){12}{\line(1,0){.059028}}
%\end
%\dashline{1}(27,134.75)(26.75,134.5)
\multiput(26.93,134.68)(-.03125,-.03125){4}{\line(0,-1){.03125}}
%\end
%\dashline{1}(26.75,134.5)(27.25,135.25)
\multiput(26.68,134.43)(.03125,.046875){8}{\line(0,1){.046875}}
%\end
%\dashline{1}(27.5,136)(34.25,141)
\multiput(27.43,135.93)(.045,.0333333){15}{\line(1,0){.045}}
\multiput(28.78,136.93)(.045,.0333333){15}{\line(1,0){.045}}
\multiput(30.13,137.93)(.045,.0333333){15}{\line(1,0){.045}}
\multiput(31.48,138.93)(.045,.0333333){15}{\line(1,0){.045}}
\multiput(32.83,139.93)(.045,.0333333){15}{\line(1,0){.045}}
%\end
%\dashline{1}(26.75,133.25)(36.5,141)
\multiput(26.68,133.18)(.0416667,.0331197){18}{\line(1,0){.0416667}}
\multiput(28.18,134.372)(.0416667,.0331197){18}{\line(1,0){.0416667}}
\multiput(29.68,135.564)(.0416667,.0331197){18}{\line(1,0){.0416667}}
\multiput(31.18,136.757)(.0416667,.0331197){18}{\line(1,0){.0416667}}
\multiput(32.68,137.949)(.0416667,.0331197){18}{\line(1,0){.0416667}}
\multiput(34.18,139.141)(.0416667,.0331197){18}{\line(1,0){.0416667}}
\multiput(35.68,140.334)(.0416667,.0331197){18}{\line(1,0){.0416667}}
%\end
%\dashline{1}(27,131)(27,130.5)
\put(26.93,130.93){\line(0,-1){.25}}
%\end
%\dashline{1}(26.75,131)(39,141.5)
\multiput(26.68,130.93)(.0378086,.0324074){18}{\line(1,0){.0378086}}
\multiput(28.041,132.096)(.0378086,.0324074){18}{\line(1,0){.0378086}}
\multiput(29.402,133.263)(.0378086,.0324074){18}{\line(1,0){.0378086}}
\multiput(30.763,134.43)(.0378086,.0324074){18}{\line(1,0){.0378086}}
\multiput(32.124,135.596)(.0378086,.0324074){18}{\line(1,0){.0378086}}
\multiput(33.485,136.763)(.0378086,.0324074){18}{\line(1,0){.0378086}}
\multiput(34.846,137.93)(.0378086,.0324074){18}{\line(1,0){.0378086}}
\multiput(36.207,139.096)(.0378086,.0324074){18}{\line(1,0){.0378086}}
\multiput(37.569,140.263)(.0378086,.0324074){18}{\line(1,0){.0378086}}
%\end
%\dashline{1}(39,141.5)(39,142)
\put(38.93,141.43){\line(0,1){.25}}
%\end
%\dashline{1}(27.5,129)(41.25,141)
\multiput(27.43,128.93)(.0381944,.0333333){18}{\line(1,0){.0381944}}
\multiput(28.805,130.13)(.0381944,.0333333){18}{\line(1,0){.0381944}}
\multiput(30.18,131.33)(.0381944,.0333333){18}{\line(1,0){.0381944}}
\multiput(31.555,132.53)(.0381944,.0333333){18}{\line(1,0){.0381944}}
\multiput(32.93,133.73)(.0381944,.0333333){18}{\line(1,0){.0381944}}
\multiput(34.305,134.93)(.0381944,.0333333){18}{\line(1,0){.0381944}}
\multiput(35.68,136.13)(.0381944,.0333333){18}{\line(1,0){.0381944}}
\multiput(37.055,137.33)(.0381944,.0333333){18}{\line(1,0){.0381944}}
\multiput(38.43,138.53)(.0381944,.0333333){18}{\line(1,0){.0381944}}
\multiput(39.805,139.73)(.0381944,.0333333){18}{\line(1,0){.0381944}}
%\end
%\dashline{1}(27.25,126.5)(41.25,138.75)
\multiput(27.18,126.43)(.0368421,.0322368){19}{\line(1,0){.0368421}}
\multiput(28.58,127.655)(.0368421,.0322368){19}{\line(1,0){.0368421}}
\multiput(29.98,128.88)(.0368421,.0322368){19}{\line(1,0){.0368421}}
\multiput(31.38,130.105)(.0368421,.0322368){19}{\line(1,0){.0368421}}
\multiput(32.78,131.33)(.0368421,.0322368){19}{\line(1,0){.0368421}}
\multiput(34.18,132.555)(.0368421,.0322368){19}{\line(1,0){.0368421}}
\multiput(35.58,133.78)(.0368421,.0322368){19}{\line(1,0){.0368421}}
\multiput(36.98,135.005)(.0368421,.0322368){19}{\line(1,0){.0368421}}
\multiput(38.38,136.23)(.0368421,.0322368){19}{\line(1,0){.0368421}}
\multiput(39.78,137.455)(.0368421,.0322368){19}{\line(1,0){.0368421}}
%\end
%\dashline{1}(30,126.5)(28.5,126.75)
\multiput(29.93,126.43)(-.1875,.03125){4}{\line(-1,0){.1875}}
%\end
%\dashline{1}(28.5,126.75)(27.25,126.5)
\multiput(28.43,126.68)(-.15625,-.03125){4}{\line(-1,0){.15625}}
%\end
%\dashline{1}(29.5,127)(28.75,127)
\put(29.43,126.93){\line(-1,0){.375}}
%\end
%\dashline{1}(29.5,126.5)(39.25,134.25)
\multiput(29.43,126.43)(.0416667,.0331197){18}{\line(1,0){.0416667}}
\multiput(30.93,127.622)(.0416667,.0331197){18}{\line(1,0){.0416667}}
\multiput(32.43,128.814)(.0416667,.0331197){18}{\line(1,0){.0416667}}
\multiput(33.93,130.007)(.0416667,.0331197){18}{\line(1,0){.0416667}}
\multiput(35.43,131.199)(.0416667,.0331197){18}{\line(1,0){.0416667}}
\multiput(36.93,132.391)(.0416667,.0331197){18}{\line(1,0){.0416667}}
\multiput(38.43,133.584)(.0416667,.0331197){18}{\line(1,0){.0416667}}
%\end
%\dashline{1}(31.75,126.75)(37.25,130.25)
\multiput(31.68,126.68)(.052381,.0333333){15}{\line(1,0){.052381}}
\multiput(33.251,127.68)(.052381,.0333333){15}{\line(1,0){.052381}}
\multiput(34.823,128.68)(.052381,.0333333){15}{\line(1,0){.052381}}
\multiput(36.394,129.68)(.052381,.0333333){15}{\line(1,0){.052381}}
%\end
\thicklines
%\dashline{1}(80,139.75)(80.25,139)
\multiput(79.93,139.68)(.03125,-.09375){4}{\line(0,-1){.09375}}
%\end
\put(80,143.5){\line(1,0){48.5}}
\put(79,117.75){\line(1,0){50.25}}
\put(84.75,132.125){\oval(8.5,3.75)[]}
\put(124.25,132.375){\oval(10,4.25)[]}
%\emline(81,143.5)(82.75,134.25)
\multiput(81,143.5)(.03365385,-.17788462){52}{\line(0,-1){.17788462}}
%\end
%\emline(82.75,134.25)(82.5,134.5)
\multiput(82.75,134.25)(-.03125,.03125){8}{\line(0,1){.03125}}
%\end
%\emline(81.75,143.25)(83.5,134.75)
\multiput(81.75,143.25)(.03365385,-.16346154){52}{\line(0,-1){.16346154}}
%\end
%\emline(90.5,143.25)(87.25,134.25)
\multiput(90.5,143.25)(-.033505155,-.092783505){97}{\line(0,-1){.092783505}}
%\end
%\emline(89.5,143.5)(86,134)
\multiput(89.5,143.5)(-.033653846,-.091346154){104}{\line(0,-1){.091346154}}
%\end
\put(82,129.75){\line(0,-1){12.5}}
%\emline(82,117.25)(82.25,117.75)
\multiput(82,117.25)(.03125,.0625){8}{\line(0,1){.0625}}
%\end
\put(83,130.5){\line(0,-1){13}}
\put(88.25,130.75){\line(0,-1){14}}
%\emline(88.25,116.75)(88.5,117)
\multiput(88.25,116.75)(.03125,.03125){8}{\line(0,1){.03125}}
%\end
\put(87.25,130.5){\line(0,-1){13}}
\put(87.25,117.5){\line(0,1){0}}
\put(120.75,144){\line(0,-1){10}}
\put(120.75,134){\line(0,1){0}}
\put(121.75,143.5){\line(0,-1){8.25}}
\put(126,143.75){\line(0,-1){9}}
\put(127,143.75){\line(0,-1){8.5}}
%\emline(121,130.25)(117.5,118)
\multiput(121,130.25)(-.033653846,-.117788462){104}{\line(0,-1){.117788462}}
%\end
%\emline(122.25,130.25)(118.5,117)
\multiput(122.25,130.25)(-.033482143,-.118303571){112}{\line(0,-1){.118303571}}
%\end
\put(126.25,129.5){\line(0,-1){11.25}}
\put(127.5,130.5){\line(0,-1){13.5}}
\thinlines
%\dashline{1}(81.75,140)(85.75,144)
\multiput(81.68,139.93)(.0333333,.0333333){15}{\line(0,1){.0333333}}
\multiput(82.68,140.93)(.0333333,.0333333){15}{\line(0,1){.0333333}}
\multiput(83.68,141.93)(.0333333,.0333333){15}{\line(0,1){.0333333}}
\multiput(84.68,142.93)(.0333333,.0333333){15}{\line(0,1){.0333333}}
%\end
%\dashline{1}(82.5,137.75)(89.25,143.75)
\multiput(82.43,137.68)(.0375,.0333333){18}{\line(1,0){.0375}}
\multiput(83.78,138.88)(.0375,.0333333){18}{\line(1,0){.0375}}
\multiput(85.13,140.08)(.0375,.0333333){18}{\line(1,0){.0375}}
\multiput(86.48,141.28)(.0375,.0333333){18}{\line(1,0){.0375}}
\multiput(87.83,142.48)(.0375,.0333333){18}{\line(1,0){.0375}}
%\end
%\dashline{1}(83,134)(89.25,140.75)
\multiput(82.93,133.93)(.0328947,.0355263){19}{\line(0,1){.0355263}}
\multiput(84.18,135.28)(.0328947,.0355263){19}{\line(0,1){.0355263}}
\multiput(85.43,136.63)(.0328947,.0355263){19}{\line(0,1){.0355263}}
\multiput(86.68,137.98)(.0328947,.0355263){19}{\line(0,1){.0355263}}
\multiput(87.93,139.33)(.0328947,.0355263){19}{\line(0,1){.0355263}}
%\end
%\dashline{1}(80.75,131.25)(84.5,134.5)
\multiput(80.68,131.18)(.0367647,.0318627){17}{\line(1,0){.0367647}}
\multiput(81.93,132.263)(.0367647,.0318627){17}{\line(1,0){.0367647}}
\multiput(83.18,133.346)(.0367647,.0318627){17}{\line(1,0){.0367647}}
%\end
%\dashline{1}(82.25,130.75)(87.5,134.75)
\multiput(82.18,130.68)(.04375,.0333333){15}{\line(1,0){.04375}}
\multiput(83.492,131.68)(.04375,.0333333){15}{\line(1,0){.04375}}
\multiput(84.805,132.68)(.04375,.0333333){15}{\line(1,0){.04375}}
\multiput(86.117,133.68)(.04375,.0333333){15}{\line(1,0){.04375}}
%\end
%\dashline{1}(82,128.25)(88.75,133.75)
\multiput(81.93,128.18)(.0397059,.0323529){17}{\line(1,0){.0397059}}
\multiput(83.28,129.28)(.0397059,.0323529){17}{\line(1,0){.0397059}}
\multiput(84.63,130.38)(.0397059,.0323529){17}{\line(1,0){.0397059}}
\multiput(85.98,131.48)(.0397059,.0323529){17}{\line(1,0){.0397059}}
\multiput(87.33,132.58)(.0397059,.0323529){17}{\line(1,0){.0397059}}
%\end
%\dashline{1}(81.75,125.75)(88.75,131)
\multiput(81.68,125.68)(.04375,.0328125){16}{\line(1,0){.04375}}
\multiput(83.08,126.73)(.04375,.0328125){16}{\line(1,0){.04375}}
\multiput(84.48,127.78)(.04375,.0328125){16}{\line(1,0){.04375}}
\multiput(85.88,128.83)(.04375,.0328125){16}{\line(1,0){.04375}}
\multiput(87.28,129.88)(.04375,.0328125){16}{\line(1,0){.04375}}
%\end
%\dashline{1}(82,122.75)(88.5,128.75)
\multiput(81.93,122.68)(.0361111,.0333333){18}{\line(1,0){.0361111}}
\multiput(83.23,123.88)(.0361111,.0333333){18}{\line(1,0){.0361111}}
\multiput(84.53,125.08)(.0361111,.0333333){18}{\line(1,0){.0361111}}
\multiput(85.83,126.28)(.0361111,.0333333){18}{\line(1,0){.0361111}}
\multiput(87.13,127.48)(.0361111,.0333333){18}{\line(1,0){.0361111}}
%\end
%\dashline{1}(82.25,119.5)(82.25,120)
\put(82.18,119.43){\line(0,1){.25}}
%\end
%\dashline{1}(82.25,121)(82.5,120.5)
\multiput(82.18,120.93)(.03125,-.0625){4}{\line(0,-1){.0625}}
%\end
%\dashline{1}(82.25,120.5)(87.5,125)
\multiput(82.18,120.43)(.0386029,.0330882){17}{\line(1,0){.0386029}}
\multiput(83.492,121.555)(.0386029,.0330882){17}{\line(1,0){.0386029}}
\multiput(84.805,122.68)(.0386029,.0330882){17}{\line(1,0){.0386029}}
\multiput(86.117,123.805)(.0386029,.0330882){17}{\line(1,0){.0386029}}
%\end
%\dashline{1}(82.25,118.25)(87.5,122)
\multiput(82.18,118.18)(.046875,.0334821){16}{\line(1,0){.046875}}
\multiput(83.68,119.251)(.046875,.0334821){16}{\line(1,0){.046875}}
\multiput(85.18,120.323)(.046875,.0334821){16}{\line(1,0){.046875}}
\multiput(86.68,121.394)(.046875,.0334821){16}{\line(1,0){.046875}}
%\end
%\dashline{1}(85,117.75)(85.5,118.5)
\multiput(84.93,117.68)(.03125,.046875){8}{\line(0,1){.046875}}
%\end
%\dashline{1}(85.5,118.5)(84.25,117.5)
\multiput(85.43,118.43)(-.041667,-.033333){10}{\line(-1,0){.041667}}
\multiput(84.596,117.763)(-.041667,-.033333){10}{\line(-1,0){.041667}}
%\end
%\dashline{1}(84.5,118.5)(87.75,120.25)
\multiput(84.43,118.43)(.059091,.031818){11}{\line(1,0){.059091}}
\multiput(85.73,119.13)(.059091,.031818){11}{\line(1,0){.059091}}
\multiput(87.03,119.83)(.059091,.031818){11}{\line(1,0){.059091}}
%\end
%\dashline{1}(120.75,140.25)(125.75,144)
\multiput(120.68,140.18)(.0446429,.0334821){16}{\line(1,0){.0446429}}
\multiput(122.108,141.251)(.0446429,.0334821){16}{\line(1,0){.0446429}}
\multiput(123.537,142.323)(.0446429,.0334821){16}{\line(1,0){.0446429}}
\multiput(124.965,143.394)(.0446429,.0334821){16}{\line(1,0){.0446429}}
%\end
%\dashline{1}(121.5,137.25)(126.75,142.75)
\multiput(121.43,137.18)(.0324074,.0339506){18}{\line(0,1){.0339506}}
\multiput(122.596,138.402)(.0324074,.0339506){18}{\line(0,1){.0339506}}
\multiput(123.763,139.624)(.0324074,.0339506){18}{\line(0,1){.0339506}}
\multiput(124.93,140.846)(.0324074,.0339506){18}{\line(0,1){.0339506}}
\multiput(126.096,142.069)(.0324074,.0339506){18}{\line(0,1){.0339506}}
%\end
%\dashline{1}(121,133.75)(126.75,138.75)
\multiput(120.93,133.68)(.0375817,.0326797){17}{\line(1,0){.0375817}}
\multiput(122.207,134.791)(.0375817,.0326797){17}{\line(1,0){.0375817}}
\multiput(123.485,135.902)(.0375817,.0326797){17}{\line(1,0){.0375817}}
\multiput(124.763,137.013)(.0375817,.0326797){17}{\line(1,0){.0375817}}
\multiput(126.041,138.124)(.0375817,.0326797){17}{\line(1,0){.0375817}}
%\end
%\dashline{1}(126.75,138.75)(126,138.5)
\multiput(126.68,138.68)(-.09375,-.03125){4}{\line(-1,0){.09375}}
%\end
%\dashline{1}(126,138.5)(127.5,138)
\multiput(125.93,138.43)(.09375,-.03125){8}{\line(1,0){.09375}}
%\end
%\dashline{1}(120.25,131)(127,136.75)
\multiput(120.18,130.93)(.0375,.0319444){18}{\line(1,0){.0375}}
\multiput(121.53,132.08)(.0375,.0319444){18}{\line(1,0){.0375}}
\multiput(122.88,133.23)(.0375,.0319444){18}{\line(1,0){.0375}}
\multiput(124.23,134.38)(.0375,.0319444){18}{\line(1,0){.0375}}
\multiput(125.58,135.53)(.0375,.0319444){18}{\line(1,0){.0375}}
%\end
%\dashline{1}(121.25,129.25)(126.75,135)
\multiput(121.18,129.18)(.0321637,.0336257){19}{\line(0,1){.0336257}}
\multiput(122.402,130.457)(.0321637,.0336257){19}{\line(0,1){.0336257}}
\multiput(123.624,131.735)(.0321637,.0336257){19}{\line(0,1){.0336257}}
\multiput(124.846,133.013)(.0321637,.0336257){19}{\line(0,1){.0336257}}
\multiput(126.069,134.291)(.0321637,.0336257){19}{\line(0,1){.0336257}}
%\end
%\dashline{1}(119.5,126)(128.75,133)
\multiput(119.43,125.93)(.0444712,.0336538){16}{\line(1,0){.0444712}}
\multiput(120.853,127.007)(.0444712,.0336538){16}{\line(1,0){.0444712}}
\multiput(122.276,128.084)(.0444712,.0336538){16}{\line(1,0){.0444712}}
\multiput(123.699,129.16)(.0444712,.0336538){16}{\line(1,0){.0444712}}
\multiput(125.122,130.237)(.0444712,.0336538){16}{\line(1,0){.0444712}}
\multiput(126.545,131.314)(.0444712,.0336538){16}{\line(1,0){.0444712}}
\multiput(127.968,132.391)(.0444712,.0336538){16}{\line(1,0){.0444712}}
%\end
%\dashline{1}(119.25,122.25)(118.5,122.75)
\multiput(119.18,122.18)(-.046875,.03125){8}{\line(-1,0){.046875}}
%\end
%\dashline{1}(118.5,122.75)(127.25,129.25)
\multiput(118.43,122.68)(.0428922,.0318627){17}{\line(1,0){.0428922}}
\multiput(119.888,123.763)(.0428922,.0318627){17}{\line(1,0){.0428922}}
\multiput(121.346,124.846)(.0428922,.0318627){17}{\line(1,0){.0428922}}
\multiput(122.805,125.93)(.0428922,.0318627){17}{\line(1,0){.0428922}}
\multiput(124.263,127.013)(.0428922,.0318627){17}{\line(1,0){.0428922}}
\multiput(125.721,128.096)(.0428922,.0318627){17}{\line(1,0){.0428922}}
%\end
%\dashline{1}(117.75,118.75)(127,126.75)
\multiput(117.68,118.68)(.0388655,.0336134){17}{\line(1,0){.0388655}}
\multiput(119.001,119.823)(.0388655,.0336134){17}{\line(1,0){.0388655}}
\multiput(120.323,120.965)(.0388655,.0336134){17}{\line(1,0){.0388655}}
\multiput(121.644,122.108)(.0388655,.0336134){17}{\line(1,0){.0388655}}
\multiput(122.965,123.251)(.0388655,.0336134){17}{\line(1,0){.0388655}}
\multiput(124.287,124.394)(.0388655,.0336134){17}{\line(1,0){.0388655}}
\multiput(125.608,125.537)(.0388655,.0336134){17}{\line(1,0){.0388655}}
%\end
%\dashline{1}(119.75,118.5)(126.25,124.5)
\multiput(119.68,118.43)(.0361111,.0333333){18}{\line(1,0){.0361111}}
\multiput(120.98,119.63)(.0361111,.0333333){18}{\line(1,0){.0361111}}
\multiput(122.28,120.83)(.0361111,.0333333){18}{\line(1,0){.0361111}}
\multiput(123.58,122.03)(.0361111,.0333333){18}{\line(1,0){.0361111}}
\multiput(124.88,123.23)(.0361111,.0333333){18}{\line(1,0){.0361111}}
%\end
%\dashline{1}(121.75,117.5)(126.75,123.25)
\multiput(121.68,117.43)(.0326797,.0375817){17}{\line(0,1){.0375817}}
\multiput(122.791,118.707)(.0326797,.0375817){17}{\line(0,1){.0375817}}
\multiput(123.902,119.985)(.0326797,.0375817){17}{\line(0,1){.0375817}}
\multiput(125.013,121.263)(.0326797,.0375817){17}{\line(0,1){.0375817}}
\multiput(126.124,122.541)(.0326797,.0375817){17}{\line(0,1){.0375817}}
%\end
%\dashline{1}(125,118.5)(127,120.25)
\multiput(124.93,118.43)(.0384615,.0336538){13}{\line(1,0){.0384615}}
\multiput(125.93,119.305)(.0384615,.0336538){13}{\line(1,0){.0384615}}
%\end
\thicklines
\put(104.25,130.75){\circle{8.016}}
%\dashline{1}(95.75,143.5)(96,143.75)
\multiput(95.68,143.43)(.03125,.03125){4}{\line(0,1){.03125}}
%\end
%\emline(96.25,143.75)(101.5,134)
\multiput(96.25,143.75)(.033653846,-.0625){156}{\line(0,-1){.0625}}
%\end
%\emline(97.25,143.5)(102,134.75)
\multiput(97.25,143.5)(.033687943,-.062056738){141}{\line(0,-1){.062056738}}
%\end
%\emline(101.75,143.5)(103,134.25)
\multiput(101.75,143.5)(.03289474,-.24342105){38}{\line(0,-1){.24342105}}
%\end
%\emline(103,143.25)(103.75,135.25)
\multiput(103,143.25)(.0326087,-.3478261){23}{\line(0,-1){.3478261}}
%\end
%\emline(107,143.5)(104.5,134.75)
\multiput(107,143.5)(-.03333333,-.11666667){75}{\line(0,-1){.11666667}}
%\end
\put(104.5,134.75){\line(0,1){0}}
%\emline(108.25,143.25)(105,134.5)
\multiput(108.25,143.25)(-.033505155,-.090206186){97}{\line(0,-1){.090206186}}
%\end
\put(105,134.5){\line(0,1){.25}}
%\emline(111.5,143.5)(105.75,134)
\multiput(111.5,143.5)(-.033625731,-.055555556){171}{\line(0,-1){.055555556}}
%\end
\put(105.75,134){\line(1,0){.75}}
%\emline(112.5,143)(106.5,134.25)
\multiput(112.5,143)(-.033707865,-.049157303){178}{\line(0,-1){.049157303}}
%\end
%\emline(106.5,134.25)(106.25,135)
\multiput(106.5,134.25)(-.03125,.09375){8}{\line(0,1){.09375}}
%\end
%\emline(106.25,135)(107.25,134.25)
\multiput(106.25,135)(.0434783,-.0326087){23}{\line(1,0){.0434783}}
%\end
%\emline(101.25,128)(95.75,117.75)
\multiput(101.25,128)(-.033536585,-.0625){164}{\line(0,-1){.0625}}
%\end
\put(95.75,117.75){\line(0,1){.25}}
\put(95.75,118){\line(0,1){0}}
\put(101.75,128){\line(0,-1){.5}}
%\emline(101.75,127.5)(97,118.25)
\multiput(101.75,127.5)(-.033687943,-.065602837){141}{\line(0,-1){.065602837}}
%\end
%\emline(103.5,127.5)(102.75,118.5)
\multiput(103.5,127.5)(-.0326087,-.3913043){23}{\line(0,-1){.3913043}}
%\end
\put(102.75,118.5){\line(0,1){0}}
%\emline(104.25,126.5)(104,117.5)
\multiput(104.25,126.5)(-.03125,-1.125){8}{\line(0,-1){1.125}}
%\end
\put(106.25,127){\line(1,-3){3}}
%\emline(106.75,127.75)(110.5,118.5)
\multiput(106.75,127.75)(.033482143,-.082589286){112}{\line(0,-1){.082589286}}
%\end
\put(89.25,132.5){\line(1,0){11.75}}
\put(88.5,131.75){\line(1,0){12.5}}
\put(107.5,131.75){\line(1,0){12.5}}
\put(108.5,131){\line(1,0){11.5}}
\put(120,131){\line(0,1){0}}
\put(81.5,113.75){$\Gamma(\Pi)$}
\put(119,113.75){$\Gamma(\Pi')$}
\put(94,126.75){$E$}
\put(103,130.75){$\pi$}
\put(116.5,145.75){$\bf p$}
\put(99,115.25){$\bf q$}
\put(21.25,110){$tp$-bond  $\Gamma$ at $\Pi$}
\put(87,108.75){Between $\Gamma(\Pi)$   and  $\Gamma(\Pi')$}
\end{picture}

\begin{lemma} \label{form}  For every disk $\Pi$ of $\Delta$,
the number of consecutive $tp$-bonds ($tq$-bonds) is greater than
$L/2 -4>0$.

No $q$- or $\theta$-band of $\Delta$ starts and ends on $\bf p$ or starts and ends on $\bf q$.
\end{lemma}

\proof If a $\theta$-band connects $\bf p$ with $\bf p$, then
a van Kampen subdiagram $\Gamma$ of $\Delta$
 is bounded by $\cal T$
and by a subpath $\bf p'$ of the cyclic path $\bf p$, and the word $U$ must contain
both positive and negative occurrences of $\theta$-letters. Proving by contradiction, we may assume that $\cal T$ is a unique maximal $\theta$-band of $\Gamma$
since every maximal $\theta$-band of $\Gamma$
has to start and end on $\bf p'$. Therefore the word $W\equiv Lab({\bf p'})$ has no $\theta$-letters except for the first and the last ones. Hence $W$ is a subword of a cyclic permutation of the word $U$. However such a diagram $\Gamma$ impossible since $U$ is an adapted word.

Assume now that there is a $q$-band $\cal C$ starting and ending on $\bf p$. Then there is a van Kampen diagram $\Gamma$ bounded by $\cal C$ and a subpath $\bf p'$ of $\bf p$. Proving by contradiction, we may assume that $\Gamma$ has no maximal $q$-bands stating and ending on $\bf p$ except for $\cal C$.  Suppose $\Gamma$ has disks. Since the sides of $\cal C$ have no $q$-edges, the disk $\Pi$ provided by Lemma \ref{extdisk} (2) applied to $\Gamma$, has at least $L-4$ consecutive $tp$-bonds. Since there are no $\tt$-letters between the ends
of the $\tt$-spokes defining these bonds, there is subpath $\bf x$ of $\bf p$ containing exactly $L-3$ $\tt$-edges, namely, the ends of the
$\tt$-spokes defining the $tp$-bonds. However $U$ has at least $L$
$\tt$-letters by the equalities (\ref{0i0}). Hence $Lab ({\bf x})$ is a subword of a cyclic permutation of the adapted word $U$, which contradicts
Lemma \ref{simpl} (a) since $L-3> L/2+1$. If $\Gamma$ contains no disks, the  argument from the previous paragraph leads to contradiction again.

Let us prove the first statement of the lemma. Under the assumption that $\Delta$ contains disks, Lemma \ref{extdisk} (2)
gives a disk $\Pi$ with  $L'+L''\ge L-3$
$\tt$-spokes ended on $\bf p$ and $\bf q$. If, say,
$L''<L/2 -3$, then
$L'> L/2$. However $U$ has at least $L$ $\tt$-letters by the equalities (\ref{0i0}). So by Property (b) of Lemma \ref{simpl} for an adapted word $U$, these $\tt$-spokes end on a subpath of $\bf p$ labeled by a subword of a cyclic permutation of $U$, a contradiction with Property
(a) of adapted words in Lemma \ref{simpl}.
Thus, we have
$L''\ge L/2 - 3$ and similarly, $L'\ge L-3>0$.

Now consider
a maximal set $\bf D$ of disks $\Pi$ having at least $L-5$ $\tt$-spokes ending either on ${\bf p}$ or on ${\bf q}$. Note that two neighbor $\tt$-spokes ending on ${\bf p}$ (or on ${\bf q}$) define a $tp$-bond at $\Pi$, i.e. the subdiagram $\Gamma$ bounded by these
spokes, by a part of $\bf p$ and a part of $\partial\Pi$ containing no disks, because
otherwise Lemma \ref{extdisk} (1) applied to $\Gamma$ should gives us another disks
with at least $L-4$ $tp$-bonds, which is impossible again since $L-4>L/2$.

We claim that all disks of $\Delta$ belong to $\bf D$. Indeed, let $\Gamma(\Pi)$ be a subdiagram formed by a disk $\Pi$ from $\bf D$ and all the $tp$- and $tq$-bonds at $\Pi$. Consider the maximal van Kampen subdiagram $E$ between neighbor $\Gamma(\Pi)$ and $\Gamma(\Pi')$. If $E$ contains a disk, then it has a disk $\pi$ provided by Lemma \ref{extdisk} (1). It has at least $L-3$ $\tt$-spokes in $E$. But the number of its
spokes ending either on ${\bf p}$ or on ${\bf q}$ is less than $L-5$ since $\pi$ does not belong to $\bf D$. It follows that a pair of $\tt$-spokes
connects $\pi$ with $\Pi$ or with $\Pi'$ by at least two $\tt$-spokes in a van Kampen subdiagram, which is impossible by Lemma \ref{2dis}. Thus, every disk $\Pi$ of $\Delta$ has to belong to the set $\bf D$.

The number of $\tt$ spokes at a disk $\Pi$ from $\bf D$ is $L$, and at most two
$\tt$-spokes connect it with neighbor disks from $\bf D$. So there are at least $(L-2)-2=L-4$ $tp$- and $tq$-bonds at $\Pi$. As above, the number of the $tp$-bonds at $\Pi$ is
less than $L/2$, whence the number of $tq$-bonds of it is greater than $L-4-L/2 = L/2-4>0$. Similarly, there are $>L/2 -4$ $tp$-bonds at $\Pi$; and the lemma is proved.
\endproof

\begin{lemma} \label{noq} There is a recursive function $f$ such that $k$ and $l$
do not exceed $f(||U||+||V||)$, provided the words $U$ and $V$ have no $q$-letters.
\end{lemma}

\proof By Lemma \ref{extdisk} (2) the annular diagram $\Delta$ contains no disks, and so it is a roll.
 Assume first that $\Delta$ has no $(\theta,q)$-cells.

{\bf 1.} If  the words $U,V$ have no $\theta$-letters, then by Lemmas \ref{form} and \ref{NoAnnul}, every maximal $\theta$-band $\cal T$ of $\Delta$
is an annulus surrounding the hole, and has side labels of the form
$(U')^k$, where $U'$ is $U$ or, due to a superscript, a copy of $U$.
This obviously bounds the number of such labels in terms of $||U||$,
and since we may assume that different $\theta$-annuli do not copy each other, every vertex of $\bf p$ can be connected with a vertex of
$\bf q$ by a  path of bounded length. If two such pathes $\bf x_1$ and $\bf x_2$ define a van Kampen subdiagram with boundary $\bf x_1y_1x_2^{-1}y_2^{-1}$, where $Lab ({\bf x_1})\equiv Lab ({\bf x_2})$,
$Lab({\bf y_1})=V^{l'}$, $Lab({\bf y_2})=U^{k'}$ with $|k'|<k$ and $|l'|<l$, then we obtain a contradiction with the choice of $k$ and $l$.
But the absence of pairs of such cuts $\bf x_1, x_2$ bounds the exponents $k$ and $l$ in terms of $||U||+||V||$ since the labels of such cuts belong to the bounded set.

The dual argument works if the words $U, V$ have no $Y$-letters.

We may also assume that $\Delta$ has no $Y$-bands starting and terminating on ${\bf p}$ (on $\bf q$). Indeed otherwise there is a rim $a$-band, and removing it we replace $U$ (or $V$) with a a conjugate word $\bar U$, such that $|\bar U|_a<|U|_a$ and $|\bar U|_{\theta}=|U|_{\theta}$; this replacement is effective.

{\bf 2.}  As in item {\bf 1}, $\Delta$ has no $(\theta,q)$-cells, but now $U$ contains $\theta$-letters. By Lemma \ref{form}, it remains to assume that every maximal $\theta$-band of $\Delta$ connects the
contours $\bf p$ and $\bf q$. The same is true for maximal $Y$-bands as we notices in the previous paragraph. It follows that $k|U|_{Y}=|U^k|_{Y}=|V^l|_{Y}= l|V|_{Y}$ and therefore
\begin{equation}\label{kl}
\frac{|U|_{Y}}{|V|_{Y}} =\frac lk
\end{equation}
Since the numbers $|U|_{Y}$ and $|V|_{Y}$ are less than $||U||+||V||$, it follows from (\ref{kl}) that $k$ and $l$ have a common
divisor $d$ such that $k=dk'$, $l=dl'$, where $k',l' <||U||+||V||$.
If $d=1$ then $k, l <||U||+||V||$, i.e. we get a desired upper bound.
Proving by contradiction, we assume now that $d>1$.
The words $U'\equiv U^{k'}$ and $V' \equiv V^{l'}$ have equal %$\theta$-lengths (
equal $Y$-length since their $d$-th powers $U^k$ and $V^l$ have equal  $Y$-length.

Without loss of generality, we may assume that $U$ starts with an $Y$-letter $a$. Let $\bf p'$ be a subpath of $\bf p$ labeled by
$U'a$. Let us denote by ${\cal T}_1$ and ${\cal T}_2$ the maximal $Y$-bands starting with the first and the last edges of $\bf p'$.
They end on $\bf q$, and we get a van Kampen subdiagram $\Gamma$
bounded by $\bf p'$, by a subpath $\bf q'$ of $\bf q$ and by the sides of ${\cal T}_1$
and ${\cal T}_2$.  Since all maximal $Y$-bands of $\Gamma$ connect
$\bf p'$ with $\bf q'$, we have $|\bf q'|_{a}=|\bf p'|_{a}$,
and so $Lab ({\bf q'}) \equiv V'a$ (Here  we may replace the word $V'$ with a cyclic permutation of it.)

The boundary label of $\Gamma\backslash {\cal T}_2$ is $T_1V'T_2^{-1}(U')^{-1}$, where $T_1$ and $T_2$ are side labels of
${\cal T}_1$ and ${\cal T}_2$, resp., and so we obtain
\begin{equation}\label{Ga}
T_2 = (U')^{-1}T_1V'
\end{equation}
in $G$.  Also, cutting $\Delta$ along the side of
${\cal T}_1$, we have  in $G$:
\begin{equation}\label{con12}
T_1^{-1}(U')^dT_1 = (V')^d
\end{equation}
because the paths $\bf p$ and $\bf q$ are labeled by $(U')^d$ and $(V')^d$, resp.

Both diagram $\Gamma$ and $\Delta$ contains only $(\theta,a)$-cells, and they are diagrams over the group $G_{\theta a}$ generated by $Y$-letters and $\theta$-letters only. The form of the $(\theta,a)$-relations
implies the existence of the homomorphism $\nu$ of $G_{\theta a}$ onto
the free group $F$ generated by $\theta$-letters: $\nu$ is identical
on $\theta$-letters and trivial on $Y$-letters. On the one hand, the equality (\ref{con12}) gives us
$$(\nu(T_1)^{-1}\nu(U')\nu(T_1))^d = \nu(V')^d,$$ which implies
$\nu(T_1)^{-1}\nu(U')\nu(T_1) = \nu(V')$ in the free group $F$, and so\\
$\nu(T_1) =\nu(U')^{-1}\nu(T_1)\nu(V')$. On the other hand, we get
$\nu(T_2) =\nu(U')^{-1}\nu(T_1)\nu(V')$ from (\ref{Ga}).
Therefore we have $\nu(T_2)=\nu(T_1)$ in $F$, whence $T_2=T_1$ in $G$, because $T_1$ and $T_2$ contain only $\theta$-letters,
Now the equalities $T_2=T_1$ and (\ref{Ga}) gives us the conjugation of the words $U'=U^{k'}$ and $V' = V^{l'}$
in $G$, where $k'<k$ and $l'<l$, which contradicts
the choice of the pair $(k,l)$.

{\bf 3.} If $\Delta$ contains $(\theta,q)$-cells, then by Lemma \ref{NoAnnul},
 every maximal $q$-band of $\Delta$ is a $q$-annulus surrounding the hole of $\Delta$ since a roll has no $q$-edges in the boundary.
 By the same lemma, every maximal $\theta$-band crossing a $q$-annulus
 $\cal C$ connects  $\bf p$ and $\bf q$ and cannot intersect $\cal C$ twice. Therefore all $q$-annuli have length $|U|_{\theta}^k$, and the boundary label of each of them is a  $k$-th power with the length of base bounded
 from above by $||U||$. Since one may assume that two different
 $q$-annuli do not copy each other, the number of $q$-annuli is effectively bounded. Therefore the solution
 of power conjugation is reduced to the annular diagrams between
 the annuli, where there are no $(\theta,q)$-cells. Since the number
 of such annular diagrams is bounded, the problem is reduced to
 the case considered in item {\bf 2}, because one can use the transitiveness: if $U^k$ is a conjugate of $V^l$ and $V^r$ is a conjugate of $W^s$, then $U^{kr}$ is a conjugate of $W^{sl}$.

 The lemma is proved. \endproof

 For any disk $\Pi$ of the diagram $\Delta$, we have a $tq$-bond $\Gamma$ at $\Pi$ by Lemma \ref{form}, because $L/2-4 >0$. If there is
 a $\theta$-band of $\Gamma$ connecting the two spokes bounding $\Gamma$,
 then there is such a $\theta$-band $\cal T$ closest to $\Pi$.
 Let $E$ be the subdiagram formed by $\Pi$ and $\cal T$. One may
 apply Lemma \ref{moving} (1)  and replace $E$ with a diagram $E'$
 formed by a new disk $\Pi'$ and a $\theta$-band $\cal T'$. This
 transformation
 replaces $\Gamma$ with the $tq$-bond $\Gamma'=\Gamma\backslash\cal T$ at $\Pi'$. The iteration of such transformation replaces
 the $tq$-bond $\Gamma$ with a $tq$-bond $\Gamma_0$ at a disk $\Pi_0$, where there are no $\theta$-bands connecting $\tt$-spokes $\cal C$ and $\cal C'$ at $\Pi_0$.

 \begin{lemma}\label{dob} The perimeter of $\Pi_0$, the lengths of
 the $\tt$-spokes $\cal C$, $\cal C'$, and the length $||{\bf r}||$ of some path ${\bf r}$ of $\theta$-length $0$ connecting $\Pi_0$ and $\bf q$   in $\Gamma_0$ are effectively
 bounded from above in terms of $||V||$.
 \end{lemma}

 \proof The quadratic upper bounds for the lengths $\cal C$, $\cal C'$,
 and $\partial\Pi_0$ in terms of the length of the subpath
 $\bf x $ of $\bf q$ connecting $\cal C$ and $\cal C'$  is given by   Lemma \ref{malo}. However we have  $||{\bf x}||<||V||$
 since the equality $\mu(V)=0$ implies that the word $V$ contains at least $L$ $\tt$-letters, but ${\bf x}$ has only $2$ $\tt$-edges by Lemma \ref{form}.  It remains to define the path $\bf r$. This path starts from $\Pi_0$, where the $q$-band $\cal C$ starts, but
 it is a side of a maximal $\theta$-band ${\cal T}_0$ of $\Gamma_0$. ${\cal T}_0$ must end on $\bf q$ by the definition of $\Gamma_0$,
 The length of $\bf r$ is bounded by Lemma \ref{malo} since the perimeter of $\Gamma_0$ is bounded and the number of cells in $\Gamma_0$ is also effectively bounded by Lemma \ref{quadr}.
 \endproof

 By Lemma \ref{form}, all disks of $\Delta$ can be moved toward $\bf q$ in the same  way we have moved $\Pi$. So we obtain an annular diagram $\tilde\Delta$, where by Lemma \ref{dob}, each disk $\Pi$ has effectively bounded perimeter and
 connected with $\bf q$ by a path ${\bf r = r}(\Pi)$ having linealy bounded length and $|{\bf r}|_{\theta}=0$.
 The obtained annular diagram $\tilde\Delta$ has the same boundary labels as $\Delta$, but it is not necessarily minimal. Every disk $\Pi$ can be removed from $\tilde\Delta$ if one makes the cut along ${\bf r}^{-1}$,
 around $\Pi$ and back along $\bf r$. After removal of all the disks,
 we obtain a diskless annular diagram $\Delta_0$.

 We may keep notation $\bf p$ and $\bf q$ for the boundary components of $\tilde\Delta$, where $\bf p$ is also
 the outer boundary component of $\Delta_0$. If $\Delta_0$ is not reduced, we replace
 it with a reduced annular diagram with the same boundary labels. So we will assume that $\Delta_0$ is a reduced annular diagram and $\tilde\Delta$ is built of $\Delta_0$ and disks. The inner contour $q_0$ of $\Delta_0$ is obtained from $\bf q$ by inserting pathes ${\bf z=z}(\Pi)$ for every disk $\Pi$, where $|{\bf z}|_{\theta}=0$ and the length $||{\bf z}||$ is effectively bounded in terms of $||V||$.

 % This is a LaTeX picture output by TeXCAD.
% File name: [delta0.pic].
% Version of TeXCAD: 4.3
% Reference / build: 30-Jun-2012 (rev. 105)
% For new versions, check: http://texcad.sf.net/
% Options on the following lines.
%\grade{\on}
%\emlines{\off}
%\epic{\off}
%\beziermacro{\on}
%\reduce{\on}
%\snapping{\off}
%\pvinsert{% Your \input, \def, etc. here}
%\quality{8.000}
%\graddiff{0.005}
%\snapasp{1}
%\zoom{4.0000}
\unitlength 1mm % = 2.845pt
\linethickness{0.4pt}
\ifx\plotpoint\undefined\newsavebox{\plotpoint}\fi % GNUPLOT compatibility
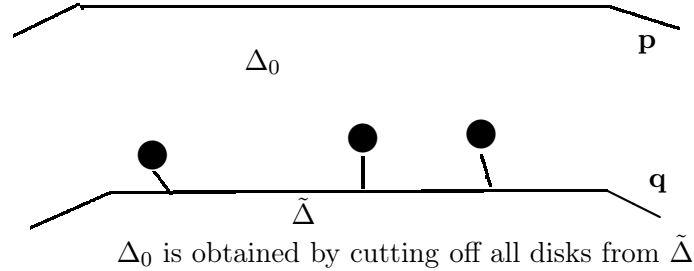
\begin{picture}(108.5,45.5)(0,20)
%\emline(29.25,54.75)(29,55.25)
\multiput(29.25,54.75)(-.03125,.0625){8}{\line(0,1){.0625}}
%\end
\thicklines
\put(29,55.25){\line(1,0){70}}
%\emline(99,55.25)(108.5,52.25)
\multiput(99,55.25)(.10674157,-.03370787){89}{\line(1,0){.10674157}}
%\end
%\emline(28.75,55.5)(20,51.25)
\multiput(28.75,55.5)(-.069444444,-.033730159){126}{\line(-1,0){.069444444}}
%\end
%\emline(22.25,25.75)(33,30.5)
\multiput(22.25,25.75)(.076241135,.033687943){141}{\line(1,0){.076241135}}
%\end
%\emline(33,30.5)(99,30.75)
\multiput(33,30.5)(8.25,.03125){8}{\line(1,0){8.25}}
%\end
\put(99,30.75){\line(2,-1){7}}
\put(66.5,35.25){\line(0,-1){4.25}}
\put(83.25,32.5){\line(-1,0){.5}}
%\emline(82.25,35.5)(83.5,31.25)
\multiput(82.25,35.5)(.03289474,-.11184211){38}{\line(0,-1){.11184211}}
%\end
%\emline(38.75,33.25)(41,30.25)
\multiput(38.75,33.25)(.03358209,-.04477612){67}{\line(0,-1){.04477612}}
%\end
\put(38.5,35.5){\circle*{3.808}}
\put(66.5,37.75){\circle*{4.031}}
\put(82.25,38.25){\circle*{3.808}}
\put(50.75,46.75){$\Delta_0$}
\put(57,26.25){$\tilde\Delta$}
\put(103,49.75){$\bf p$}
\put(104.5,31.25){$\bf q$}
\put(33.75,21.25){$\Delta_0$ is obtained by cutting off all disks from $\tilde\Delta$}
\end{picture}

 \begin{lemma}\label{not} There is no $\theta$-band in $\Delta_0$
 which starts and ends on ${\bf q}_0$ or starts and ends on $\bf p$.
 \end{lemma}
 \proof Every paths ${\bf z=z}(\Pi)$ has no $\theta$-edges. Therefore a $\theta$-band $\cal T$ starting and ending on ${\bf q_0}$
 has to start and end on $\bf q$. So the word $V$ has $\theta$-letters,
 and there is van Kampen subdiagram $\Gamma$ in $\tilde\Delta$,
 where the boundary of $\Gamma$ has form $\bf uv$, where $\bf u$ is a side of $\cal T$ and $Lab({\bf v})$  has no $\theta$-letters except for the first and the last letter; whence $Lab({\bf v})$ is a subword of a cyclic permutation of $V^{\pm 1}$. The diagram $\Gamma$ can be replaced with a minimal diagram with the same boundary label whose two $\theta$-edges have to be connected by a $\theta$-band. But this is not possible for the adapted word $V$. The $\bf p$-version of the lemma admits similar proof.
 \endproof

 \begin{lemma} \label{ring} There is a recursive function $f$ such that $k$ and $l$
do not exceed \\ $f(||U||+||V||)$ provided the path $\bf p$ has no $\theta$-edges.
\end{lemma}

\proof It follows from Lemmas \ref{NoAnnul} and \ref{not} that every maximal $\theta$-band of $\Delta_0$ is an annulus
crossing every maximal $q$-band starting on $\bf p$ exactly once. Therefore
all maximal $q$-bands starting on $\bf p$ have equal histories. The history and the one-letter base
determine side labels of a $q$-band up to superscripts. If we have two maximal $q$-bands $\cal C$ and $\cal C'$ starting with two edges $\bf e$ and $\bf e'$ of a subpath $\bf efe'$ of $\bf p$
and the length $||ef||$ is a multiple of $||U||$, then the corresponding superscripts must be equal by Remark \ref{pm1} since $\mu(U)=0$ in (\ref{0i0}),
that is $\cal C$ and $\cal C'$ have equal side labels.
So there is a set $\bf S$ of different  sides with equal labels, where $\#({\bf S})\ge k$.

Arbitrary path $\bf s$ from $\bf S$ either connects $\bf p$ and $\bf q$ or ends on  a
disks $\Pi$ of $\tilde\Delta$. In the latter case the path $\bf s$ can be extended by a subpath $\bf x$ of  the path ${\bf z}(\Pi)$. The extension $\bf s'$ connects $\bf p$ and $\bf q$. The lengths of all ${\bf z}(\Pi)$ and so the length of the extending paths
${\bf x}$ were effectively bounded  in terms of $||V||$ in Lemma \ref{dob}.  Hence there is a set of paths $\bf S'$ with equal labels, connecting $\bf p$
and $\bf q$, where $\#({\bf S'})> c'k$ and the
positive constant $c'$ depends on the number of words of bounded lengths in the generators of the group $G$.

 Arbitrary  path $\bf s'\in \bf S'$ starts with a vertex of $\bf p$, which decomposes the period $U$ of $Lab(\bf p)$ as $U\equiv U_1U_2$.
 similarly, the end of $\bf s'$ gives a factorisation $V\equiv V_1V_2$.
 If two cuts $\bf s_1, s_2 \in S'$ define the same factorizations of the words $U$ and $V$,
 we say this cuts are {\it compatible}. Since the number of factorizations of the words $U$ and $V$ are bounded, there is a set of pairwise
 compatible paths $\bf S''\subset S'$ with $\#({\bf S''})> c''k$,
 where the positive constant $c''$ is effectively bounded from below.
 However two different compatible cuts from $\bf S''$ together with parts of $\bf p$ and $\bf q$ bound a simply connected diagram with the label
 $T(U')^{k'}T^{-1}(V')^{-l'}$, where $T$ is the label of these cuts, $U'$ and $V'$ are cyclic permutations of the words $U$ and $V$, resp., and $k'<k$, $l'<l$.
 It follow that the powers $U^{k'}$ and $V^{l'}$ are conjugate in the group $G$ contrary to the choice of $k$ and $l$. Hence $c''k\le 1$,
 which effectively bounds $k$ from above. Lemma \ref{form} linearly bounds the $q$-length of the path $\bf q$ in terms of $|{\bf p}|_q$. Therefore
 the exponent $l$ is also effectively bounded.
 \endproof

 \begin{lemma} \label{spiral} There is a recursive function $f$ such that $k$ and $l$
do not exceed $f(||U||+||V||)$, provided the path ${\bf p}$ has $\theta$-edges and $q$-edges.
\end{lemma}

\proof {\bf 1.} Let $\cal C$ be a maximal $q$-band of $\Delta_0$ starting on $\bf p$. As in Lemma \ref{form}, it  ends on ${\bf q}_0$ since the word $U$ is adapted. If a $\theta$-band $\cal T$ starting from $\bf p$ crosses $\cal C$ from left to right, then it follows from Lemma \ref{NoAnnul} that it cannot cross
$\cal C$ again, but from right to left. Also there is no other $\theta$-band $\cal T'$ starting on $\bf p$ and crossing the $q$-band $\cal C$
from right to left since both $\cal T$ and $\cal T'$ cannot cross each other but both should end on $\bf q_0$.
Therefore the maximal $\theta$-bands consequently crossing $\cal C$
and starting from $\bf p$, all cross $\cal C$ from left to right (or all
cross it from right to left). It follows that these $\theta$-bands start with consecutive $\theta$-edges of $\bf p$, and so the history of $\cal C$ is a periodic word whose period is the $\theta$-projection of $U$ because $Lab({\bf p})\equiv U^k$. Moreover, the history of all maximal $q$-bands
starting with $\bf p$ are periodic words with the same period $H$, where $0<||H||<||U||$.

Furthermore, a side label of $\cal C$ is a periodic word with a period $u$, where
$|u|_{\theta}=|U|_{\theta}$. To prove this, one should show that the cell $\pi$ number
$a$ in $\cal C$ (counting from $\bf p$) is a copy of the cell $\pi'$ having number $a+|U|_{\theta}$ in $\cal C$. Indeed, if a $\theta$-band $\cal T$ ($\theta$-band $\cal T'$) starts
on $\bf p$ and crosses $q$-bands $b$ (resp. $b'$) times before it crosses $\cal C$, then $b-b' =|U|_{\theta}$. Since $\cal T$ and $\cal T'$ have the first $\theta$-edges labeled by the same letter, by Remark \ref{pm1}, we have equal superscripts when $\cal T$ and $\cal T'$ cross $\cal C$ at $\pi$ and $\pi'$, resp.,
because $\mu(U)=0$ in (\ref{0i0}).

{\bf 2.} As in Lemma \ref{ring}, a side $\bf y$ of every maximal $q$-band
admits a continuation $\bf x= yz$ in $\tilde\Delta$, where the length of $\bf z$ is bounded, and we have a set $\bf S$ of such compatible cutting paths
%with equal labels of the first $||H||$ edges
$\bf x_1, x_2,\dots, x_r$ starting with different vertices of $\bf p$, and so, all the beginnings $\bf y_1, y_2,\dots, y_r$ are the side labels of $q$-bands ${\cal C}_1,\dots, {\cal C}_r$
starting with the edges of $\bf p$ with the same base letter $q_0$.
We add the additional requirement that the prefixes of length $||H||$ of all words $Lab({\bf y_1}),\dots, Lab({\bf y_r})$ are equal (say, the histories of the corresponding $q$-band ${\cal C}_i$-s start with $H$), and still have $r > ck$, where the positive constant $c^{-1}$ is recursively bounded from above in terms of $||U||+||V||$.

Since the side label of ${\cal C}_i$-s are compatible  and $\mu (U)=0$,  we have $Lab ({\bf x})\equiv u^sv$, for every $\bf x \in S$, where $s=s({\bf x})$ and the word $v$ has bounded length. So changing the constant $c$ effectively, one may assume that the suffices $v$ are the same for
every $\bf x \in S$.
Then it follows that we have
sufficiently many different pairs of different paths $({\bf x', x''})$ from $\bf S$, where the starting vertices of the paired paths
are ''close'' to each other; more precisely, the number of pairs, where
the subpaths of $\bf p$ connecting the origins of paired paths have length $\le 3c^{-1}$ is greater than $r/2 -1$.  Let $\bf P$ be the set of such pairs.

{\bf 3.} We want to bound from above the lengths $||{\bf x'}||, ||{\bf x''}||$ for arbitrary pair $\bf (x',x'')\in P$. Thereby the number  of different labels of
the paths from such pair will be effectively bounded. However two compatible cutting paths from $\bf S$ cannot have equal labels, since as in Lemma \ref{ring}, this would lead to a contradiction with the minimality of the pair $(k,l)$

Let $E$ be a van Kampen subdiagram of $\tilde\Delta$ with boundary path $\bf x'q'(x'')^{-1}(p')^{-1}$, where $({\bf x',x'')\in P}$, ${\bf p'}$ and $\bf q'$  are subpaths of
$\bf p$ and $\bf q$, resp., and so $||{\bf p'}||\le 3c$ and $|{\bf p'}|_q\le 3c$.
This implies that $|{\bf q'}|_q\le cL$ since every maximal $q$-band starting with $\bf q$ ends either on $\bf p$ or on a disk, which also connected with $\bf p$ by $q$-bands.

Replacing the words $U$ and $V$ with cyclic permutations, we may assume that $Lab({\bf p'})\equiv U^a$ and $Lab({\bf q'})\equiv V^b$ for some $a,b>0$.

{\bf 4.} Recall that $\bf x'=y'z'$, where the length of $\bf z'$ is bounded and
$\bf y'$ is a side of a maximal $q$-band $\cal C'$ stating on $\bf p$.
Similarly, we have $\cal C''$ and $\bf x'' = y''z''$. If $E$ has a $\theta$-band connecting $\cal C'$ and $\cal C''$, we have a trapezium $\Gamma$ of maximal height formed by such $\theta$-bands and parts
of $\cal C'$ and $\cal C''$. Two  components of $E\backslash \Gamma$ (just one if $\Gamma$ is empty) have maximal $q$-subbands of bounded lengths
since maximal $\theta$-bands crossing them have at least one end on $\bf p'$ or $\bf q'$. Thus, it remains to bound the height $h$ of $\Gamma$.

By Lemma \ref{simul} (1), the top and the bottom labels $W_0$ and $W_h$
of $\Gamma$ are the first and the last permissible words of a computation
${\cal W}: W_0^{\emptyset}\to\dots\to W_h^{\emptyset}$ with periodic history having period $H$. Therefore
by Lemma \ref{perio} the height $h$ is recursively bounded in terms of
$||W_0||$ $||W_h||$, and $||H||$, provided there is no subcomputation
$W_i^{\emptyset}\to\dots\to W_j^{\emptyset}$ of $\cal W$ with history $H$ and with $W_i^{\emptyset}\equiv W_j^{\emptyset}$. Then it follows that $h$ is also effectively bounded in terms of
$||U||+||V||$, as desired.
Thus to complete the proof by contradiction, we assume now that $\cal W$
contains a subcomputation $W_i^{\emptyset}\to\dots\to W_j^{\emptyset}$  with history $H$ and with $W_i^{\emptyset}\equiv W_j^{\emptyset}$.

{\bf 5.} It follows from Lemma \ref{simul} (2) that the trapezium
$\Gamma$ contains
 a subtrapezium $\Gamma'$ corresponding to the subcomputation ${\cal W'}: W_i^{\emptyset}\to\dots\to W_j^{\emptyset}$. Since $W_i^{\emptyset}\equiv W_j^{\emptyset}$, we have $W_i\equiv W_j$, because $\Gamma$ is bounded by
 subbands of $\cal C$ and $\cal C'$, which are copies of each other,
 and so the corresponding letters of $W_i$ and $W_j$ have equal superscripts
 by Remark \ref{pm1}.

Consider now the following auxiliary surgery. Since $W_i\equiv W_j$, one can make a cut along a side of a $\theta$-band of $\Gamma$ labeled by $W_i$ and insert
a trapezium $\Gamma^{(n)}$ with history $H^n$, where $n>1$.
 The obtained
trapezium $\Gamma_n$ has the same top/bottom labels, has $H$-periodic history, but $h_n - h = n||H||$, where $h_n$ is the height of $\Gamma_n$. This surgery also replaces the diagram $E=E_0$ with a diagram $E_n$.

% This is a LaTeX picture output by TeXCAD.
% File name: [EEn.pic].
% Version of TeXCAD: 4.3
% Reference / build: 30-Jun-2012 (rev. 105)
% For new versions, check: http://texcad.sf.net/
% Options on the following lines.
%\grade{\on}
%\emlines{\off}
%\epic{\off}
%\beziermacro{\on}
%\reduce{\on}
%\snapping{\off}
%\pvinsert{% Your \input, \def, etc. here}
%\quality{8.000}
%\graddiff{0.005}
%\snapasp{1}
%\zoom{4.0000}
\unitlength 1mm % = 2.845pt
\linethickness{0.4pt}
\ifx\plotpoint\undefined\newsavebox{\plotpoint}\fi % GNUPLOT compatibility
\begin{picture}(117.75,60.25)(0,30)
%\emline(55.5,76.75)(54.75,76.5)
\multiput(55.5,76.75)(-.09375,-.03125){8}{\line(-1,0){.09375}}
%\end
%\vector[middle](21,76.5)(21,55.5)
\put(21,66){\vector(0,-1){.07}}\put(21,76.5){\line(0,-1){21}}
%\end
%\vector[middle](21,55.5)(55.5,55.5)
\put(38.25,55.5){\vector(1,0){.07}}\put(21,55.5){\line(1,0){34.5}}
%\end
%\vector[middle](21.5,76)(55,80)
\put(38.25,78){\vector(1,0){.07}}\multiput(21.5,76)(.281512605,.033613445){119}{\line(1,0){.281512605}}
%\end
%\vector[middle](55,80)(55,55.5)
\put(55,67.75){\vector(0,-1){.07}}\put(55,80){\line(0,-1){24.5}}
%\end
%\vector[middle](80,80)(115.25,85)
\put(97.625,82.5){\vector(4,1){.07}}\multiput(80,80)(.236577181,.033557047){149}{\line(1,0){.236577181}}
%\end
%\vector[middle](115.25,85)(115.25,57.5)
\put(115.25,71.25){\vector(0,-1){.07}}\put(115.25,85){\line(0,-1){27.5}}
%\end
%\vector[middle](79.5,56.75)(115.5,56.75)
\put(97.5,56.75){\vector(1,0){.07}}\put(79.5,56.75){\line(1,0){36}}
%\end
%\vector[middle](80.25,80.25)(80.25,56.5)
\put(80.25,68.375){\vector(0,-1){.07}}\put(80.25,80.25){\line(0,-1){23.75}}
%\end
\put(32.25,81.25){$U^a$}
\put(30.5,52){$V^b$}
\put(23.5,64.5){$u^sv$}
\put(57.25,64.5){$u^rv$}
\put(91.5,85.25){$U^a$}
\put(89.25,53.25){$V^b$}
\put(83,65){$u^{s+n}$}
\put(117.75,66.25){$u^{r+n}$}
\put(38.75,67.25){$E$}
\put(97.25,70.25){$E_n$}
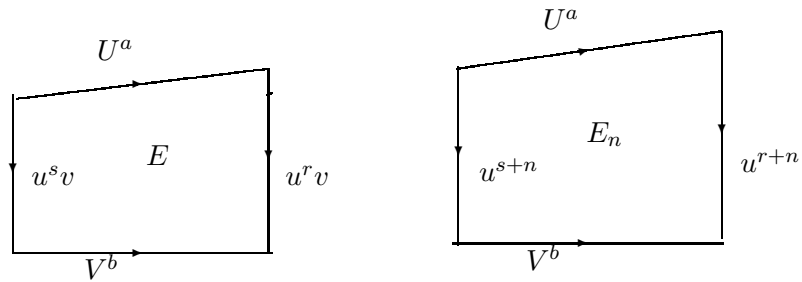
\put(43,48){Subdiagram $E$ and diagram $E_n$}
\end{picture}

Recall that by the definition of the set of cuts $\bf S$, both words $Lab ({\bf x'}) $ and $Lab ({\bf x''})$ are equal to $u^sv$ and $u^rv$ with bounded length of $v$, $Lab({\bf p'}) \equiv U^a$, $Lab({\bf q'}) \equiv V^b$. Since $a<k$ and $b< l$, we have $r\ne s$, and without loss of generality, we may assume that $r>s$. Thus, the boundary
label of $E$ gives us the equality $u^rv = U^{-a}u^sv V^b$  in $G$.
For $n =(s-r)$, the diagram $E_n$ provides us with the
equality $u^{r+n}v = U^{-a}u^{s+n}v V^b$, i.e. $u^sv = U^{-a}u^{s+n}vV^b$. Similarly, from $E_{2n}, E_{3n},\dots$,
we obtain
$$u^rv = U^{-a}u^{r+n}v V^b,\;u^{r+n}v=U^{-a}u^{r+2n}vV^b,\dots, u^{r+(l-1)n} v =U^{-a}u^{r+nl}vV^b$$
On the one hand, it follows that
\begin{equation}\label{ab}
u^rv=U^{-a}u^{r+n}v V^b=U^{-2a}u^{r+2n}vV^{2b}=\dots= U^{-la}u^{r+ln}v V^{lb}
\end{equation}
in $G$. On the other hand, cutting $\Delta$ along the path $\bf x'$,
we obtain a diagram, whose boundary label gives us $u^rv= U^{-k}u^{r}v V^{l}$ in $G$, whence
$u^rv= U^{-kb}u^{r}v V^{lb}$, which together with (\ref{ab}) gives
\begin{equation}\label{uU}
u^{ln}=U^{kb-la}
\end{equation}

{\bf 6.} To obtain the final contradiction, it remains to show that the equality
(\ref{uU}) is impossible in $G$.

The word $u$ is a label of a side of a reduced $q$-band. Therefore its label
is a  word with non-empty cyclically reduced $\nu$-projection onto
the free group generated by $\theta$-letters. If $kb-la =0$,
then by Lemmas \ref{mnogospits} and \ref{NoAnnul}, the minimal diagram for the equality $u^{ln}=1$ has neither
disks nor $(\theta, q)$-cells. So it is a diagram over a group generated
by $\theta$- and $Y$-letters. Then the homomorphism $\nu$  gives the equality $\nu(u)^{ln}=1$ in the free group, a contradiction.

If $kb-la \ne 0$, Then the van Kampen diagram $\Delta'$ corresponding to (\ref{uU}) has no disks. Indeed, otherwise by Lemma \ref{extdisk} (1), we have a disk with  $l\ge L-3$ consecutive $\tt$-spokes ${\cal C}_1,\dots,{\cal C}_r$ ending on the boundary subpath
labeled by $U^{kb-la}$, because $u$ has no $q$-letters. If there are no other disks between neighbor  ${\cal C}_i$ and ${\cal C}_{i+1}$ ($i=1,\dots, r-1$), then we have a contradiction with the property that $U$ is an adapted word. If there is a disk in a diagram $\Gamma_i$, between  some ${\cal C}_i$ and ${\cal C}_{i+1}$, then again Lemma \ref{extdisk} (1) provides
us with a disk $\pi$ in $\Gamma_i$, contrary to the definition of adapted word.
Every maximal $q$-band of $\Delta'$  has to start and end on
the boundary subpath labeled by the power of $U$, and so there is a $q$-band starting and ending on a subpath labeled by a cyclic permutation of $U^{\pm 1}$, which is impossible since the worqd $U$ is adapted.
Hence $U$ cannot contain $q$-letters, contrary to the assumption of the lemma. \endproof

\medskip

{\bf Proof of Theorem \ref{pc}(1)}. To decide if some powers $U^k$ and $V^l$ with non-zero exponents are conjugate in $G$, we may assume by Lemma \ref{enough} that the words $U$ and $V$ represent elements of infinite order. Also it can be assumed that the equality (\ref{0i0} holds and
that the words $U$ and $V$ are adapted according to Lemma \ref{simpl}.
If both $U$ and $V$ have no $q$-letters, then the exponents $k,l$ can be effectively bounded in terms of $||U||+||V||$ by Lemma \ref{noq}. Otherwise the recursive bound for $k$ and $l$ are given by Lemmas \ref{ring} and \ref{spiral}. This reduces the power conjugacy to the conjugacy of words of bounded length. Since the conjugacy problem for pairs of words of infinite orders is decidable
by Theorem \ref{c}, the power conjugacy problem is decidable in $G$.
The group $G$ has undecidable conjugacy problem and quadratic
Dehn function by Lemma \ref{conj} and \ref{quadr} if the machine ${\bf M_0}$ is choosen with nonrecursive language of accepted input words. Thus, Theorem \ref{pc} (1) is proved. $\Box$

{\bf Proof of Theorem \ref{pc}(2).} Let us start with McCool's group
$$\Pi_2 =\langle y_n, z_n (n=1,2,\dots)\mid y_nz_n = z_ny_n, y_{\phi(n)}=z_{\phi(n)}^n (n=1,2,\dots)\rangle,$$
where $\phi$ is a recursive one-to one function with a non-recursive range.
This group has decidable word problem \cite{McC}, and so it has decidable conjugacy
problem, being a free product of abelian groups. It follows from the relations
that some powers of $y_i$ and $z_i$ are conjugate if and only if they are equal,
and we can obtain such an equality if and only if $i$ belongs to the range
of the function $\phi$. Since this range is not recursive, the power conjugacy
problem is undecidable in the group $\Pi_2$.

By Theorem 3 from \cite{OS05}, the countable group $\Pi_2$ with decidable
conjugacy problem embeds in a 2-generated group $K$ with decidable conjugacy
problem. Moreover, by Lemma 8(6) from \cite{OS05}, this embedding has Frattini
property, i.e. two elements from the subgroup $\Pi_2$ are conjugate in $K$ if and only if they are conjugate in $\Pi_2$. Hence the power conjugation problem
is undecidable in $K$ too.

Finally, by Theorem 1.1 from \cite{OS04} the finitely generated group $K$
having decidable conjugacy problem Frattini embeds in a finitely presented
group with decidable conjugacy problem. Thus, the power conjugacy problem
is undecidable in $H$ too, and Theorem \ref{pc} (2) is proved. $\Box$

\end{document}